\documentclass[preprint,12pt,authoryear]{elsarticle}

\usepackage{subcaption}

\usepackage{amssymb}
\usepackage{amsmath,bm}
\usepackage{natbib}
\setcitestyle{numbers,open={[},close={]},order}

\usepackage{soul}
\usepackage{color}

\usepackage{hyperref}
\usepackage{appendix}

\begin{document}

%\tableofcontents

\sloppy

\begin{frontmatter}\if


\title{Lattice Boltzmann framework for multiphase flows by Eulerian-Eulerian Navier-Stokes equations} %% Article title


\author[label1]{Matteo Maria Piredda}
\author[label1,label2]{Pietro Asinari}
\cortext[cor1]{Corresponding author}
\ead{pietro.asinari@polito.it}

\affiliation[label1]{organization={Dipartimento di Energia, Politecnico di Torino},%Department and Organization
            addressline={Corso Duca degli Abruzzi 24}, 
            city={10129},
            postcode={Torino(TO)}, 
            country={Italy}}

\affiliation[label2]{organization={Istituto Nazionale di Ricerca Metrologica},%Department and Organization
            addressline={Strada delle Cacce 91}, 
            city={10135},
            postcode={Torino(TO)}, 
            country={Italy}}

%% Abstract
\begin{abstract}
Although Lattice Boltzmann Method (LBM) is relatively straightforward, it demands a well-crafted framework to handle the complex partial differential equations involved in multiphase flow simulations. For the first time to our knowledge, this work proposes a novel LBM framework to solve Eulerian-Eulerian multiphase flow equations, without any finite difference correction, including very large density ratios and also a realistic relation for the drag coefficient. The proposed methodology and all reported LBM formulas can be applied to any dimension. This opens a promising venue for simulating multiphase flows on large High Performance Computing (HPC) facilities and on novel parallel hardware. This LBM framework consists of six coupled LBM schemes - running on the same lattice - ensuring an efficient implementation in large codes with minimum effort. The preliminary numeral results agree in an excellent way with the a reference numerical solution obtained by a traditional finite difference solver.
\end{abstract}

%\includegraphics[scale=0.1]{polito_logo_2021_blu.jpeg}

%\item Mathematical modeling
%\item Computational fluid dynamics
%\item Lattice Boltzmann Method
%\item Multiphase flows

%% Keywords
\begin{keyword}
Computational Fluid Dynamics \sep multiphase flows \sep Eulerian-Eulerian Navier-Stokes equations \sep Lattice Boltzmann Method
\end{keyword}

\end{frontmatter}

%\hl{VERSION v4 !} 
            


\section{Introduction}

Computational Fluid Dynamics (CFD) of multiphase flows is crucial in the energy sector, particularly in oil and gas, as it enables the detailed simulation and optimization of complex fluid interactions, such as those in bubble column reactors. These reactors are used extensively in refining processes, gas treatment, and chemical synthesis, where accurate modeling of gas-liquid interactions can lead to improved efficiency, reduced energy consumption, and enhanced safety. Using CFD, industry can better understand the flow dynamics, optimize reactor designs, and ultimately reduce costs and emissions, making it a vital tool for the advancement of sustainable energy solutions.

A multiphase flow system consists of the simultaneous flow of materials with different phases or states of matter in regions bounded by moving interfaces. A possible approach to modeling multiphase flows is the mixture approach, which essentially treats the multiphase system as a single mixture with properties that are weighted averages of the different phases \cite{vtt1996}. Instead of solving separate momentum equations for each phase, a single momentum equation is solved for the mixture. For this reason, it is a simple and computationally convenient approach, but with a lower detail level because a proper model must be provided for the velocity difference between the two phases. It needs small relative velocities, so it is more suitable for closely coupled phases moving together (slurry, sediment).

This work, on the other hand, focuses on developing an Eulerian-Eulerian lattice Boltzmann method (LBM) to simulate multiphase flows. This implies the application of a volume averaging procedure to both the continuity and momentum equation for every phase and it requires to define a volume fraction as follows:
\begin{equation}
    \alpha_\varphi = \frac{V_\varphi}{V}, \hspace{2cm} \sum_{\varphi=1}^{P}\alpha_\varphi=1,
\end{equation}
where V is the volume under consideration, the subscript $\varphi$ stays for a generic phase and P for the total number of phases. In the following, for the sake of simplicity and without loss of generality, let us focus on two phases, namely $P=2$, where two phases are identified by $l$ for the liquid phase and $g$ for the gas phase, respectively. In case $P=2$, if $\varphi$ stays for a generic phase, then $\overline{\varphi}$ stays for the other one. Writing down the governing equation for both phases constitutes the Eulerian-Eulerian approach (or two-fluid approach), and is a framework used to model multiphase flows where each phase is treated as a continuum, and separate sets of equations are solved for each phase in the same Eulerian (fixed) reference frame. They are used in fluid dynamics to describe the motion of a fluid while accounting for the effects of spatial variations over a finite region \cite{book2007}. 
Let's start by recalling the original Eulerian-Eulerian Navier-Stokes equations for multiphase flows, in the isothermal limiting case, used in traditional CFD for both phases \cite{Maniscalco2021}:
\begin{equation}\label{cont-l}
    \frac{\partial (\rho_l \alpha_l)}{\partial t} + \nabla \cdot (\alpha_l \rho_l \textbf{u}_l) = 0,
\end{equation}
\begin{equation}
    \frac{\partial (\rho_l \alpha_l \textbf{u}_l)}{\partial t} + \nabla \cdot (\alpha_l \rho_l \textbf{u}_l \textbf{u}_l) = - \alpha_l \nabla p + \rho_l \alpha_l \textbf{g} + \nabla \cdot (\alpha_l \rho_l \bm{\sigma}_l) + \mathbf{F}_{lg} + \mathbf{F}_{l},
\end{equation}
\begin{equation}\label{cont-g}
    \frac{\partial (\rho_g \alpha_g)}{\partial t} + \nabla \cdot (\alpha_g \rho_g \textbf{u}_g) = 0,
\end{equation}
\begin{equation}
    \frac{\partial (\rho_g \alpha_g \mathbf{u}_g)}{\partial t} + \nabla \cdot (\alpha_g \rho_g \mathbf{u}_g \mathbf{u}_g) = - \alpha_g \nabla p + \rho_g \alpha_g \mathbf{g} + \nabla \cdot (\alpha_g \rho_g \bm{\sigma}_g) + \mathbf{F}_{gl} + \mathbf{F}_{g}.
\end{equation}
For each phase (the subscript $\varphi$ can be either $l$ or $g$), $\alpha_\varphi$ is the previously mentioned volume fraction, $\rho_\varphi$ is the density, $\mathbf{u}_\varphi$ the velocity vector, $\bm{\sigma}_\varphi$ the viscous stress tensor per unit volume, $p$ the common pressure, $\mathbf{F}_{gl}=-\mathbf{F}_{lg}$ is the interphase momentum exchange term, which is a cumulative force resulting from the summation of the single interfacial forces acting between the two phases (drag force, lift force, wall lubrication, virtual mass, etc) per unity volume and $\mathbf{F}_{\varphi}$ is the force per unit volume of the specific phase $\varphi$. Concerning the interphase momentum exchange, let us suppose to limit the following discussion to the drag contribution only, which can be formulated as
\begin{equation}\label{interphasemomentum}
    \mathbf{F}_{gl} = \rho_g\,K_I\,\|\mathbf{u}_l-\mathbf{u}_g\|\,
    \left(\mathbf{u}_l-\mathbf{u}_g\right)=-\mathbf{F}_{lg},
\end{equation}
where $K_I$ is the effective drag coefficient and is evaluated from empirical or semi-empirical correlations. The viscous stress tensor per unit volume is defined as follows:
\begin{equation}\label{stresstensor}
    \bm{\sigma}_\varphi = \nu_\varphi\left(\nabla \mathbf{u}_\varphi + \nabla \mathbf{u}_\varphi^T\right) 
    +\left(\xi_\varphi- \frac{2}{3}\,\nu_\varphi\right)\,(\nabla \cdot \mathbf{u}_\varphi)\,\mathbf{I},
\end{equation}
where $\nu_\varphi$ is the effective kinematic viscosity for the generic phase, defined as the ratio between dynamic viscosity and density, i.e. $\nu_\varphi = \mu_\varphi/\rho_\varphi$, and $\xi_\varphi$ is the effective kinematic bulk viscosity. In general, the effective kinematic viscosity should also include turbulent effects. It is important to realize that, in multiphase flows, even under the incompressible limit which we will discuss in the following, $\nabla \cdot \mathbf{u}_\varphi$ for the individual phase can be different from zero in general. 

The previous equations define a proper system of equations for the following variables: $\alpha_g$, $p$, $\mathbf{u}_g$ and $\mathbf{u}_l$, as far as the equation of state for the dispersed phase, i.e. $\rho_g=\rho_g(p)$, and the equation of state for the liquid phase, i.e. $\rho_l=\rho_l(p)$ in the isothermal limiting case, are provided.

In case of the incompressible limit (i.e. low Mach number limit), one can assume the following equations of state: $\rho_g=\rho_g(p)=\rho_g^0$ and $\rho_l=\rho_l(p)=\rho_l^0$, where $\rho_g^0$ and $\rho_l^0$ are proper constants consisting in the average density of each phase. In this case, the previous equations reduce to the following:
\begin{equation}\label{continuity_l}
    \frac{\partial \alpha_l}{\partial t} + \nabla \cdot (\alpha_l \textbf{u}_l) = 0,
\end{equation}
\begin{equation}
    \frac{\partial (\alpha_l \textbf{u}_l)}{\partial t} + \nabla \cdot (\alpha_l \textbf{u}_l \textbf{u}_l) = - \frac{\alpha_l}{\rho_l^0} \nabla p + \alpha_l \textbf{g} + \nabla \cdot (\alpha_l \bm{\sigma}_l) + \frac{1}{\rho_l^0} \left(\mathbf{F}_{lg}+\mathbf{F}_{l}\right),
\end{equation}
\begin{equation}\label{voidfraction}
    \frac{\partial \alpha_g}{\partial t} + \nabla \cdot (\alpha_g \textbf{u}_g) = 0,
\end{equation}
\begin{equation}
    \frac{\partial (\alpha_g \mathbf{u}_g)}{\partial t} + \nabla \cdot (\alpha_g \mathbf{u}_g \mathbf{u}_g) = - \frac{\alpha_g}{\rho_g^0} \nabla p + \alpha_g \mathbf{g} + \nabla \cdot (\alpha_g \bm{\sigma}_g) + \frac{1}{\rho_g^0} \left(\mathbf{F}_{gl}+\mathbf{F}_{g}\right).
\end{equation}
From the numerical point of view, the previous formulation is not convenient and it is better to derive an equivalent system of equations. Summing up continuity equations for both phases in the incompressible limit yields:
\begin{equation}\label{continuitysum}
     \nabla \cdot (\alpha_g \textbf{u}_g + \alpha_l \textbf{u}_l) = 0.
\end{equation}
Summing up momentum equations for both phases in the incompressible limit yields:
\begin{multline}
    \frac{\partial (\alpha_g \mathbf{u}_g+\alpha_l \mathbf{u}_l)}{\partial t} + \nabla \cdot (\alpha_g \mathbf{u}_g \mathbf{u}_g+\alpha_l \mathbf{u}_l \mathbf{u}_l) = - \left(\frac{\alpha_g}{\rho_g^0}+\frac{\alpha_l}{\rho_l^0}\right) \nabla p + \mathbf{g} +\dots\\ 
    \dots+\nabla \cdot (\alpha_g \bm{\sigma}_g) + 
    \nabla \cdot (\alpha_l \bm{\sigma}_l)
    + \frac{1}{\rho_g^0} \left(\mathbf{F}_{gl}+\mathbf{F}_{g}\right)
    + \frac{1}{\rho_l^0} \left(\mathbf{F}_{lg}+\mathbf{F}_{l}\right).
\end{multline}
Applying the divergence operator to both sides of the equation and swapping the order of the derivatives in the first term yield:
\begin{multline}
    \frac{\partial}{\partial t} \left[\nabla\cdot(\alpha_g \mathbf{u}_g+\alpha_l \mathbf{u}_l)\right] + \nabla \cdot \nabla \cdot (\alpha_g \mathbf{u}_g \mathbf{u}_g+\alpha_l \mathbf{u}_l \mathbf{u}_l) = - \nabla \cdot\left[\left(\frac{\alpha_g}{\rho_g^0}+\frac{\alpha_l}{\rho_l^0}\right) \nabla p\right] \dots\\ 
    +\nabla \cdot\left[ \nabla \cdot (\alpha_g \bm{\sigma}_g) + 
     \nabla \cdot (\alpha_l \bm{\sigma}_l)
    + \frac{1}{\rho_g^0} \left(\mathbf{F}_{gl}+\mathbf{F}_{g}\right)
    + \frac{1}{\rho_l^0} \left(\mathbf{F}_{lg}+\mathbf{F}_{l}\right)\right].
\end{multline}
Applying Eq. (\ref{continuitysum}), it is possible to derive the Poisson equation, which is typically solved for computing the common pressure field $p$, namely
\begin{multline}\label{poisson}
    \nabla \cdot \nabla \cdot (\alpha_g \mathbf{u}_g \mathbf{u}_g+\alpha_l \mathbf{u}_l \mathbf{u}_l) = - \nabla \cdot\left[\left(\frac{\alpha_g}{\rho_g^0}+\frac{\alpha_l}{\rho_l^0}\right) \nabla p\right] \dots\\ 
    +\nabla \cdot\left[ \nabla \cdot (\alpha_g \bm{\sigma}_g) + \nabla \cdot (\alpha_l \bm{\sigma}_l)
    + \frac{1}{\rho_g^0} \left(\mathbf{F}_{gl}+\mathbf{F}_{g}\right)
    + \frac{1}{\rho_l^0} \left(\mathbf{F}_{lg}+\mathbf{F}_{l}\right)\right].
\end{multline}
From the numerical point of view, it is better to reformulate the momentum equations in the incompressible limit. This is done by expanding the derivatives in the left-hand side of the equations and applying Eq. (\ref{continuity_l}) and Eq. (\ref{voidfraction}) to the momentum equation of the liquid and dispersed phase, respectively:
\begin{equation}\label{dispersemom}
    \frac{\partial \mathbf{u}_g}{\partial t} + \mathbf{u}_g \cdot\nabla \mathbf{u}_g = - \frac{1}{\rho_g^0} \nabla p + \mathbf{g} 
    + \frac{1}{\alpha_g} \nabla \cdot (\alpha_g \bm{\sigma}_g)
    + \frac{1}{\alpha_g\rho_g^0} \left(\mathbf{F}_{gl}
    +\mathbf{F}_{g}\right),
\end{equation}
\begin{equation}\label{liquidmom}
    \frac{\partial \mathbf{u}_l}{\partial t} + \mathbf{u}_l\cdot\nabla \mathbf{u}_l = - \frac{1}{\rho_l^0} \nabla p + \mathbf{g} + \frac{1}{\alpha_l}\nabla \cdot (\alpha_l \bm{\sigma}_l)
    + \frac{1}{\alpha_l\rho_l^0}\left(\mathbf{F}_{lg}+\mathbf{F}_{l}\right).
\end{equation}
The singularity of the forces depend on the terms $\nabla\ln{(\alpha_\varphi)}$ (where $\varphi$ can be either $g$ or $l$). There are well-established techniques in standard computational fluid dynamics for multiphase flows to handle exactly this kind of terms  \cite{Oliveira2003}. For this term to be well-behaved as $\alpha_\varphi \to 0$ it is necessary for the gradient $\nabla \alpha_\varphi$ to approach zero faster than $\alpha_\varphi$. Numerically, in standard CFD methods, it is straightforward to discretize this term in a way that prevents division by zero: this can be achieved by representing $\alpha_\varphi$ in the denominator using a proper volumetric average and/or by applying a proper slope limiter \cite{laney1998}.

The system of equations defined by Eq. (\ref{voidfraction}), Eq. (\ref{poisson}), Eq. (\ref{dispersemom}) and Eq. (\ref{liquidmom}) in terms of quantities $\alpha_g$, $p$, $\mathbf{u}_g$ and $\mathbf{u}_l$ represents a promising starting point for numerics in most of the existing software for solving multiphase flows by the Eulerian-Eulerian approach (e.g. OpenFOAM). 

In spite of the existence of feasible numerical way, still there are remaining complexities of the Eulerian-Eulerian approach which must be faced. The solution of the multi-fluid set of equations presents many challenges:
\begin{itemize}
    \item possible singularities of the phase momentum equations;
    \item coupling between the phases, which could lead to instabilities of the numerical procedure;
    \item phase volume fraction needs to be bounded between 0 and 1;
    \item sharp profiles of phase volume fractions;
    \item extension of iterative solution procedures to the co-located grid arrangement for avoiding checkerboard instability patterns.
\end{itemize}
For example, the open-source code OpenFOAM implements a numerical iterative solution procedure proposed by Passalacqua et al. \cite{Passalacqua2011}, originally developed  for fluid-particle flows. It consists in a finite volume method (FVM) discretizing the Navier-Stokes equations Eq. (\ref{voidfraction}), Eq. (\ref{poisson}), Eq. (\ref{dispersemom}) and Eq. (\ref{liquidmom}). It uses face fluxes and velocity fluxes to overcome the problems described above along with deriving the pressure equation and the dispersed phase continuity equation. Despite its effectiveness in FVM, this procedure is clearly unfeasible for LBM, e.g. because LBM cannot easily solve Eq. (\ref{poisson}) and cannot straightforwardly impose the same pressure gradient $\nabla p$ to both phases.

\section{Eulerian-Eulerian lattice Boltzmann method (LBM) to simulate multiphase flows}\label{LBMdesign}

The Lattice Boltzmann Method (LBM) is highly promising for simulating multiphase flows on High Performance Computing (HPC) systems due to its unique computational structure, which is naturally parallel and localized. Unlike conventional CFD methods that rely on solving the Navier-Stokes equations through techniques like finite difference (FD), finite volume method (FVM), or finite element method (FEM), LBM simulates fluid flows by modeling the fluid as discrete particles moving and colliding on a lattice \cite{LALLEMAND2021109713}. This approach enables LBM to perform computations locally at each grid point, making it well-suited for parallel execution across HPC platforms.

For LBM to fully exploit HPC capabilities, it is essential to maintain the standard, unmodified LBM formulation, avoiding FD corrections or other modifications that introduce non-local dependencies. These corrections, sometimes used to address stability or accuracy, create dependencies that disrupt the purely local and independent computations that are a key advantage of LBM. Such non-local adjustments increase computational costs, reduce parallel efficiency, and complicate memory access patterns, which diminishes the performance of LBM on large HPC systems.

Therefore, adhering closely to the standard LBM formulation preserves its localized computation advantage, reducing communication overhead and allowing efficient scaling across numerous computational nodes. This efficiency enables the simulation of high-resolution, large-domain multiphase flows crucial for energy sector applications, like bubble column reactors in oil and gas, where detailed fluid dynamics insights are essential for optimizing processes.

This paper aims to derive a LBM framework to simulate Eulerian-Eulerian equations for multiphase flows. The overall scheme consists of six LBM schemes for the two phases, namely:
\begin{itemize}
  \item two LBM schemes for artificially compressible continuity equations and momentum equations, in terms of particle distribution functions $f_g$ and $f_l$,
  \item two LBM schemes for phase volume fractions, in terms of $f_{\alpha g}$ and $f_{\alpha l}$;
  \item two LBM schemes for phase continuity sources, in terms of $f_{\beta g}$ and $f_{\beta l}$.
\end{itemize}
The proposed LBM framework is derived and analyzed in the rest of this section in a general, way which does not pose any constraint on the physical dimensionality, while the following section reports a preliminary numerical validation.

\subsection{The key point: Artificially compressible continuity equation and momentum equation for each phase}

The LBM uses the artificial compressibility concept \cite{asinariOhwada} and hence it cannot solve Eq. (\ref{poisson}) directly. Hence a pseudo-compressible system of equations for each phase is needed as a starting point for deriving the Eulerian-Eulerian LBM framework to simulate multiphase flows. Applying the artificial compressibility concept requires a system of equations made of (a) a momentum equation and (b) an artificially compressible continuity equation, which is obviously missing in the incompressible formulation, for each phase. Let us suppose to start with the dispersed phase and let's start by modifying Eq. (\ref{continuitysum}) in the following way
\begin{equation}\label{acm}
     \frac{\partial \epsilon_g}{\partial t}+     
     \nabla \cdot (\alpha_g \textbf{u}_g + \alpha_l \textbf{u}_l) = 0,
\end{equation}
where $\epsilon_g$ is a function computed by flow quantities. If the artificial compressibility is used as a shortcut to find out only the steady state solution, then it is enough that $\partial \epsilon_g/\partial t$ is small enough with regards to the other terms. However, in this work, the idea is to recover the right incompressible dynamics as well and hence $\partial \epsilon_g/\partial t$ must be asymptotically small during the entire dynamics, according to the order of convergence of the numerical scheme. For example, this term could be $\partial \epsilon_g/\partial t = O(h^2)$ in case of a second-order method, where $h$ is the mesh spacing. The last condition can be recovered by a proper scaling, namely a proper choice of the simulation parameters during mesh refinement. We will enforce this condition in the following, after discussing the asymptotic analysis of the proposed methodology. Let's manipulate the previous equation as follows
\begin{equation}\label{lbm-cont-g}
    \frac{\partial \epsilon_g}{\partial t}+
    \nabla \cdot \textbf{u}_g = 
    \nabla \cdot [\alpha_l (\textbf{u}_g - \textbf{u}_l)]
    \equiv S_{g}.
\end{equation}
The previous equation is an artificially compressible continuity equation designed for the dispersed phase velocity field $\textbf{u}_g$. We have to derive a similar equation for the liquid phase velocity field $\textbf{u}_l$. Hence, this time, let's modify Eq. (\ref{continuitysum}) as follows
\begin{equation}\label{lbm-cont-l}
    \frac{\partial \epsilon_l}{\partial t}+
    \nabla \cdot \textbf{u}_l = 
    \nabla \cdot [\alpha_g (\textbf{u}_l - \textbf{u}_g)]
    \equiv S_{l}.
\end{equation}
Different strategies are possible in choosing the pair of functions $\epsilon_g$ and $\epsilon_l$. 
\begin{itemize}
    \item Because both previous equations are derived from the same Eq. (\ref{continuitysum}), then one could assume $\epsilon_g=\epsilon_l=\epsilon(p)$. The simplest choice would be $\epsilon(p)=p$ or, by looking at the Poisson equation for the mixture velocity given by Eq. (\ref{poisson}), one could use instead $\epsilon(p)=(\alpha_g/\rho_g^0+\alpha_l/\rho_l^0)\,p$, introducing a further dependence on the volume fraction in the pseudo-compressibility. 
    \item The problem with the previous approach is that it forces both phases, not only to have the same asymptotic target given by Eq. (\ref{continuitysum}), but also to have the same approaching dynamics. Because the two phases are subject to different forces, this could lead to over constraining. Hence one could assume instead two independent functions $\epsilon_g=\epsilon_g(p)$ and $\epsilon_l=\epsilon_l(p)$. In this work, we choose this second approach. The explicit expressions of these functions will be provided in the following.
\end{itemize}

Concerning the momentum equations, let's manipulate Eq. (\ref{dispersemom}) and Eq. (\ref{liquidmom}) as follows
\begin{equation}\label{dispersemom2}
    \frac{\partial \mathbf{u}_g}{\partial t} + 
    \nabla\cdot (\mathbf{u}_g\mathbf{u}_g)
    -\mathbf{u}_g\nabla\cdot \mathbf{u}_g = - \frac{1}{\rho_g^0} \nabla p + \mathbf{g} + \frac{1}{\alpha_g} \nabla \cdot (\alpha_g \bm{\sigma}_g)
    + \frac{1}{\alpha_g\rho_g^0} \left(\mathbf{F}_{gl}+\mathbf{F}_{g}\right),
\end{equation}
\begin{equation}\label{liquidmom2}
    \frac{\partial \mathbf{u}_l}{\partial t} + 
    \nabla\cdot (\mathbf{u}_l\mathbf{u}_l)
    -\mathbf{u}_l\nabla\cdot \mathbf{u}_l = - \frac{1}{\rho_l^0} \nabla p + \mathbf{g} + \frac{1}{\alpha_l}\nabla \cdot (\alpha_l \bm{\sigma}_l)
     + \frac{1}{\alpha_l\rho_l^0}\left(\mathbf{F}_{lg}+\mathbf{F}_{l}\right).
\end{equation}
Using Eq. (\ref{lbm-cont-g}) and Eq. (\ref{lbm-cont-l}), coherently with the artificial compressibility approach, $\nabla\cdot \mathbf{u}_g= S_{g} + O(h^2)$ and $\nabla\cdot \mathbf{u}_l= S_{l} + O(h^2)$. Let us assume that $\nabla\cdot \mathbf{u}_g\approx S_{g}$ and $\nabla\cdot \mathbf{u}_l\approx S_{l}$, which can be used to simplify the previous equations. For the dispersed phase, this yields
\begin{equation}\label{lbm-mom-g}
    \frac{\partial \mathbf{u}_g}{\partial t} + \nabla \cdot (\mathbf{u}_g \mathbf{u}_g) = - \nabla\left(\frac{p}{\rho_g^0}\right) + \nabla \cdot \bm{\sigma}_g + \mathbf{G}_{g},
\end{equation}
where
\begin{equation}\label{Gg}
    \mathbf{G}_{g} \equiv S_{g}\mathbf{u}_g + 
    \frac{1}{\alpha_g}\bm{\sigma}_g \cdot \nabla\alpha_g + \mathbf{g} + \frac{1}{\alpha_g\rho_g^0}(\mathbf{F}_{gl}+\mathbf{F}_{g}).
\end{equation}
Similarly, it holds:
\begin{equation}\label{lbm-mom-l}
    \frac{\partial \mathbf{u}_l}{\partial t} + \nabla \cdot (\mathbf{u}_l \mathbf{u}_l) = - \nabla\left(\frac{p}{\rho_l^0}\right) + \nabla \cdot \bm{\sigma}_l + \mathbf{G}_{l},
\end{equation}
where
\begin{equation}\label{Gl}
    \mathbf{G}_{l} \equiv S_{l}\mathbf{u}_l + \frac{1}{\alpha_l}\bm{\sigma}_l \cdot \nabla\alpha_l+
    \mathbf{g} +
    \frac{1}{\alpha_l\rho_l^0}(\mathbf{F}_{lg}+\mathbf{F}_{l}).
\end{equation}

For a given dispersed phase volume fraction $\alpha_g$, Eq. (\ref{lbm-cont-g}) and Eq. (\ref{lbm-mom-g}) may be target equations for a LBM scheme for the dispersed phase. Similarly Eq. (\ref{lbm-cont-l}) and Eq. (\ref{lbm-mom-l}) may be target equations for a LBM scheme for the liquid phase. However, it is important to highlight the following peculiarity: Eq. (\ref{lbm-cont-g}) and Eq. (\ref{lbm-cont-l}) are driven by the same pressure time derivative because $\partial_t\,\epsilon_g=(\partial\epsilon_g/\partial p)\,\partial_t\,p$ and $\partial_t\,\epsilon_l=(\partial\epsilon_l/\partial p)\,\partial_t\,p$ (time pressure coupling) and, moreover, Eq. (\ref{lbm-mom-g}) and Eq. (\ref{lbm-mom-l}) are coupled by the same pressure gradient $\nabla p$ (space pressure coupling). It is important to remind at this point that the system of equations for both phases can be closed only by providing the additional equation for $\alpha_g$ (which is discussed in the following sections).

Before entering into the details of the LBM schemes, it is worth to discuss first both the time and the space pressure coupling, mentioned above, in the LBM context. Let us recall first some preliminaries about LBM. Let us define $\epsilon$ as the zero-th order moment of the distribution function, namely $\epsilon=\sum_q f(q)$, where $q$ identifies the probability distribution function $f(q)$ corresponding to the lattice velocity $\mathbf{v}_q\in \mathbb{L}$ in the velocity lattice/set $\mathbb{L}$, designed for the $D$-dimensional physical space and for the $Q$-dimensional velocity space. In the LBM framework, as it will be clearer in the following, it is possible to impose an equation of state where the pressure is $p = c_s^2\,\phi\,\rho\,\epsilon$, where $c_s\sqrt{\phi}$ is the (artificial) sound speed and $\phi$ is a tunable function (close to $1$ for stability reasons). Let us now assume that Eq. (\ref{lbm-cont-g}) and Eq. (\ref{lbm-mom-g}) are solved by a LBM scheme in terms of $f_g$ and, for the sake of simplicity, let us assume $\phi_g=1$. This allows one to identify the pressure $p$ as 
\begin{equation}\label{lbm-pressure-g}
    p \equiv c_s^2\rho_g^0 \epsilon_g,
\end{equation}
where $\epsilon_g=\sum_q f_g(q)$ or equivalently
\begin{equation}\label{sigmag}
    \epsilon_g(p) = \frac{p}{c_s^2\rho_g^0},
\end{equation}
which identifies the first function appearing in Eq. (\ref{lbm-cont-g}).

Let us now assume that Eq. (\ref{lbm-cont-l}) and Eq. (\ref{lbm-mom-l}) are solved by another LBM scheme in terms of $f_l$ and that this time $\phi_l\neq1$, which requires a specific equilibrium distribution function. Because Eq. (\ref{lbm-mom-g}) and Eq. (\ref{lbm-mom-l}) are coupled, one has to ensure that the same pressure gradient (space pressure coupling) will drive also the evolution of this second LBM scheme for the liquid phase as well. The target Eq. (\ref{lbm-mom-l}) involves the following term, which is ruled by the generalized equation of state, namely
\begin{equation}
    \nabla\left(\frac{p}{\rho_l^0}\right) = \nabla\left(\phi_l\,c_s^2\epsilon_l\right),
\end{equation}
where $\epsilon_l=\sum_l f_l(q)$. Substituting Eq. (\ref{lbm-pressure-g}) into the previous one yields
\begin{equation}
     \nabla\left(\phi_l\,\epsilon_l\right) = 
     \frac{1}{R}\,\nabla\epsilon_g,
\end{equation}
where $R=\rho_l^0/\rho_g^0> 1$. Solving the previous equations yields
\begin{equation}
     \phi_l\,\epsilon_l = 
     \frac{1}{R}\,\epsilon_g+k,
\end{equation}
where $k$ is an arbitrary constant. In order to choose the proper constant $k$, let us remember that $\epsilon_\varphi=\sum_q f_\varphi(q)$ is close to $1$ in the incompressible limit and that one wants $\phi_l$ being close to $1$ for stability reasons. Assuming $k=-1/R+1$ yields
\begin{equation}\label{lbm-phi-l}
     \phi_l = \phi \equiv \frac{1}{\epsilon_l}\,
     \left[1+\frac{1}{R}\,\left(\epsilon_g-1\right)\right],
\end{equation}
or equivalently
\begin{equation}\label{sigmal}
     \epsilon_l(p) = \frac{1}{\phi}\,
     \left[1+\frac{1}{R}\,\left(\frac{p}{c_s^2\rho_g^0}-1\right)\right],
\end{equation}
which identifies the second function appearing in Eq. (\ref{lbm-cont-l}). Clearly $\epsilon_g(p)$ given by Eq. (\ref{sigmag}) and $\epsilon_l(p)$ given by Eq. (\ref{sigmal}) are driven by the same pressure dynamics (time pressure coupling), but they are different functions. 

Let us recap the work done so far. Eq. (\ref{lbm-cont-g}) with source $S_{g}$ and Eq. (\ref{lbm-mom-g}) with force $\mathbf{G}_{g}$ can be solved by a LBM scheme in terms of $f_g$ with a standard equation of state given by Eq. (\ref{lbm-pressure-g}). On the other hand, Eq. (\ref{lbm-cont-l}) with source $S_{l}$ and Eq. (\ref{lbm-mom-l}) with force $\mathbf{G}_{l}$ can be solved by a LBM scheme in terms of $f_l$ with a generalized equation of state, i.e. with $\phi_l$ given by Eq. (\ref{lbm-phi-l}). In the next section, we will discuss the complete LBM schemes proposed for solving these equations.

\subsubsection{LBM schemes for solving the key equations}

Let us consider the following two LBM schemes for solving the artificially compressible continuity equations and the momentum equations for the two phases, formulated in terms of the corresponding particle distribution functions $f_g$ and $f_l$, namely
\begin{equation}\label{lbmg}
f_g(\hat{\mathbf{x}}+\mathbf{v}_q,
\hat{t}+1,q)=
f_g(q)+\hat{\Omega}_g(q)+
\hat{\Omega}_g^F(q),
\end{equation}
\begin{equation}\label{lbml}
f_l(\hat{\mathbf{x}}+\mathbf{v}_q,
\hat{t}+1,q)=
f_l(q)+\hat{\Omega}_l(q)+
\hat{\Omega}_l^F(q),
\end{equation}
where all functions in the right hand side are computed in $(\hat{\mathbf{x}},\hat{t})$ locally, $\hat{\mathbf{x}}$ is the space coordinate divided by the distance $\lambda$ between two neighboring lattices nodes (mean free path), $\hat{t}$ is the physical time divided by the time $\tau$ between two consecutive lattice collisions (mean collision time), $\mathbf{v}_q$ is the generic particle velocity divided by the average particle velocity $c=\lambda/\tau$, in the velocity lattice/set $\mathbb{L}$ designed for the $D$-dimensional physical space and for the $Q$-dimensional velocity space. These normalizations are consistent with the usual Boltzmann scaling adopted in the LBM numerical codes, which will be denoted by the hat notation $\hat{\cdot}$ from here on, in order to simplify the discussion about how to tune the simulation parameters for solving the target equations.  In the following, we discuss the terms $\hat{\Omega}_\varphi$ and $\hat{\Omega}_\varphi^F$. 

The term $\hat{\Omega}_\varphi$ is the collisional operator, defined as 
\begin{equation}\label{colloperator}
\hat{\Omega}_\varphi(q)=
\omega_\varphi
\left[f^{eq}_\varphi(q) - f_\varphi(q)\right],
\end{equation}
where $\omega_\varphi$ is the relaxation frequency divided by the frequency $1/\tau$ for the phase $\varphi$, which drives the relaxation of the distribution function towards the equilibrium $f^{eq}_\varphi$. The equilibrium $f^{eq}_\varphi$ is defined by means of the so-called incompressible equilibrium $f_{I}^e(\phi,\hat{\epsilon},\hat{\mathbf{u}})$. This equilibrium is called incompressible because all terms of the equilibrium moments which depend on the velocity $\hat{\mathbf{u}}$ do not depend on the zero-th order moment $\hat{\epsilon}$, e.g. $\sum_q \mathbf{v}_q f_I^e(\phi,\hat{\epsilon},\hat{\mathbf{u}},q)=\hat{\mathbf{u}}$. The equilibrium distribution is different for the two phases, even though it is formulated by means of the same functional form in the incompressible limit \cite{KruegerLBM}, namely
%
%\begin{equation}\label{equilibrium}
%f^{eq}(\rho,\mathbf{u},q) = \rho\,w_q\,
%\left[1+\frac{\mathbf{v}_q\cdot\mathbf{u}}{c_s^2}
%+\frac{(\mathbf{v}_q\mathbf{v}_q-c_s^2\,\mathbf{I}):
%\mathbf{u}\mathbf{u}}{2 c_s^4}\right]
%\end{equation}
%
%\begin{equation}\label{equilibriumIrho}
%f_{I}^{eq}(\phi,\hat{\rho},\hat{\mathbf{u}},q) = w_q\,
%\left[\hat{\rho}\,\eta_q(\phi)+
%\frac{\mathbf{v}_q\cdot\hat{\mathbf{u}}}{c_s^2}
%+\frac{(\mathbf{v}_q\mathbf{v}_q-%c_s^2\,\mathbf{I}):
%\hat{\mathbf{u}}\hat{\mathbf{u}}}{2 c_s^4}\right]
%\end{equation}
%
\begin{equation}\label{equilibriumI}
f_{I}^{eq}(\phi,\hat{\epsilon},\hat{\mathbf{u}},q) = w_q\,
\left[\hat{\epsilon}\,\eta_q(\phi)+
\frac{\mathbf{v}_q\cdot\hat{\mathbf{u}}}{c_s^2}
+\frac{(\mathbf{v}_q\mathbf{v}_q-c_s^2\,\mathbf{I}):
\hat{\mathbf{u}}\hat{\mathbf{u}}}{2 c_s^4}\right],
\end{equation}
where
\begin{equation}\label{eta}
\eta_q(\phi)=\delta_q\left[1+\frac{(1-w_0)}{w_0}\,(1-\phi)\right]
+(1-\delta_q)\,\phi,
\end{equation}
and
\begin{equation}\label{delta}
\delta_q = 
    \begin{cases}
            1, &         \text{if } q=0,\\
            0, &         \text{if } q\neq 0.
    \end{cases}
\end{equation}
The two equilibrium distributions differ from each other because $f^{eq}_g \equiv f_{I}^{eq}(1,\hat{\epsilon}_g,\hat{\mathbf{u}}_g)$ for the dispersed phase, while $f^{eq}_l \equiv f_{I}^{eq}(\phi,\hat{\epsilon}_l,\hat{\mathbf{u}}_l)$ for the liquid phase, where $\phi$ is given by Eq. (\ref{lbm-phi-l}). The quantities $\hat{\epsilon}_\varphi$ and $\hat{\mathbf{u}}_\varphi$ can be computed by the function $f_\varphi$ which acts as an auxiliary vector, namely
\begin{equation}\label{momentzero}
\hat{\epsilon}_\varphi = \sum_{q = 0}^{Q-1} f_\varphi(q),
\end{equation}
\begin{equation}\label{momentfirst}
\hat{\mathbf{u}}_\varphi = \sum_{q = 0}^{Q-1} 
\mathbf{v}_q f_\varphi(q).
\end{equation}
It is important to highlight the difference between $\hat{\mathbf{u}}_\varphi$ and $\mathbf{u}_\varphi$ discussed in the previous sections. The latter is the velocity field normalized by the characteristic flow speed $U=L/T$, where $L\gg\lambda$ is the characteristic length scale of the flow field and $T\gg\tau$ is the characteristic time scale of the flow field. Hence the following relation holds: $\hat{\mathbf{u}}_\varphi = \mathbf{u}_\varphi\,(L/\lambda)\,(\tau/T)$. Hence it is not granted that the code output $\hat{\mathbf{u}}_\varphi$ converges to the physical solution $\mathbf{u}_\varphi$ when the mesh is refined, namely when $\lambda/L\to 0$. The strategy for tuning the code parameters such that this convergence is enforced is called the scaling. Typically the acoustic scaling, i.e. $\tau/T=\lambda/L$, and the diffusive scaling, i.e. $\tau/T=(\lambda/L)^2$, are the two most popular examples. Similar considerations hold also for $\hat{\epsilon}_\varphi$, but they require to proceed with the formal asymptotic expansion, which is done in the following. 

The term $\hat{\Omega}^F_\varphi$ is the forcing operator and it can be defined by the same linearized functional form, namely
\begin{equation}\label{forcingoperator}
f^{eq}_L(\psi,\hat{S},\hat{\mathbf{G}},q)=
w_q\left[\hat{S}\,\eta(\psi)+
\frac{\mathbf{v}_q\cdot\hat{\mathbf{G}}}{c_s^2}
\right].
\end{equation}
In particular, $\hat{\Omega}^F_g \equiv f^{eq}_L(\psi_g,\hat{S}_{g},\hat{\mathbf{G}}_{g})$ for the dispersed phase and $\hat{\Omega}^F_l=f^{eq}_L
(\psi_l,\hat{S}_{l},\hat{\mathbf{G}}_{l})$. In general, $\psi_g$ is different from $\psi_l$ and is different from $\phi$. Intuitively $\hat{S}_\varphi$ and $\hat{\mathbf{G}}_\varphi$ are the expressions computed using $\hat{\mathbf{u}}_\varphi$ (Boltzmann scaling) in the code.  

\subsubsection{Asymptotic analysis by the equivalent moment system}

Eqs. (\ref{lbmg}, \ref{lbml}) do not solve directly the target fluid equations, in the sense that the computed quantities, indicated with the hat notation $\hat{\cdot}$, do not converge automatically to the target fluid quantities if a proper scaling is not provided. The scaling is the set of rules used to update the input parameters once the mesh is refined. Hence an asymptotic analysis is needed in order to find out under which conditions these equations at least approximate the target fluid equations. Here we use the asymptotic method based on the equivalent moment system \cite{Junk2005676}, inspired by the moment method, which was first introduced to gas kinetic theory by H. Grad. The previous equations share the same structure, namely %
\begin{equation}\label{lbm}
f_\varphi(\hat{\mathbf{x}}+\mathbf{v}_q,
\hat{t}+1,q)=
f_\varphi(q)+\hat{\Omega}_\varphi(q)+
\hat{\Omega}_\varphi^F(q).
\end{equation}
Let us apply a Taylor expansion, namely
\begin{equation}\label{taylor}
f_\varphi(\hat{\mathbf{x}}+\mathbf{v}_q,
\hat{t}+1,q)-
f_\varphi(\hat{\mathbf{x}},
\hat{t},q)
=
\frac{\partial f_\varphi}{\partial \hat{t}}+
\mathbf{v}_q\cdot\hat{\nabla} f_\varphi+
\frac{1}{2}\,(\mathbf{v}_q\cdot\hat{\nabla})^2 f_\varphi
+\frac{1}{2}\,
\frac{\partial^2 f_\varphi}{\partial \hat{t}^2}+
\dots.
\end{equation}
In the previous equation, the unit of the space coordinate is the distance $\lambda$ between two neighboring lattices nodes (mean free path) and the unit of the time evolution is the time $\tau$ between two consecutive lattice collisions (mean collision time). Obviously, they are not appropriate as the characteristic scales for the flow field in the continuum limit, namely $L$ and $T$. Let us introduce the parameter $\varepsilon = \lambda/L \ll 1$ (mesh spacing) and let us reduce the corresponding time step quadratically, namely $\tau/T = \varepsilon^2$ (diffusive scaling), which increases quadratically the number of time steps, as well as the computational time needed to solve the problem. By using these assumptions, the following relations hold $\hat{t}=t/\varepsilon^2$ and $\hat{\mathbf{x}}=x/\varepsilon$, which lead to 
\begin{equation}\label{taylor2}
f_\varphi(\hat{\mathbf{x}}+\mathbf{v}_q,
\hat{t}+1,q)-
f_\varphi(\hat{\mathbf{x}},
\hat{t},q)
=
\varepsilon\,\mathbf{v}_q\cdot\nabla f_\varphi+
\frac{\varepsilon^2}{2}\,(\mathbf{v}_q\cdot\nabla)^2 f_\varphi
+\varepsilon^2\,\frac{\partial f_\varphi}{\partial t}
+O(\varepsilon^3).
\end{equation}
Neglecting terms $O(\varepsilon^3)$ yields
\begin{equation}\label{taylor3}
\varepsilon^2\,\frac{\partial f_\varphi}{\partial t}+
\varepsilon\,\mathbf{v}_q\cdot\nabla f_\varphi+
\frac{\varepsilon^2}{2}\,(\mathbf{v}_q\cdot\nabla)^2 f_\varphi
=
\hat{\Omega}_\varphi(q)+
\hat{\Omega}_\varphi^F(q).
\end{equation}
Taking the zero-th order, first order and second order moments of the previous equations yields
\begin{equation}\label{taylor-zeroth}
\varepsilon^2\,\frac{\partial \hat{\epsilon}_\varphi}{\partial t}+
\varepsilon\,\nabla\cdot\hat{\mathbf{u}}_\varphi+
\frac{\varepsilon^2}{2}\,\nabla\cdot\nabla\cdot\hat{\mathbf{\Pi}}_\varphi
= \hat{S}_\varphi,
\end{equation}
\begin{equation}\label{taylor-first}
\varepsilon^2\,\frac{\partial \hat{\mathbf{u}}_\varphi}{\partial t}+
\varepsilon\,\nabla\cdot\hat{\mathbf{\Pi}}_\varphi+
\frac{\varepsilon^2}{2}\,\nabla\cdot\nabla\cdot\hat{\mathbf{\Phi}}_\varphi
= \hat{\mathbf{G}}_\varphi,
\end{equation}
\begin{equation}\label{taylor-second}
\varepsilon^2\,\frac{\partial \hat{\mathbf{\Pi}}_\varphi}{\partial t}+
\varepsilon\,\nabla\cdot\hat{\mathbf{\Phi}}_\varphi+
\frac{\varepsilon^2}{2}\,\nabla\cdot\nabla\cdot
\sum_q \mathbf{v}_q\mathbf{v}_q\mathbf{v}_q\mathbf{v}_qf_\varphi(q)
= \omega_\varphi(\hat{\mathbf{\Pi}}^{eq}_\varphi-\hat{\mathbf{\Pi}}_\varphi)+\psi_\varphi c_s^2\,\hat{S}_\varphi\,\mathbf{I},
\end{equation}
where $\hat{\mathbf{\Pi}}_\varphi=\sum_{q = 0}^{Q-1} 
\mathbf{v}_q\mathbf{v}_q f_\varphi(q)$ is the second order tensor of the distribution function and $\hat{\mathbf{\Phi}}_\varphi=\sum_{q = 0}^{Q-1}\mathbf{v}_q\mathbf{v}_q\mathbf{v}_q f_\varphi(q)$ is the third order tensor. Moreover, by the definition given by Eq. (\ref{equilibriumI}), it is possible to compute
\begin{equation}\label{secondtensoreq}
\hat{\mathbf{\Pi}}^{eq}_\varphi = \hat{\phi}_\varphi c_s^2\,\hat{\epsilon}_\varphi\,\mathbf{I}+
\hat{\mathbf{u}}_\varphi\hat{\mathbf{u}}_\varphi,
\end{equation}
and 
\begin{equation}\label{thirdtensoreq-th}
(\hat{\mathbf{\Phi}}^{eq}_{\varphi})_{ijk} = 
c_s^2\left(
\hat{u}_{\varphi\,i}\delta_{jk}+
\hat{u}_{\varphi\,j}\delta_{ik}+
\hat{u}_{\varphi\,k}\delta_{ij}
\right).
\end{equation}

%\hl{The property} $\sum_q \mathbf{v}_q\mathbf{v}_q \hat{\Omega}^F_\varphi(q) =
%\sum_q \mathbf{v}_q\mathbf{v}_q f^{eq}_L(0,q)=\mathbf{0}$ has been used. 

In order to analyze this system of equations, it is now time to understand the impact of the previous assumptions on the scaling of the moments, namely how changing the mesh spacing impacts on the numerical values of the lattice moments (in Boltzmann scaling) computed by the code. The relation $\hat{\mathbf{u}}_\varphi = \mathbf{u}_\varphi\,(L/\lambda)\,(\tau/T)$ becomes $\hat{\mathbf{u}}_\varphi = \varepsilon\,\mathbf{u}_\varphi$. Consequently $\hat{S}_\varphi=\varepsilon^2\,S_\varphi$ because the source is the divergence of a velocity vector combination. Recalling the definitions given by Eq. (\ref{Gg}) and Eq. (\ref{Gl}), there are four terms in $\hat{\mathbf{G}}_{\varphi}$. According to the previous scaling, $\hat{S}_{\varphi}\hat{\mathbf{u}}_\varphi$ and $\hat{\bm{\sigma}}_\varphi \cdot \hat{\nabla}\alpha_g$ are automatically $O(\varepsilon^3)$. Hence it makes sense to scale in the same way also the two remaining terms. Adopting $\hat{\mathbf{g}}=\varepsilon^3\mathbf{g}$ means that, if $\varepsilon$ becomes half, then $\hat{\mathbf{g}}$ must become $1/2^3$ of the original value. Similarly, we adopt $\hat{\mathbf{F}}_{\varphi\overline{\varphi}}=\varepsilon^3\mathbf{F}_{\varphi\overline{\varphi}}$, where $\overline{\varphi}$ is the other phase in relation with phase $\varphi$, and $\hat{\mathbf{F}}_{\varphi}=\varepsilon^3\mathbf{F}_{\varphi}$. The first decision implies $\hat{K}_I=\varepsilon\,K_I$ because in Eq. (\ref{interphasemomentum}) there is already a quadratic dependence on some velocity. Putting together all these assumptions yield $\hat{\mathbf{G}}_\varphi=\varepsilon^3\,\mathbf{G}_\varphi$, which is consistent with the fact that a force induces an acceleration, namely $\hat{\mathbf{G}}_\varphi \propto \partial\hat{\mathbf{u}}_\varphi/\partial\hat{t}=\varepsilon^3\,\partial\mathbf{u}_\varphi/\partial t$. 

A system of moments can be truncated if we have some expectations about the high order moments. Actually, on a discrete lattice, the system of moments is automatically truncated because of the limited number of independent degrees of freedoms (typically up to some components of the fourth order moment) \cite{Junk2005676}. Let us now imagine the equation for the third order moment $\hat{\mathbf{\Phi}}_\varphi$, which is very similar to the previous equations. At the left hand side of this equation, there is a term which is proportional to $\varepsilon^2$ multiplied by the double divergence of the fifth order moment, which is ``odd'' with regard to the power of $\mathbf{v}_q$ and hence it scales as $\hat{\mathbf{u}}_\varphi=\varepsilon{\mathbf{u}}_\varphi$: altogether this term at the left hand side scales as $O(\varepsilon^3)$. At the right hand side, we have $\hat{\mathbf{\Phi}}^{eq}_\varphi - \hat{\mathbf{\Phi}}_\varphi$ and the third order moment of the forcing operator, which is proportional to $\varepsilon^3\,\mathbf{G}_\varphi$ and hence again $O(\varepsilon^3)$. Putting pieces together, the equation for the third order moment looks like $\hat{\mathbf{\Phi}}_\varphi = \hat{\mathbf{\Phi}}^{eq}_\varphi+O(\varepsilon^3)$. It is clear from Eq. (\ref{thirdtensoreq-th}) that $\hat{\mathbf{\Phi}}^{eq}_\varphi=\varepsilon{\mathbf{\Phi}}^{eq}_\varphi$ and consequently
\begin{equation}\label{taylor-third2}
\mathbf{\Phi}_\varphi = \mathbf{\Phi}^{eq}_\varphi+O(\varepsilon^2).
\end{equation}
Substituting these considerations in Eq. (\ref{taylor-first}) yields
\begin{equation}\label{taylor-first2}
\varepsilon^3\,\frac{\partial \mathbf{u}_\varphi}{\partial t}+
\varepsilon\,\nabla\cdot\hat{\mathbf{\Pi}}_\varphi+
\frac{\varepsilon^3}{2}\,\nabla\cdot\nabla\cdot\mathbf{\Phi}^{eq}_\varphi
= \varepsilon^3\,\mathbf{G}_\varphi+O(\varepsilon^5),
\end{equation}
which implies $\nabla\cdot\hat{\mathbf{\Pi}}_\varphi=\varepsilon^2\,\nabla\cdot\mathbf{\Pi}_\varphi$. Consequently
\begin{equation}\label{taylor-first3}
\frac{\partial \mathbf{u}_\varphi}{\partial t}+
\nabla\cdot\mathbf{\Pi}_\varphi+
\frac{1}{2}\,\nabla\cdot\nabla\cdot\mathbf{\Phi}^{eq}_\varphi
= \mathbf{G}_\varphi+O(\varepsilon^2).
\end{equation}
Similarly, the equation for the second order moment looks like $\hat{\mathbf{\Pi}}_\varphi = \hat{\mathbf{\Pi}}^{eq}_\varphi+O(\varepsilon^2) = \hat{\phi}_\varphi c_s^2\,\hat{\epsilon}_\varphi\,\mathbf{I} + O(\varepsilon^2)$. Using again the last relation implies $\nabla\cdot\hat{\mathbf{\Pi}}_\varphi=\varepsilon^2\,\nabla\cdot\mathbf{\Pi}_\varphi=O(\varepsilon^2)$ and then $\nabla(\hat{\phi}_\varphi c_s^2\,\hat{\epsilon}_\varphi)=O(\varepsilon^2)$ and equivalently $\nabla(\hat{\phi}_\varphi\hat{\epsilon}_\varphi)=O(\varepsilon^2)$. Moreover, let us assume that, in some portion of the domain boundary (typically the outflow), namely for $x\in \Omega_B$, the following boundary conditions hold: $\hat{\epsilon}_g\left.\right|_{B}=1$, $\hat{\epsilon}_l\left.\right|_{B}=1$ and consequently $\hat{\phi}\left.\right|_{B}=1$. These boundary conditions allow one to integrate the previous relation, which yields 
\begin{equation}\label{pressurescaling-g}
\hat{\epsilon}_g=
1+\varepsilon^2\,{\epsilon}_g,
\end{equation}
\begin{equation}\label{pressurescaling-phil}
\hat{\phi}\,\hat{\epsilon}_l=
1+(1/R)\,\left(\hat{\epsilon}_g-1\right)=
1+\varepsilon^2\,{\epsilon}_g/R.
\end{equation}
Assuming $\phi_g=1$ and $\phi_l=\left({\epsilon}_g/{\epsilon}_l\right)/R$, the relations given by Eq. (\ref{pressurescaling-g}) and Eq. (\ref{pressurescaling-phil}) can be expressed by one unique formula, namely 
\begin{equation}\label{pressurescaling}
\hat{\phi}_\varphi\hat{\epsilon}_\varphi=
1+\varepsilon^2\,\phi_\varphi{\epsilon}_\varphi.
\end{equation}
Concerning $\hat{\epsilon}_l$, introducing the relation $\nabla\cdot\hat{\mathbf{\Pi}}_l=O(\varepsilon^2)$ into Eq. (\ref{taylor-zeroth}) implies that $\partial \hat{\epsilon}_l/\partial t=O(1)$ or $\partial \hat{\epsilon}_l/\partial t=O(\varepsilon^2)$. Let us focus on the second case by assuming that the boundary conditions change smoothly in time. Taking into account the boundary condition $\hat{\epsilon}_l\left.\right|_{B}=1$ and integrating the relation $\partial \hat{\epsilon}_l/\partial t=O(\varepsilon^2)$ yield
\begin{equation}\label{pressurescaling-l}
\hat{\epsilon}_l=
1+\varepsilon^2\,{\epsilon}_l,
\end{equation}
and consequently
\begin{equation}\label{pressurescaling-phi}
\hat{\phi}=\frac{1+\varepsilon^2\,{\epsilon}_g/R}
{1+\varepsilon^2\,{\epsilon}_l}=
1 + \varepsilon^2\,\left({\epsilon}_g/R - {\epsilon}_l\right)+O(\varepsilon^4).
\end{equation}
Introducing Eq. (\ref{pressurescaling}), Eq. (\ref{pressurescaling-g}) and Eq. (\ref{pressurescaling-l}) into Eq. (\ref{taylor-zeroth}) yields
\begin{equation}\label{taylor-zeroth2}
\nabla\cdot\mathbf{u}_\varphi
= S_\varphi+O(\varepsilon^2).
\end{equation}

Eq. (\ref{taylor-first3}) and Eq. (\ref{taylor-zeroth2}) prove that the proposed scaling converges asymptotically with second order of convergence to some equations, but we still have to prove that they coincide with the target equations. First of all, we search for an approximation of $\mathbf{\Pi}_\varphi$. Recalling Eq. (\ref{taylor-second}), taking into account that the term with sequence of divergences of the fourth order moment scales as $\varepsilon^2\,O(\nabla(\hat{\phi}_\varphi\hat{\epsilon}_\varphi))=O(\varepsilon^4)$ because of Eq. (\ref{pressurescaling}) and recalling Eq. (\ref{taylor-zeroth2}), it follows
%
%\begin{equation}\label{taylor-second2}
%\varepsilon^2\,\nabla\cdot\mathbf{\Phi}^{eq}_\varphi
%= \omega_\varphi(\hat{\mathbf{\Pi}}^{eq}_\varphi-%\hat{\mathbf{\Pi}}_\varphi)+\varepsilon^2\psi_\varphi %c_s^2\,S_\varphi\,\mathbf{I}+O(\varepsilon^4).
%\end{equation}
%
%Using Eq. (\ref{thirdtensoreq}) and Eq. (\ref{taylor-zeroth2}) yield
%
\begin{equation}\label{taylor-second3}
\varepsilon^2\,\nabla\cdot\mathbf{\Phi}^{eq}_\varphi
= \omega_\varphi(\hat{\mathbf{\Pi}}^{eq}_\varphi-\hat{\mathbf{\Pi}}_\varphi)+\varepsilon^2\psi_\varphi c_s^2\,\nabla\cdot\mathbf{u}_\varphi\,\mathbf{I}+O(\varepsilon^4).
\end{equation}
Recalling that $\nabla\cdot\hat{\mathbf{\Pi}}_\varphi=\varepsilon^2\,\nabla\cdot\mathbf{\Pi}_\varphi$ yields
\begin{equation}\label{taylor-second4}
\nabla\cdot\mathbf{\Pi}_\varphi
= \nabla\cdot\mathbf{\Pi}^{eq}_\varphi-
\nabla\cdot\left[
\frac{1}{\omega_\varphi}\,\nabla\cdot\mathbf{\Phi}^{eq}_\varphi
-\frac{c_s^2\psi_\varphi }{\omega_\varphi}\,
\nabla\cdot\mathbf{u}_\varphi\,\mathbf{I}
\right]+O(\varepsilon^2).
\end{equation}
Substituting Eq. (\ref{taylor-second4}) into Eq. (\ref{taylor-first3}) yields
\begin{equation}\label{taylor-first4}
\frac{\partial \mathbf{u}_\varphi}{\partial t}+
\nabla\cdot\mathbf{\Pi}^{eq}_\varphi
=
\nabla\cdot\left[
\left(\frac{1}{\omega_\varphi}-\frac{1}{2}\right)\,\nabla\cdot\mathbf{\Phi}^{eq}_\varphi
-\frac{c_s^2\psi_\varphi }{\omega_\varphi}\,
\nabla\cdot\mathbf{u}_\varphi\,\mathbf{I}
\right]
+\mathbf{G}_\varphi+O(\varepsilon^2).
\end{equation}
From the definition given by Eq. (\ref{thirdtensoreq-th}) it is possible to prove that
\begin{equation}\label{thirdtensoreq}
\nabla\cdot\mathbf{\Phi}^{eq}_\varphi = 
c_s^2\,\left(\nabla\mathbf{u}_\varphi
+\nabla\mathbf{u}_\varphi^T
+\nabla\cdot\mathbf{u}_\varphi\,\mathbf{I}
\right).
\end{equation}
Let us add another ingredient by defining the kinematic viscosity as
\begin{equation}\label{effective-viscosity}
\nu_\varphi=
c_s^2
\left(\frac{1}{\omega_\varphi}-\frac{1}{2}\right),
\end{equation}
and by choosing the tunable parameter $\psi_\varphi$ in order to recover the right kinematic bulk viscosity $\xi_\varphi$ as
\begin{equation}\label{effective-bulk}
\psi_\varphi =\frac{\omega_\varphi}{c_s^2}
\left(\frac{5}{3}\,\nu_\varphi
-\xi_\varphi\right).
\end{equation}
It is important to note that both previous viscosity coefficients are normalized by the characteristic scales for the flow field in the continuum limit, namely by $L^2/T$ (in other words, the viscosity coefficients are the inverses of the corresponding Reynolds numbers). Substituting the previous expression into Eq. (\ref{taylor-first4}) yields
\begin{equation}\label{taylor-first5}
\frac{\partial \mathbf{u}_\varphi}{\partial t}+
\nabla\cdot\mathbf{\Pi}^{eq}_\varphi
=\nabla \cdot \bm{\sigma}_\varphi+ \mathbf{G}_\varphi+O(\varepsilon^2).
\end{equation}
Clearly Eq. (\ref{taylor-first5}) and Eq. (\ref{taylor-zeroth2}) approach the solution of the target equations in the asymptotic limit of $\varepsilon\to0$ with second order accuracy, thanks to the adopted diffusive scaling. 

These LBM schemes are not only useful to recover an approximated solution of the target equations, but they can be used to compute relevant quantities which will be used in the following. For example, it is possible to compute the viscous stress tensor per unit volume without further computational efforts. Comparing Eq. (\ref{taylor-first5}) and Eq. (\ref{taylor-first4}), it is possible to derive
\begin{equation}\label{lbm-sigma}
\bm{\sigma}_\varphi=
\left(\frac{1}{\omega_\varphi}-\frac{1}{2}\right)\,\nabla\cdot\mathbf{\Phi}^{eq}_\varphi
-\frac{c_s^2\psi_\varphi }{\omega_\varphi}\,
\nabla\cdot\mathbf{u}_\varphi\,\mathbf{I}
+O(\varepsilon^2).
\end{equation}
From the numerical point of view, it is possible to compute $\nabla\cdot\mathbf{\Phi}^{eq}_\varphi$ by Eq. (\ref{taylor-second3}) and $\nabla\cdot\mathbf{u}_\varphi$ by Eq. (\ref{taylor-zeroth2}), namely
\begin{equation}\label{lbm-sigma2}
\bm{\sigma}_\varphi=
\left(\frac{1}{\omega_\varphi}-\frac{1}{2}\right)\,\left[\omega_\varphi\,\left(\mathbf{\Pi}^{eq}_\varphi-\mathbf{\Pi}_\varphi\right)+\psi_\varphi\,c_s^2\,S_\varphi\,\mathbf{I}\right]
-\frac{c_s^2\psi_\varphi }{\omega_\varphi}\,
S_\varphi\,\mathbf{I}
+O(\varepsilon^2),
\end{equation}
or equivalently
\begin{equation}\label{lbm-sigma3}
\bm{\sigma}_\varphi=
\left(\frac{1}{\omega_\varphi}-\frac{1}{2}\right)\,\omega_\varphi\,\left(\mathbf{\Pi}^{eq}_\varphi-\mathbf{\Pi}_\varphi\right)
-\frac{c_s^2\psi_\varphi }{2}\,
S_\varphi\,\mathbf{I}
+O(\varepsilon^2),
\end{equation}
which do not require to compute explicitly additional space derivatives. 

%\begin{equation}\label{taylor-second4bis}
%\left(\frac{1}{\omega_\varphi}-\frac{1}{2}\right)\left[-\omega_\varphi\,\left(\mathbf{\Pi}_\varphi-\mathbf{\Pi}^{eq}_\varphi\right)+\psi_\varphi\,c_s^2\,\nabla\cdot\mathbf{u}_\varphi\,\mathbf{I}\right]
%-\frac{c_s^2\psi_\varphi }{\omega_\varphi}\,
%\nabla\cdot\mathbf{u}_\varphi\,\mathbf{I}=\nu_\varphi\,\left[\nabla\mathbf{u}_\varphi
%+\nabla\mathbf{u}_\varphi^T
%+\nabla\cdot\mathbf{u}_\varphi\,\mathbf{I}
%\right]-\frac{c_s^2\psi_\varphi }{\omega_\varphi}\,
%\nabla\cdot\mathbf{u}_\varphi\,\mathbf{I}
%+O(\varepsilon^2).
%\end{equation}

\subsection{LBM schemes for solving the volume fraction equation for each phase and their asymptotic analysis}

In addition of the two LBM schemes discussed above for the momentum equations, one needs to solve the equation for the volume fraction, but ensuring that the volume fraction is bounded between zero and one, namely $0\leq\alpha_\varphi\leq 1$. There are many advanced numerical techniques to ensure this condition. One simple alternative may be that proposed by Spalding  \cite{Spalding1985}. Essentially the idea is to solve both Eq. (\ref{voidfraction}) and Eq. (\ref{continuity_l}) for $\alpha_g$ and $\alpha_l$, and then perform a proper renormalisation. In the LBM context, let us consider the following two LBM schemes for solving the volume fraction equations given by Eq. (\ref{voidfraction}) and Eq. (\ref{continuity_l}), formulated in terms of the corresponding particle distribution functions $f_{\alpha g}$ and $f_{\alpha l}$, namely
\begin{equation}\label{lbmalpha}
f_{\alpha\varphi}(\hat{\mathbf{x}}+\mathbf{v}_q,
\hat{t}+1,q)=
f_{\alpha\varphi}(q)+
\hat{\omega}_{\alpha\varphi}
\left[f^{eq}_{\alpha\varphi}(q) - f_{\alpha\varphi}\right],
\end{equation}
where $\hat{\omega}_{\alpha\varphi}$ (for the volume phase fraction $\alpha_\varphi$) is the relaxation frequency, divided by the frequency $1/\tau$, which drives the relaxation of the distribution function towards the equilibrium $f^{eq}_{\alpha\varphi}$. The fact that we report the relaxation frequency $\hat{\omega}_{\alpha\varphi}$ with the hat notation $\hat{\cdot}$ means that we allow the possibility to adjust the relaxation frequency during the mesh refinement according the adopted scaling. The equilibrium distribution is formulated by means of the same functional form, namely $f^{eq}_{\alpha\varphi} \equiv f^{eq}(\hat{\alpha}_{\varphi},\hat{\mathbf{u}}_\varphi)$, where
\begin{equation}\label{equilibrium}
f^{eq}(\hat{\alpha},\hat{\mathbf{u}},q) = w_q\,\hat{\alpha}\,
\left[1+
\frac{\mathbf{v}_q\cdot\hat{\mathbf{u}}}{c_s^2}
+\frac{(\mathbf{v}_q\mathbf{v}_q-c_s^2\,\mathbf{I}):
\hat{\mathbf{u}}\hat{\mathbf{u}}}{2 c_s^4}\right],
\end{equation}
and
\begin{equation}\label{momentzeroapha}
\hat{\alpha}_\varphi = \sum_{q = 0}^{Q-1} f_{\alpha\varphi}(q).
\end{equation}
It is very important to highlight that the velocity field $\hat{\mathbf{u}}_\varphi$ is computed by the LBM schemes discussed in the previous section, and hence the first order moment $\hat{\alpha}_\varphi\hat{\bm{\upsilon}}_\varphi$ of the distribution $f_{\alpha\varphi}$ is not conserved, namely
\begin{equation}\label{momentfirstalpha}
\hat{\mathbf{u}}_\varphi = \sum_{q = 0}^{Q-1} 
\mathbf{v}_q f_\varphi(q)\neq
\sum_{q = 0}^{Q-1} 
\mathbf{v}_q f_{\alpha\varphi}(q)=
\hat{\alpha}_\varphi\hat{\bm{\upsilon}}_\varphi.
\end{equation}

As discussed in the previous section, an asymptotic analysis is needed in order to find out under which conditions, these equations at least approximate the target volume fraction equations. Again we use here the asymptotic method based on the equivalent moment system \cite{Junk2005676}. Let us apply a Taylor expansion and let us assume a diffusive scaling, namely $\hat{t}=t/\varepsilon^2$ and $\hat{\mathbf{x}}=x/\varepsilon$, which leads to 
\begin{equation}\label{taylor2alpha}
f_{\alpha\varphi}(\hat{\mathbf{x}}+\mathbf{v}_q,
\hat{t}+1,q)-
f_{\alpha\varphi}(\hat{\mathbf{x}},
\hat{t},q)
=
\varepsilon\,\mathbf{v}_q\cdot\nabla f_{\alpha\varphi}+
\frac{\varepsilon^2}{2}\,(\mathbf{v}_q\cdot\nabla)^2 f_{\alpha\varphi}
+\varepsilon^2\,\frac{\partial f_{\alpha\varphi}}{\partial t}
+O(\varepsilon^3).
\end{equation}
Neglecting terms $O(\varepsilon^3)$ yields
\begin{equation}\label{taylor3alpha}
\varepsilon^2\,\frac{\partial f_{\alpha\varphi}}{\partial t}+
\varepsilon\,\mathbf{v}_q\cdot\nabla f_{\alpha\varphi}+
\frac{\varepsilon^2}{2}\,(\mathbf{v}_q\cdot\nabla)^2 f_{\alpha\varphi}
=
\hat{\omega}_{\alpha\varphi}
\left[f^{eq}_{\alpha\varphi}(q) - f_{\alpha\varphi}\right].
\end{equation}

Taking the zero-th order, first order and second order moments of the previous equations yields 
\begin{equation}\label{taylor-zerothalpha}
\varepsilon^2\,\frac{\partial \hat{\alpha}_\varphi}{\partial t}+
\varepsilon\,\nabla\cdot\left(\hat{\alpha}_\varphi\hat{\bm{\upsilon}}_\varphi\right)+
\frac{\varepsilon^2}{2}\,\nabla\cdot\nabla\cdot\hat{\mathbf{\Pi}}_{\alpha\varphi}
= 0,
\end{equation}
\begin{equation}\label{taylor-firstalpha}
\varepsilon^2\,\frac{\partial \left(\hat{\alpha}_\varphi\hat{\bm{\upsilon}}_\varphi\right)}{\partial t}+
\varepsilon\,\nabla\cdot\hat{\mathbf{\Pi}}_{\alpha\varphi}+
\frac{\varepsilon^2}{2}\,\nabla\cdot\nabla\cdot\hat{\mathbf{\Phi}}_{\alpha\varphi}
= \hat{\omega}_{\alpha\varphi}\,\hat{\alpha}_\varphi\left(\hat{\mathbf{u}}_\varphi-\hat{\bm{\upsilon}}_\varphi\right),
\end{equation}
\begin{equation}\label{taylor-secondalpha}
\varepsilon^2\,\frac{\partial \hat{\mathbf{\Pi}}_{\alpha\varphi}}{\partial t}+
\varepsilon\,\nabla\cdot\hat{\mathbf{\Phi}}_{\alpha\varphi}+
O(\varepsilon^2)
= \hat{\omega}_{\alpha\varphi}(\hat{\mathbf{\Pi}}^{eq}_{\alpha\varphi}-\hat{\mathbf{\Pi}}_{\alpha\varphi}),
\end{equation}
where $\hat{\mathbf{\Pi}}_{\alpha\varphi}=\sum_{q = 0}^{Q-1} 
\mathbf{v}_q\mathbf{v}_q f_{\alpha\varphi}(q)$ is the second order tensor of the distribution function and $\hat{\mathbf{\Phi}}_{\alpha\varphi}=\sum_{q = 0}^{Q-1}\mathbf{v}_q\mathbf{v}_q\mathbf{v}_q f_{\alpha\varphi}(q)$ is the third order tensor. Moreover, by the definition given by Eq. (\ref{equilibrium}), it is possible to compute
\begin{equation}\label{secondtensoreqalpha}
\hat{\mathbf{\Pi}}^{eq}_{\alpha\varphi} =  c_s^2\,\hat{\alpha}_{\varphi}\,\mathbf{I}+
\hat{\alpha}_{\varphi}\,\hat{\mathbf{u}}_\varphi\hat{\mathbf{u}}_\varphi,
\end{equation}
and 
\begin{equation}\label{thirdtensoreq-thalpha}
(\hat{\mathbf{\Phi}}^{eq}_{\alpha\varphi})_{ijk} = 
c_s^2\,\hat{\alpha}_{\varphi}\,\left(
\hat{u}_{\varphi\,i}\delta_{jk}+
\hat{u}_{\varphi\,j}\delta_{ik}+
\hat{u}_{\varphi\,k}\delta_{ij}
\right).
\end{equation}

As discussed in the previous section, in order to analyze this system of equations, one needs to understand the impact of the diffusive scaling on the scaling of the moments. Let us suppose not to scale the volume fraction, namely $\hat{\alpha}_\varphi={\alpha}_\varphi$, which has an impact on all even moments of the distribution function. On the other hand, let us keep the same scaling as before for the odd moments, because they depend on the imposed velocity $\hat{\mathbf{u}}_\varphi = \varepsilon\,\mathbf{u}_\varphi$. This implies
\begin{equation}\label{secondtensoreqalpha2}
\hat{\mathbf{\Pi}}^{eq}_{\alpha\varphi} =  c_s^2\,{\alpha}_{\varphi}\,\mathbf{I}+
\varepsilon^2{\alpha}_{\varphi}\,{\mathbf{u}}_\varphi{\mathbf{u}}_\varphi,
\end{equation}
and 
\begin{equation}\label{thirdtensoreq-thalpha2}
(\hat{\mathbf{\Phi}}^{eq}_{\alpha\varphi})_{ijk} = 
\varepsilon\,c_s^2\,{\alpha}_{\varphi}\,\left(
{u}_{\varphi\,i}\delta_{jk}+
{u}_{\varphi\,j}\delta_{ik}+
{u}_{\varphi\,k}\delta_{ij}
\right).
\end{equation}

A system of moments can be truncated if we have some expectations about the high order moments. Let us now imagine the equation for the third order moment $\hat{\mathbf{\Phi}}_{\alpha\varphi}$, which can be analyzed with very similar arguments discussed in the previous section. These arguments lead to conclude that the equation for the third order moment looks like $\hat{\mathbf{\Phi}}_{\alpha\varphi} = \hat{\mathbf{\Phi}}^{eq}_{\alpha\varphi}+O(\varepsilon^3)$ and consequently $\hat{\mathbf{\Phi}}_{\alpha\varphi} = \varepsilon{\mathbf{\Phi}}^{eq}_{\alpha\varphi}+O(\varepsilon^3)$, because of Eq. (\ref{thirdtensoreq-thalpha2}). Substituting the last result in Eq. (\ref{taylor-secondalpha}) and taking into account Eq. (\ref{secondtensoreqalpha2}) yields
\begin{equation}\label{taylor-secondalpha2}
O(\varepsilon^2)
-\hat{\omega}_{\alpha\varphi}
\varepsilon^2{\alpha}_{\varphi}\,{\mathbf{u}}_\varphi{\mathbf{u}}_\varphi
= \hat{\omega}_{\alpha\varphi}\left(c_s^2\,{\alpha}_{\varphi}\,\mathbf{I}-\hat{\mathbf{\Pi}}_{\alpha\varphi}\right),
\end{equation}
or equivalently $\hat{\mathbf{\Pi}}_{\alpha\varphi}=c_s^2\,{\alpha}_{\varphi}\,\mathbf{I}+O(\varepsilon^2)$. Substituting this result in Eq. (\ref{taylor-zerothalpha}) yields
\begin{equation}\label{taylor-zerothalpha2}
\varepsilon^2\,\frac{\partial \hat{\alpha}_\varphi}{\partial t}+
\varepsilon\,\nabla\cdot\left(\hat{\alpha}_\varphi\hat{\bm{\upsilon}}_\varphi\right)+
\frac{\varepsilon^2}{2}\,c_s^2\,\nabla^2\alpha_{\varphi}
= O(\varepsilon^4),
\end{equation}
which proves that $\hat{\bm{\upsilon}}_\varphi=\varepsilon{\bm{\upsilon}}_\varphi$. Substituing this result into Eq. (\ref{taylor-firstalpha}) yields
\begin{equation}\label{taylor-firstalpha2}
c_s^2\,\nabla\alpha_{\varphi}
= \hat{\omega}_{\alpha\varphi}\,{\alpha}_\varphi\left({\mathbf{u}}_\varphi-{\bm{\upsilon}}_\varphi\right)+
O(\varepsilon^2).
\end{equation}
Substituting the previous equation into Eq. (\ref{taylor-zerothalpha2}) yields
\begin{equation}\label{taylor-zerothalpha3}
\frac{\partial {\alpha}_\varphi}{\partial t}+
\nabla\cdot\left({\alpha}_\varphi{\mathbf{u}}_\varphi\right)
= \hat{\chi}_{\alpha\varphi}\,\nabla^2\alpha_{\varphi}
+O(\varepsilon^2),
\end{equation}
where
\begin{equation}\label{diffusivity}
\hat{\chi}_{\alpha\varphi}=
c_s^2\,\left(\frac{1}{\hat{\omega}_{\alpha\varphi}}
-\frac{1}{2}\right).
\end{equation}
Adopting the following scaling 
\begin{equation}\label{diffusivity2}
\hat{\omega}_{\alpha\varphi} = 
\frac{2}{1+2\,\varepsilon^2\,\chi_{\alpha\varphi}/c_s^2},
\end{equation}
implies $\hat{\chi}_{\alpha\varphi}=\varepsilon^2\, \chi_{\alpha\varphi}$ and consequently
\begin{equation}\label{taylor-zerothalpha4}
\frac{\partial {\alpha}_\varphi}{\partial t}+
\nabla\cdot\left({\alpha}_\varphi{\mathbf{u}}_\varphi\right)
= O(\varepsilon^2).
\end{equation}
Clearly the previous equation approaches the solution of the target equations given by Eq. (\ref{voidfraction}) and Eq. (\ref{continuity_l}) in the asymptotic limit of $\varepsilon\to0$ with second order accuracy, thanks to the
adopted diffusive scaling. 

These LBM schemes are not only useful to recover an approximated solution of the target equations, but they can be also used to compute other relevant quantities. For example, it is possible to compute the volume fraction gradient without further computational efforts using Eq. (\ref{taylor-firstalpha2}), namely
\begin{equation}\label{gradlnalpha}
\frac{1}{\alpha_{\varphi}}\nabla\alpha_{\varphi}
=\frac{\hat{\omega}_{\alpha\varphi}}{c_s^2}
\,\left({\mathbf{u}}_\varphi-{\bm{\upsilon}}_\varphi\right)+
O(\varepsilon^2).
\end{equation}
In particular Eq. (\ref{lbm-sigma3}) and Eq. (\ref{gradlnalpha}) can be used to compute the second term of the force $\mathbf{G}_g$ given by Eq. (\ref{Gg}), as well as the second term of the force $\mathbf{G}_l$ given by Eq. (\ref{Gl}), namely
\begin{equation}\label{forcealpha}
\frac{1}{\alpha_\varphi}\,\bm{\sigma}_\varphi \cdot \nabla\alpha_\varphi=\frac{\hat{\omega}_{\alpha\varphi}}{c_s^2}
\left[\left(1-\frac{\omega_\varphi}{2}\right)
\left(\mathbf{\Pi}^{eq}_\varphi-\mathbf{\Pi}_\varphi\right)
-\frac{c_s^2\psi_\varphi }{2}\,
S_\varphi\,\mathbf{I}\right]\cdot
\left({\mathbf{u}}_\varphi-{\bm{\upsilon}}_\varphi\right).
\end{equation}

One last remark is about the need of ensuring that volume fraction $\alpha_\varphi$ is bounded between zero and one, namely $0\leq\alpha_\varphi\leq 1$. Among many available numerical techniques to ensure this condition, we follow here that proposed by Spalding  \cite{Spalding1985}. Essentially the idea is to solve both Eq. (\ref{voidfraction}) and Eq. (\ref{continuity_l}) by the discussed LBM schemes for $f_{\alpha\varphi}$, compute the raw volume fraction ${\alpha}_\varphi = \sum_{q = 0}^{Q-1} f_{\alpha\varphi}(q)$ and then perform the following renormalisation at the beginning of every collision \& streaming cycle 
\begin{equation}\label{spalding}
    {\alpha}_\varphi^B \equiv \frac{B(\alpha_\varphi)}{B(\alpha_g)+B(\alpha_l)},
\end{equation}
where $B(x) = x\,H(x)\,H(1-x)+H(x-1)$ and $H(x)$ is the Heaviside step function. Of course, ${\alpha}_l^B = 1-{\alpha}_g^B$. 

Similarly, we propose a second correction for ensuring that $\nabla{\alpha}_g = -\nabla{\alpha}_l$, when using Eq. (\ref{gradlnalpha}). Let us suppose that the vector $\bm{\delta} = \nabla{\alpha}_g+\nabla{\alpha}_l$ has a non-zero modulus, i.e. $\|\bm{\delta}\|\neq 0$. In this case, the second Spalding-like correction is
\begin{equation}\label{spalding2}
    \nabla{\alpha}_\varphi^G \equiv \nabla{\alpha}_\varphi-\bm{\delta}/2.
\end{equation}
Of course $\nabla{\alpha}_l^G + \nabla{\alpha}_g^G = 0$. In the following, we will drop the superscript $B$ in Eq. (\ref{spalding}), as well as the superscript $G$ in Eq. (\ref{spalding2}), for the sake of simplicity and without risk of confusion.

\subsection{LBM schemes for solving the phase continuity source for each phase and their asymptotic analysis}

In addition of the two LBM schemes discussed above for the momentum equations and the two LBM schemes for the volume fractions, one needs to compute the phase continuity sources, namely $S_g$ in Eq. (\ref{lbm-cont-g}) and $S_l$ in Eq. (\ref{lbm-cont-l}). Let us consider the following two schemes for computing the phase continuity sources, formulated in terms of the corresponding particle distribution functions $f_{\beta g}$ and $f_{\beta l}$, namely
\begin{equation}\label{lbmbeta}
f_{\beta\varphi}(\hat{\mathbf{x}}+\mathbf{v}_q,
\hat{t}+1,q)=
f^{eq}_{\beta\varphi}(\hat{\mathbf{x}},
\hat{t},q),
\end{equation}
where $f^{eq}_{\beta\varphi} \equiv f^{eq}(\hat{\beta}_{\varphi},\hat{\mathbf{c}}_\varphi)$ and $f^{eq}$ is the functional form given by Eq. (\ref{equilibrium}). The zero order moment $\hat{\beta}_\varphi$ is given by 
\begin{equation}\label{momentzerobeta}
\hat{\beta}_\varphi = \sum_{q = 0}^{Q-1} f_{\beta\varphi}(q)
=\sum_{q = 0}^{Q-1} f_{\beta\varphi}^{eq}(q),
\end{equation}
and the imposed velocity field is given by
\begin{equation}\label{velocitybeta}
\hat{\mathbf{c}}_\varphi = \frac{1}{\hat{\beta}_\varphi}\,
\hat{\alpha}_{\overline{\varphi}} \left(\hat{\mathbf{u}}_\varphi
-\hat{\mathbf{u}}_{\overline{\varphi}}\right).
\end{equation}
The velocity field $\hat{\mathbf{c}}_\varphi$ is computed by the LBM schemes discussed in the previous sections and, by definition, the following relation holds
\begin{equation}\label{momentfirstbetaeq}
\sum_{q = 0}^{Q-1} 
\mathbf{v}_q f_{\alpha\varphi}^{eq}(q)=
\hat{\beta}_\varphi\hat{\bm{c}}_\varphi=
\hat{\alpha}_{\overline{\varphi}} \left(\hat{\mathbf{u}}_\varphi
-\hat{\mathbf{u}}_{\overline{\varphi}}\right).
\end{equation}
It is worth noting that the continuity sources which we are trying to compute are related to the fluxes $\hat{\beta}_\varphi\hat{\bm{c}}_\varphi$, namely $\hat{S}_\varphi=\hat{\nabla}\cdot(\hat{\beta}_\varphi\hat{\bm{c}}_\varphi)=\hat{\nabla}\cdot\left[\hat{\alpha}_{\overline{\varphi}} \left(\hat{\mathbf{u}}_\varphi
-\hat{\mathbf{u}}_{\overline{\varphi}}\right)\right]$.

Eq. (\ref{lbmbeta}) is the lattice kinetic scheme (LKS) \cite{Inamuro2002477} and is a particular case of the link-wise artificial compressibility method (link-wise ACM) \cite{Asinari20125109}. As in LBM schemes, also a link-wise ACM scheme proceeds by a sequence of collision \& streaming cycles. Eq. (\ref{lbmbeta}) assumes collision first and then streaming, but the scheme can also be rationalized by inverting this order, namely
\begin{equation}\label{lbmbeta2}
f_{\beta\varphi}(\hat{\mathbf{x}},
\hat{t}+1,q)=
f^{eq}_{\beta\varphi}(\hat{\mathbf{x}}-\mathbf{v}_q,\hat{t},q).
\end{equation}
Let us apply the Taylor expansion to the previous expression, namely
\begin{equation}\label{lbmbetataylor}
f_{\beta\varphi}(\hat{\mathbf{x}},
\hat{t}+1,q)-f^{eq}_{\beta\varphi}(\hat{\mathbf{x}},\hat{t},q)=
-\varepsilon\,\mathbf{v}_q\cdot\nabla f_{\beta\varphi}^{eq}+
\frac{\varepsilon^2}{2}\,(\mathbf{v}_q\cdot\nabla)^2 f_{\beta\varphi}^{eq}
-\frac{\varepsilon^3}{6}\,(\mathbf{v}_q\cdot\nabla)^3 f_{\beta\varphi}^{eq}
+O(\varepsilon^4).
\end{equation}
Computing the zero-th order moment of the previous expression and applying the scaling to the known quantities yields
\begin{equation}\label{lbmbetataylorzeroth}
\hat{\beta}_{\varphi}(\hat{\mathbf{x}},
\hat{t}+1)-\hat{\beta}_{\varphi}(\hat{\mathbf{x}},\hat{t})=
-\varepsilon^2\,S_\varphi+
\frac{\varepsilon^2}{2}\,
{\nabla}\cdot{\nabla}\cdot\hat{\mathbf{\Pi}}_{\beta\varphi}^{eq}
-\frac{\varepsilon^3}{6}\,
{\nabla}\cdot{\nabla}\cdot{\nabla}\cdot
\hat{\mathbf{\Phi}}^{eq}_{\beta\varphi}
+O(\varepsilon^4),
\end{equation}
where 
\begin{equation}\label{secondtensoreqbeta2}
\hat{\mathbf{\Pi}}^{eq}_{\beta\varphi} =  \hat{\beta}_{\varphi}\,\left(c_s^2\,\mathbf{I}+
\varepsilon^2\,{\mathbf{u}}_\varphi{\mathbf{u}}_\varphi
\right)=\hat{\beta}_{\varphi}\,c_s^2\,\mathbf{I}
+O(\varepsilon^2),
\end{equation}
and 
\begin{equation}\label{thirdtensoreq-thbeta2}
(\hat{\mathbf{\Phi}}^{eq}_{\beta\varphi})_{ijk} = 
\varepsilon\,c_s^2\,\hat{\beta}_{\varphi}\left(
{u}_{\varphi\,i}\delta_{jk}+
{u}_{\varphi\,j}\delta_{ik}+
{u}_{\varphi\,k}\delta_{ij}
\right).
\end{equation}
Substituting the previous expressions into Eq. (\ref{lbmbetataylorzeroth}) yields
\begin{equation}\label{lbmbetataylorzeroth2}
\hat{\beta}_{\varphi}(\hat{\mathbf{x}},
\hat{t}+1)-\hat{\beta}_{\varphi}(\hat{\mathbf{x}},\hat{t})=
-\varepsilon^2\,S_\varphi+
\varepsilon^2\,\frac{c_s^2}{2}\,
{\nabla}^2\hat{\beta}_{\varphi}(\hat{\mathbf{x}},\hat{t})
+O(\varepsilon^4).
\end{equation}
It is important to remind that $\hat{\beta}_{\varphi}$ are just auxiliary functions without physical interest, with the exception of their derivatives, as it will be clarified in the following. Assuming $\hat{\beta}_{\varphi}(\hat{\mathbf{x}},\hat{t})=1$ at the beginning of every time step yields

\begin{equation}\label{lbmbetataylorzeroth3}
\hat{\beta}_{\varphi}(\hat{\mathbf{x}},
\hat{t}+1)-1=-\varepsilon^2\,S_\varphi
+O(\varepsilon^4),
\end{equation}
or equivalently
\begin{equation}\label{lbmbetataylorzeroth4}
\frac{\partial \hat{\beta}_{\varphi}}{\partial t} = -S_\varphi
+O(\varepsilon^2).
\end{equation}
In particular, Eq. (\ref{lbmbetataylorzeroth3}) can be reformulated as
\begin{equation}\label{lbmbetataylorzeroth5}
\hat{S}_\varphi=1-\hat{\beta}_{\varphi}(\hat{\mathbf{x}},
\hat{t}+1)
+O(\varepsilon^4),
\end{equation}
which is particularly useful in this context, because it allows one to compute $\hat{S}_\varphi$ as coded into the LBM schemes, without explicit finite difference formulas.

\section{Numerical validation}

\subsection{One dimensional test case}

We consider a one-dimensional vertical tube with a generic two-phase flow as a test case to validate the proposed computational methodology. Although simplified, the one-dimensional vertical tube offers a foundational model for analyzing multiphase flow behavior in systems such as bubble column reactors, where gas-liquid interactions are key to reactor performance. Additionally, this configuration is pertinent for studying flow dynamics in natural circulation loops driven by gas injection, a mechanism critical in systems like molten salt nuclear reactors. This test case provides a controlled yet insightful context, enabling us to investigate flow patterns, pressure gradients, and phase interactions relevant to both chemical and nuclear reactor environments.

Let us suppose that the vertical tube is aligned with the axis identified by the unit vector $\mathbf{e}_x$, where $\|\mathbf{e}_x\|=1$. In this case, the gravitational acceleration field is $\mathbf{g}=-g\,\mathbf{e}_x$, where $g$ is the standard acceleration of gravity. It is sometimes convenient to remove the hydrostatic pressure from the pressure in order to make the buoyancy effect explicit, by defining a new quantity 
\begin{equation}\label{pressureprime}
p'=p+\rho_l^0\,g\,x,
\end{equation}
where $x$ is the coordinate along the unit vector $\mathbf{e}_x$. Let us introduce the kinematic pressure $p_k$ defined as
\begin{equation}\label{kinematicpressure}
p_k = \frac{p'}{\rho_g^0} = \frac{p+\rho_l^0\,g\,x}{\rho_g^0}
 = \frac{p}{\rho_g^0}+R\,g\,x,
\end{equation}
and consequently
\begin{equation}\label{kinematicpressure2}
\frac{1}{R}\,p_k = \frac{p}{\rho_l^0}+g\,x.
\end{equation}
The previous definitions lead to $-\partial_x(p/\rho_g^0)-g=-\partial_x p_k+(R-1)\,g$ for the dispersed phase and $-\partial_x(p/\rho_l^0)-g=-(1/R)\,\partial_x p_k$ for the liquid phase. Hence, using the kinematic pressure $p_k$, instead of the original pressure $p$, allows one to consider only one gravitational acceleration in the dispersed phase, i.e. the term $(R-1)\,g$, which is the buoyancy acceleration. 

Reformulating the artificially compressible equations in terms of the kinematic pressure implies that Eq. (\ref{sigmag}) must be redefined as
\begin{equation}\label{sigmag1D}
    \epsilon_g(p) = \frac{p'}{c_s^2\rho_g^0} = \frac{p_k}{c_s^2},
\end{equation}
without modifying the proposed methodology. As already pointed out, for having Eq. (\ref{lbm-cont-g}) valid also during the dynamics, the following condition must hold $\partial_t \epsilon_g = O(h^2)$. This means that the first derivative in Eq. (\ref{lbm-cont-g}) is not mesh independent. Assuming $h\equiv\varepsilon$ is possible to make explicit such dependence on the mesh spacing by intending $\partial_t \epsilon_g$ more precisely as $\partial_t \hat{\epsilon}_g$ into Eq. (\ref{lbm-cont-g}), since $\partial_t \hat{\epsilon}_g = \varepsilon^2\,\partial_t \epsilon_g$ according to Eq. (\ref{pressurescaling-g}). Recalling that $\hat{\epsilon}_g=\hat{p}_k/c_s^2$ where $\hat{p}_k = c_s^2+\varepsilon^2\,p_k$, Eq. (\ref{lbm-cont-g}), in the one-dimensional case, becomes
\begin{equation}\label{lbm-cont-g1D}
    \frac{\varepsilon^2}{c_s^2}\,\frac{\partial {p}_k}{\partial t}+
    \frac{\partial {u}_g}{\partial x} = 
    \frac{\partial }{\partial x} [\alpha_l ({u}_g - {u}_l)]
    = S_{g}.
\end{equation}
Substituting the previous quantities into Eq. (\ref{lbm-mom-g}), Eq. (\ref{Gg}), Eq. (\ref{lbm-mom-l}) and Eq. (\ref{Gl}), in the one-dimensional case, yields
\begin{equation}\label{lbm-mom-g1D}
    \frac{\partial {u}_g}{\partial t} + \frac{\partial}{\partial x}({u}_g {u}_g) = - \frac{\partial p_k}{\partial x}  + \frac{\partial {\sigma}_g}{\partial x} + {G}_{g},
\end{equation}
where
\begin{equation}\label{Gg1D}
    {G}_{g} = S_{g}\,{u}_g + 
    \frac{1}{\alpha_g}\,{\sigma}_g \frac{\partial \alpha_g}{\partial x} + (R-1)\,{g} + \frac{1}{\alpha_g\rho_g^0}({F}_{gl}+{F}_{g}).
\end{equation}
Similarly, it holds:
\begin{equation}\label{lbm-mom-l1D}
    \frac{\partial {u}_l}{\partial t} + \frac{\partial}{\partial x} ({u}_l {u}_l) = - \frac{1}{R}\,\frac{\partial p_k}{\partial x} + \frac{\partial {\sigma}_l}{\partial x} + {G}_{l},
\end{equation}
where
\begin{equation}\label{Gl1D}
    {G}_{l} = S_{l}\,{u}_l + \frac{1}{\alpha_l}\,{\sigma}_l \frac{\partial \alpha_l}{\partial x}+
    \frac{1}{\alpha_l\rho_l^0}({F}_{lg}+{F}_{l}).
\end{equation}
In the one-dimensional case, the viscous stress tensor per unit volume defined by Eq. (\ref{stresstensor}), assuming $\xi_\varphi=-(1/3)\,\nu_\varphi$ (which is a clear indication of the degeneracy of the one-dimensional case because it implies negative kinematic bulk viscosity), becomes
\begin{equation}\label{stresstensor1D}
    {\sigma}_\varphi = 
    \left(\frac{4}{3}\,\nu_\varphi+\xi_\varphi\right)\,
    \frac{\partial {u}_\varphi}{\partial x}=
    \nu_\varphi\,
    \frac{\partial {u}_\varphi}{\partial x}.
\end{equation}
Concerning the interphase momentum exchange, in the one-dimensional case, Eq. (\ref{interphasemomentum}) becomes
\begin{equation}\label{interphasemomentum1D}
    {F}_{gl} = \rho_g^0\,K_I\,|{u}_l-{u}_g|
    \left({u}_l-{u}_g\right)=-{F}_{lg}.
\end{equation}
Concerning the momentum exchange with the wall, which cannot be properly modeled in the one-dimensional case because we have no curvature of the velocity field in axes other than $\mathbf{e}_x$, let us assume $F_g=0$ and 
\begin{equation}\label{wallmomentum1D}
    {F}_{l} = \rho_g^0\,K_W\,|-{u}_l|
    \left(-{u}_l\right).
\end{equation}
It is clear the analogy between Eq. (\ref{interphasemomentum1D}) and Eq. (\ref{wallmomentum1D}), where the latter assumes zero wall velocity. Finally, the dispersed phase volume fraction equation given by Eq. (\ref{voidfraction}), in the one-dimensional case, becomes
\begin{equation}\label{voidfraction1D}
    \frac{\partial \alpha_g}{\partial t} + \frac{\partial}{\partial x} (\alpha_g {u}_g) = 0,
\end{equation}
The artificially compressible continuity equation given by Eq. (\ref{lbm-cont-g1D}), Eq. (\ref{lbm-mom-g1D}), Eq. (\ref{lbm-mom-l1D}) and Eq. (\ref{voidfraction1D}) define a proper system of equations in the one-dimensional case for $p_k$, $u_g$, $u_l$ and $\alpha_g$.

Before proceeding the numerical results, let us search for a special analytical solution indicated by the notation $\overline{\cdot}$. Let us consider a flow regime, where all flow quantities have zero gradients at steady state, but the gradient of the kinematic pressure is constant, namely
\begin{equation}\label{lbm-mom-g1Danalytical}
0 = - \partial_x \overline{p}_k + (R-1)\,{g} + \frac{K_I}{\overline{\alpha}_g}\,|\overline{{u}}_l-\overline{{u}}_g|
    \left(\overline{{u}}_l-\overline{{u}}_g\right),
\end{equation}
\begin{equation}\label{lbm-mom-l1Danalytical}
0 = - \frac{1}{R}\,\partial_x \overline{p}_k + \frac{1}{\overline{\alpha}_l R}\left[K_I\,|\overline{{u}}_g-\overline{{u}}_l|
    \left(\overline{{u}}_g-\overline{{u}}_l\right)-K_W\,|\overline{{u}}_l|
    \overline{{u}}_l\right].
\end{equation}
Eliminating the kinematic pressure gradient from the previous equations yields
\begin{equation}\label{analytical}
\frac{K_I}{\overline{\alpha}_g\,\overline{\alpha}_l}\,
|\overline{{u}}_l-\overline{{u}}_g|
    \left(\overline{{u}}_l-\overline{{u}}_g\right)
    +\frac{K_W}{\overline{\alpha}_l}\,|\overline{{u}}_l|
    \overline{{u}}_l+(R-1)\,{g}=0,
\end{equation}
which can be interpreted as a condition for $\overline{{u}}_l$ as a function of $\overline{{u}}_g$ in this flow regime. In particular, let us search for $\overline{{u}}_g^0$ such that $\overline{{u}}_l^0=0$ (which corresponds to the bubble column case), namely
\begin{equation}\label{analytical0}
(R-1)\,{g}=\frac{K_I}{\overline{\alpha}_g\,\overline{\alpha}_l}\,
|\overline{{u}}_g^0|
    \overline{{u}}_g^0.
\end{equation}
Substituting the previous expression back into the original equation yields
\begin{equation}\label{analytical2}
\frac{K_W}{\overline{\alpha}_l}\,|\overline{{u}}_l|
    \overline{{u}}_l=
    \frac{K_I}{\overline{\alpha}_g\,\overline{\alpha}_l}\,
|\overline{{u}}_g-\overline{{u}}_l|
    \left(\overline{{u}}_g-\overline{{u}}_l\right)
    -\frac{K_I}{\overline{\alpha}_g\,\overline{\alpha}_l}\,
|\overline{{u}}_g^0|
    \overline{{u}}_g^0.
\end{equation}
Let us assume $\overline{{u}}_l<\overline{{u}}_g$ (which corresponds to the natural circulation loops driven by gas injection), then the previous condition becomes
\begin{equation}\label{analytical3}
r \left(\overline{{u}}_l\right)^2=
    \left(\overline{{u}}_g-\overline{{u}}_l\right)^2
    -(\overline{{u}}_g^0)^2,
\end{equation}
where $r=(K_W/K_I)\,\overline{\alpha}_g$, which admits the following analytical solution
\begin{equation}\label{analyticalsolution}
\overline{{u}}_l = 
\frac{\overline{{u}}_g-\sqrt{(\overline{{u}}_g)^2
-(1-r)[(\overline{{u}}_g)^2-(\overline{{u}}_g^0)^2]}}
{1-r}.
\end{equation}
It is worth discussing a couple of limiting cases. If $r=0$, then the only possible coupling between the two phases is by the pressure gradient $\partial_x \overline{p}_k$: the relation becomes $\overline{{u}}_l=\overline{{u}}_g-\overline{{u}}_g^0$. If $r=1$, then the coupling between the phases is ruled by the same proportionality than that between the liquid and the wall: in this case, the relation becomes $\overline{{u}}_l=(\overline{{u}}_g-\overline{{u}}_g^0)\,(1+\overline{{u}}_g^0/\overline{{u}}_g)/2$.

\subsection{Very large density ratios}\label{largedensityratios}

In the limit of very large density ratios, namely $R\gg 1$, some terms proportional to $1/R$ in the momentum equation of the liquid phase become very small and comparable with numerical errors, leading to numerical instabilities. Hence, special numerical ingredients are needed to overcome these instabilities. In this section, we analyze two of such ingredients.

In Ref. \cite{asinariOhwada}, Ohwada suggested a simple but effective ingredient for the suppression of spurious acoustic mode, in the context of single-phase artificial compressibility method. The basic idea here is the introduction of a similar dissipation term into the continuity equation for the gas, given by Eq. (\ref{lbm-cont-g1D}), but we also automatically ensure the right asymptotic target equation in the context of the LBM. Using again the usual Boltzmann scaling adopted in the LBM numerical codes, this idea reads
\begin{equation}\label{lbm-cont-g1Dhat}
    \frac{\partial \hat{\epsilon}_g}{\partial \hat{t}}+
    \frac{\partial \hat{u}_g}{\partial \hat{x}} = \hat{S}_{g}-\gamma\,
    \hat{u}_g^2\,
    \left(\hat{\epsilon}_g-1\right)\equiv \hat{S}_{g}',
\end{equation}
where $\gamma$ is a tunable constant. From a coding point of view, the previous ingredient is equivalent to using the source term $\hat{S}_{g}'$ instead of $\hat{S}_{g}$. The previous equation can be reformulated as
\begin{equation}\label{lbm-cont-g1Dhat2}
    \left(\frac{\partial }{\partial \hat{t}}+
    \gamma\,
    \hat{u}_g^2\right)
    \left(\hat{\epsilon}_g-1\right)+
    \frac{\partial \hat{u}_g}{\partial \hat{x}} = \hat{S}_{g},
\end{equation}
or equivalently
\begin{equation}\label{lbm-cont-g1D-ohwada}
    \frac{\varepsilon^2}{c_s^2}\left(\frac{\partial p_k}{\partial t}+
    \gamma\,{u}_g^2\,p_k
    \right)
    +\frac{\partial {u}_g}{\partial x} = S_{g},
\end{equation}
where it is clear that the term multiplied by $\gamma$ acts as a dashpot in a simple mechanical oscillation system. Concerning the continuity equation for the liquid phase, some options are possible, which can be discussed by the following generic expression:
\begin{equation}\label{lbm-cont-l1Dhat}
    \frac{\partial \hat{\epsilon}_l}{\partial \hat{t}}+
    \frac{\partial \hat{u}_l}{\partial \hat{x}} = 
    \hat{S}_{l}-\gamma\,
    \hat{u}_l^2\,
    \hat{\Gamma}_l\equiv \hat{S}_{l}',
\end{equation}
where $\hat{\Gamma}_l=\epsilon^2\,{\Gamma}_l$ is a function to be specified, or equivalently 
\begin{equation}\label{lbm-cont-l1D-ohwada}
    \varepsilon^2\left(\frac{\partial {\epsilon}_l}{\partial {t}}
    +\gamma\,
    {u}_l^2\,
    \Gamma_l\right)
    +\frac{\partial {u}_l}{\partial {x}} = 
    {S}_{l}.
\end{equation}
Different strategies are possible in choosing the function $\Gamma_l=\Gamma_l(p_k)$. 
\begin{itemize}
    \item The strategy most similar to the one adopted for the gas phase implies $\hat{\Gamma}_l=\hat{\epsilon}_l-1$ or equivalently $\Gamma_l={\epsilon}_l$ because of Eq. (\ref{pressurescaling-l}), which would ensure a proper dashpot for the liquid phase as well. The problem is that ${\epsilon}_l$ is not simply related to $p_k$ and therefore, for the finite-difference engine (considered for comparison), this strategy would require solving an additional equation (precisely for ${\epsilon}_l$).
    \item The second strategy is based on assuming that $p_k$ is finally what really matters for both phases and therefore it assumes $\hat{\Gamma}_l=\hat{\epsilon}_g-1$ for the liquid phase as well, which implies $\Gamma_l(p_k)={\epsilon}_g=p_k/c_s^2$ because of Eq. (\ref{pressurescaling-g}). This strategy has definitively the advantage of simplicity, but there is the drawback that different quantities are involved in the first two terms of Eq. (\ref{lbm-cont-l1D-ohwada}).    
    \item The third strategy aims to improve the consistency by using liquid phase quantities in computing $\hat{\Gamma}_l$ by including also $\hat{\phi}$, namely $\hat{\Gamma}_l=\hat{\phi}\,\hat{\epsilon}_l-1$, which implies $\Gamma_l(p_k)={\epsilon}_g/R=(p_k/c_s^2)\,/R$ because of Eq. (\ref{pressurescaling-phil}). The advantage is that $p_k/R$ is exactly the effective pressure acting on the liquid phase, which makes sense to use in the design of the liquid dashpot as well. On the other hand, the first two terms of Eq. (\ref{lbm-cont-l1D-ohwada}) are now formulated for ${\epsilon}_l$ and $p_k$ respectively, which is different from what happens in canonical modeling of a mechanical dashpot.  
\end{itemize}

As far as LBM is concerned, these three strategies are substantially equivalent. In this work, for the above reasons, we choose to adopt the third strategy.

The second ingredient for dealing with very large density ratios consists in updating the sources of the continuity equations only once on a while to let the fluid equations gradually accommodate the changes. In particular, it consists in updating them only at times $\hat{t}_\gamma$ which are multiple of the natural number $n_\gamma$, namely
\begin{equation}\label{lbm-cont-1D-update}
\hat{S}_\varphi''\left(\hat{{x}},
\hat{t}\right)=\hat{S}_\varphi'\left(\hat{{x}},
\hat{t}_\gamma\right),
\end{equation}
where $\hat{t}_\gamma$ is the highest multiple of $n_\gamma$ but still smaller than $\hat{t}$, namely $\hat{t}_\gamma = r_\gamma\,n_\gamma$ where $r_\gamma$ is another natural number and $r_\gamma\,n_\gamma\leq\hat{t}<(r_\gamma+1)\,n_\gamma$. Of course $n_\gamma$ is a free tunable parameter such that $n_\gamma\ll N_t$, where $N_t$ is the total number of time steps. It is worth the effort to estimate the error introduced by the approximation given by Eq. (\ref{lbm-cont-1D-update}), which means to estimate the difference $|\hat{S}_\varphi'(\hat{{x}},
\hat{t})-\hat{S}_\varphi'(\hat{{x}},
\hat{t}_\gamma)|$ or equivalently
\begin{equation}\label{lbm-cont-1D-update-error}
\left|\hat{S}_\varphi'\left(\hat{{x}},
\hat{t}\right)-\hat{S}_\varphi''\left(\hat{{x}},
\hat{t}\right)\right|\approx
(\hat{t}-\hat{t}_\gamma)
\left|\frac{\partial \hat{S}_\varphi'}{\partial \hat{t}}\right|<n_\gamma\left|\frac{\partial \hat{S}_\varphi'}{\partial \hat{t}}\right|
=O(\epsilon^4).
\end{equation}
This proves that the approximation given by Eq. (\ref{lbm-cont-1D-update}) does not spoil the order of convergence of the proposed methodology. Even though this ingredient is simple and computationally very cheap, it ensures excellent stability, but it smooths out the fastest dynamics (which is not compatible anyway with the incompressible limit).

\subsection{Realistic phenomenological relation for the drag force}

In order to better analyze the proposed methodology in practical applications, let us consider a more realistic phenomenological relation for the drag force. In particular, let us consider the model by Clift, Grace and Weber (1978) \cite{PFLEGER19995091}. In this model, the effective drag coefficient $K_I$ in Eq. (\ref{interphasemomentum1D}) of the interphase momentum exchange can be expressed as $K_I=\kappa_I\,\alpha^2\,\rho^0/\rho_g^0$, where $\kappa_I=(3\,C_d)/(4\,d)$, $C_d$ is the drag coefficient which depends on the relative Reynolds number, $d$ is the characteristic size of the individual entity of the dispersed phase, $\alpha=\sqrt{\alpha_g\,\alpha_l}$ and $\rho^0=\alpha_g\rho_g^0+\alpha_l\rho_l^0$. Combining these definitions yields
\begin{equation}\label{CliftGraceWeber}
K_I = \kappa_I \left(\alpha_g\,\alpha_l\right)
\left(\alpha_g+\alpha_l\,R\right).
\end{equation}
The model by Clift, Grace and Weber (CGW) introduces a further dependence on the volume fraction into the effective drag coefficient. This dependence can be analyzed by setting $K_I = \kappa_I\,\Lambda$ where 
\begin{equation}\label{Lambda}
\Lambda(\alpha_g) \equiv 
\left(\alpha_g\,\alpha_l\right)
\left(\alpha_g+\alpha_l\,R\right).
\end{equation}
The function $\Lambda(\alpha_g)$ goes to zero for both $\alpha_g=0$ and $\alpha_g=1$, with a positive maximum in between. In the limiting case that the gas is the dispersed phase (namely that $\alpha_g$ is small enough), the function $\Lambda(\alpha_g)$ can be approximated by its tangent at $\alpha_g=0$, namely $\alpha^2\approx \alpha_g$ and $\rho^0\approx\rho_l^0$, which imply $\Lambda \approx \Lambda_0 = \alpha_g \,R$ \cite{PFLEGER19995091} \cite{Maniscalco2021}. In case of large density ratios, the tangent $\Lambda_0 = \alpha_g \,R$ becomes quite steep and this makes the CGW model stiff because it can change significantly the drag force for moderate changes of the void fraction. For the purpose of containing the numerical instability due to this stiffness, let us update the function $\Lambda(\alpha_g)$ only once on a while to let the equations gradually accommodate the changes
\begin{equation}\label{Lambda-update}
\Lambda'\left(\hat{{x}},
\hat{t}\right)=\Lambda\left(\hat{{x}},
\hat{t}_\gamma\right),
\end{equation}
where again $\hat{t}_\gamma$ is the highest multiple of $n_\gamma$ but still smaller than $\hat{t}$, namely $\hat{t}_\gamma = r_\gamma\,n_\gamma$ where $r_\gamma$ is a natural number and $r_\gamma\,n_\gamma\leq\hat{t}<(r_\gamma+1)\,n_\gamma$. Again $n_\gamma\ll N_t$, where $N_t$ is the total number of time steps. It is possible to prove (see previous section) that this approximation also does not spoil the order of convergence of the proposed methodology.

Another important ingredient for realistic simulations by LBM is the concept of fluid-dynamic similarity. Because of its kinetic origin, LBM can not deal directly with macroscopic geometric dimensions. Fortunately, the (kinematic) similarity allows one to apply LBM to a geometrically similar setup (same shape but different sizes) with the same boundary conditions (e.g., no-slip, inlet velocity) and the same relevant dimensionless numbers. In other words, the velocity at any point in the LBM model flow is proportional by a constant scaling factor to the velocity at the same point in the real flow, while maintaining the same flow streamlines. It is important to identify this scaling factor, i.e. a proper way to convert flow quantities in physical units into lattice quantities coded by the Boltzmann scaling. Let us start with the main driving force, namely the buoyancy force, driven by the gravitational acceleration $g$, and its LBM counterpart $\hat{g}$ in Boltzmann scaling, which must be selected to ensure numerical stability. The ratio between these two quantities defines the first constraint in terms of $c/\tau$, namely
\begin{equation}\label{scaling-grav-acceleration}
c/\tau=\frac{g}{\hat{g}}.
\end{equation}
A second constraint is given by a similar ratio between the drag factor in the CGW model, i.e. $\kappa_I$, and its LBM counterpart, i.e. $\hat{\kappa}_I$, which must ensure stable simulations. It is possible to express the last ratio in terms of $c\,\tau$, namely 
\begin{equation}\label{scaling-kappaI}
c\,\tau=\frac{\hat{\kappa}_I}{\kappa_I}.
\end{equation}
The previous constraints allow to compute both $c$ and $\tau$. In particular, the so-called lattice speed, which is the average fictitious particle velocity, can be computed as
\begin{equation}\label{lattice-speed}
c = \sqrt{(c/\tau)\,(c\tau)}=\sqrt{\frac{g}{\hat{g}}\,\frac{\hat{\kappa}_I}{\kappa_I}}.
\end{equation}
This poses an upper limit to the maximum magnitude of the velocity field $u_\varphi(x,t)$ which can be simulated by stable simulations, because the latter relation requires that $\hat{u}_\varphi=u_\varphi/c$ is small enough, coherently with the incompressible limit. If these quantities are chosen consistently with regards to realistic setups, this is usually not problematic. 

In particular, it is worth to compare the two most important driving forces of the momentum equation for the gas phase given by Eq. (\ref{lbm-mom-g1D}), namely the interphase momentum exchange force and the buoyancy force. For the purpose of doing so, the relevant fields, e.g. the velocity field $u_\varphi(x,t)$, can be characterized by a reference value identified by the dagger superscript, e.g. $u_\varphi^{\dagger}$, which is typically the inlet (known) value. Using these reference values, the ratio between the magnitudes of the exchange force and the buoyancy force for the gas phase can be expressed as $M_{ex}$ and, according to the CGW model, it can be computed as
\begin{equation}\label{Mex}
M_{ex} = \frac{\kappa_I\,\Lambda^{\dagger}\left|{u}_g^{\dagger}-{u}_l^{\dagger}\right|^2/\alpha_g^{\dagger}}{(R-1)\,g}.
\end{equation}
where 
$\Lambda^{\dagger}=\Lambda(\alpha_g^{\dagger})=\alpha_g^{\dagger}\,\alpha_l^{\dagger}(\alpha_g^{\dagger}+\alpha_l^{\dagger}\,R)$. Computing the same dimensionless number by LBM quantities in Boltzmann scaling yields
\begin{equation}\label{Mex-hat}
\hat{M}_{ex} = \frac{\hat{\kappa}_I\,\Lambda^{\dagger}\left|\hat{u}_g^{\dagger}-\hat{u}_l^{\dagger}\right|^2/\alpha_g^{\dagger}}{(R-1)\,\hat{g}}=(c/\tau)\,
\frac{c\,\tau}{c^2}\,
M_{ex}=M_{ex},
\end{equation}
which clearly proves that $\hat{M}_{ex}=M_{ex}$ and hence the fluid-dynamic similarity ensures the same ratio between the two main driving forces.

\subsection{Finite difference solver for the reference solution}

In this section, we develop an algorithm based on the finite-difference (FD) method for solving the system of equations given by Eq. (\ref{lbm-cont-g1D}), Eq. (\ref{lbm-mom-g1D}), Eq. (\ref{lbm-mom-l1D}) and Eq. (\ref{voidfraction1D}), in order to obtain numerically a reference solution. The reference solution will be used to validate the LBM results. Hence, in order to simplify the comparisons between the results of the two methods, let us formulate both numerical engines in terms of $(\hat{{x}},\hat{t})$, where $\hat{{x}}$ is the space coordinate divided by the distance $\lambda$ between two neighboring lattice nodes (mean free path) and $\hat{t}$ is the physical time divided by the time $\tau$ between two consecutive lattice collisions (mean collision time). Let us define a computational domain which is a multiple $N_x$ of $\lambda$ and let us search for a solution of the previous system of equations in the nodal points $\hat{x}_i\in\{1,2,\dots, N_x\}$ where $1\leq i\leq N_x$. Let us consider the generic quantity $\hat{u}(\hat{x})$ and let us call the nodal values as $\hat{u}(\hat{x}_i)=\hat{u}_i$. Let us store the nodal values in a vector $\{\hat{u}\}=\{\hat{u}_1,\hat{u}_2,\dots,\hat{u}_{N_x}\}$ with elements $\{\hat{u}\}(i)=\hat{u}_i$ for $1\leq i\leq N_x$. It is possible include also boundary conditions (BCs) into this nomenclature. From the nodal vector of the generic quantity $\{\hat{u}\}$, by wrapping it with proper BCs, let us construct a larger nodal vector $\langle\hat{u}\rangle=\{\hat{u}_0,\hat{u}_1,\hat{u}_2,\dots,\hat{u}_{N_x},\hat{u}_{N_x+1}\}$ with elements $\langle\hat{u}\rangle(i)=\hat{u}_i$ for $0\leq i\leq (N_x+1)$, where $\hat{u}_0$ and $\hat{u}_{N_x+1}$ are some BCs which must be specified for the generic known (as a function of the other nodal values). Let us introduce also the following combinatorial rule. Given two known vectors $\langle\hat{u}\rangle$ and $\langle\hat{v}\rangle$, let us assume that $(\langle\hat{u}\rangle\langle\hat{v}\rangle)(i)=\langle\hat{u}\rangle(i)\langle\hat{v}\rangle(i)$. The same combinatorial rule applies also to vectors $\{\cdot\}$, namely $(\{\hat{u}\}\{\hat{v}\})(i)=\{\hat{u}\}(i)\{\hat{v}\}(i)$.

Moreover, it is possible to introduce the central difference operator $\hat{D}_x\,:\,\langle\hat{u}\rangle\rightarrow\{\hat{D}_x\hat{u}\}$ (contraction) defined as 
\begin{equation}\label{FD-Dx}
(\hat{D}_x\langle\hat{u}\rangle)(i)=
\{\hat{D}_x\hat{u}\}(i)=
\left[\langle\hat{u}\rangle(i+1)-\langle\hat{u}\rangle(i-1)\right]/2,
\end{equation}
where $1\leq i\leq N_x$. Analogously, by applying differencing formulas in a recursive manner, it is possible to introduce $\hat{D}_x^2\,:\,\langle\hat{u}\rangle\rightarrow\{\hat{D}_x^2\hat{u}\}$ (contraction) defined as
\begin{equation}\label{FD-DxDx}
(\hat{D}_x^2\langle\hat{u}\rangle)(i)=
\{\hat{D}_x^2\hat{u}\}(i)=
\langle\hat{u}\rangle(i+1)
-2\,\langle\hat{u}\rangle(i)
+\langle\hat{u}\rangle(i-1),
\end{equation}
where again $1\leq i\leq N_x$. These are the essential ingredients of the FD engine.

The FD engine searches for a numerical solution of the following vectors: $\{\hat{p}_k\}$, $\{\hat{u}_g\}$, $\{\hat{u}_l\}$ and $\{\hat{\alpha}_g\}$. 
 Multiplying Eq. (\ref{lbm-cont-g1D}) by $\varepsilon^2$ and applying the code scaling yields
\begin{equation}\label{lbm-cont-g1DFD}
    \frac{d \{\hat{p}_k\}}{d \hat{t}}    =
    -c_s^2\,\hat{D}_x\langle\hat{u}_g\rangle
    + c_s^2\,\hat{D}_x\left[
    \langle\hat{\alpha}_l\rangle \left(\langle\hat{{u}}_g\rangle - \langle\hat{{u}}_l\rangle\right)\right].
\end{equation}
Multiplying Eq. (\ref{lbm-mom-g1D}) and Eq. (\ref{lbm-mom-l1D}) by $\varepsilon^3$ and applying the code scaling yields
\begin{multline}\label{lbm-mom-g1DFD}
    \frac{d \{\hat{u}_g\}}{d \hat{t}}  = 
    - \hat{D}_x(\langle \hat{u}_g\rangle \langle\hat{u}_g\rangle)
    - \hat{D}_x\langle \hat{p}_k\rangle  + \nu_g\,\hat{D}_x^2\langle \hat{u}_g\rangle + \{\hat{S}_{g}\}\,\{\hat{u}_g\}+ 
    \{\hat{\sigma}_g/\hat{\alpha}_g\}
    \hat{D}_x\langle \hat{\alpha}_g\rangle\dots\\
    +(R-1)\,\hat{g} + \{\hat{K}_I/\hat{\alpha}_g\}\,|\{\hat{u}_l\}-\{\hat{u}_g\}|\left(\{\hat{u}_l\}-\{\hat{u}_g\}\right),
\end{multline}
and
\begin{multline}\label{lbm-mom-l1DFD}
    \frac{d \{\hat{u}_l\}}{d \hat{t}} = 
    - \hat{D}_x(\langle \hat{u}_l\rangle \langle\hat{u}_l\rangle)
    - (1/R)\,\hat{D}_x\langle \hat{p}_k\rangle + \nu_l\,\hat{D}_x^2\langle \hat{u}_l\rangle + 
    \{\hat{S}_{l}\}\,\{\hat{u}_l\} + 
    \{\hat{\sigma}_l/\hat{\alpha}_l\}\, \hat{D}_x\langle\hat{\alpha}_l\rangle\dots\\
    +(1/R)\,\{\hat{K}_I/\hat{\alpha}_l\}\,|\{\hat{{u}}_g\}-\{\hat{{u}}_l\}|
    \left(\{\hat{{u}}_g\}-\{\hat{{u}}_l\}\right)
    +(1/R)\,\{\hat{K}_W/\hat{\alpha}_l\}\,
    |\{\hat{{u}}_l\}|
    \{\hat{{u}}_l\},
\end{multline}
where we assumed $\varepsilon\,K_W=\hat{K}_W$ in analogy with the scaling $\hat{K}_I=\varepsilon K_I$ which was adopted in the asymptotic analysis of the LBM schemes. Finally, multiplying Eq. (\ref{voidfraction1D}) by $\varepsilon^2$ and applying the code scaling yields
\begin{equation}\label{voidfraction1DFD}
    \frac{d \{\hat{\alpha}_g\}}{d \hat{t}} = 
    -\hat{D}_x(\langle \hat{\alpha}_g\rangle 
    \langle\hat{u}_g\rangle).
\end{equation}
Eq. (\ref{lbm-cont-g1DFD}), Eq. (\ref{lbm-mom-g1DFD}), Eq. (\ref{lbm-mom-l1DFD}) and Eq. (\ref{voidfraction1DFD}) define a system of ordinary differential equations (ODEs) for $\{\hat{p}_k\}$, $\{\hat{u}_g\}$, $\{\hat{u}_l\}$ and $\{\hat{\alpha}_g\}$, as far as proper boundary conditions are specified for computing the larger vector $\langle\hat{u}\rangle$ from the generic vector $\{\hat{u}\}$. For the kinematic pressure, this mapping is realized by the following wrapping
\begin{equation}\label{pkBC}
\langle\hat{p}_k\rangle=
\{2\,\hat{p}_k(1)-\hat{p}_k(2),\hat{p}_k(1),\dots,
\hat{p}_k(N_x),
2\,c_s^2-\hat{p}_k(N_x)\},
\end{equation}
where $\hat{p}_k(0)=2\,\hat{p}_k(1)-\hat{p}_k(2)$ at the inlet is an extrapolation and $\hat{p}_k(N_x+1)=2\,c_s^2-\hat{p}_k(N_x)$ at the outlet is a Dirichlet condition corresponding to $\hat{p}_k(N_x+1/2)=c_s^2$ (half-way). Other flow quantities have the following mappings
\begin{equation}\label{ugBC}
\langle\hat{u}_g\rangle=
\{2\,\hat{u}_g^{IN}-\hat{u}_g(1),\hat{u}_g(1),\dots,
\hat{u}_g(N_x),
2\,\hat{u}_g(N_x)-\hat{u}_g(N_x-1)\},
\end{equation}
\begin{equation}\label{ulBC}
\langle\hat{u}_l\rangle=
\{2\,\hat{u}_l^{IN}-\hat{u}_l(1),\hat{u}_l(1),\dots,
\hat{u}_l(N_x),2\,\hat{u}_l(N_x)-\hat{u}_l(N_x-1)\},
\end{equation}
\begin{equation}\label{alphagBC}
\langle\hat{\alpha}_g\rangle=
\{2\,\hat{\alpha}_g^{IN}-\hat{\alpha}_g(1),
\hat{\alpha}_g(1),\dots,
\hat{\alpha}_g(N_x),
2\,\hat{\alpha}_g(N_x)-\hat{\alpha}_g(N_x-1)\},
\end{equation}
where, for the generic quantity, $\hat{u}(0) = 2\,\hat{u}^{IN}-\hat{u}(1)$ at the inlet is a Dirichlet condition corresponding to $\hat{u}(1/2)=\hat{u}^{IN}$ (half-way) and $\hat{u}(N_x+1) = 2\,\hat{u}(N_x)-\hat{u}(N_x-1)$ at the outlet is an extrapolation. In order to be consistent with the incompressible limit, the generic inlet quantity $\hat{y}^{IN}\in\{\hat{\alpha}_g^{IN}, \hat{u}_g^{IN}, \hat{u}_l^{IN}\}$ involved in the previous formulas is progressively increased by the following function
\begin{equation}\label{transient}
\hat{y}^{IN}(\hat{t}) = (\hat{y}^{max}-\hat{y}^{min})\tanh{(\hat{t}/n_t)}
+\hat{y}^{min},
\end{equation}
where $\hat{y}^{max}$ and $\hat{y}^{min}$ are quantity-specific values which depend on the considered test case. The parameter $n_t< N_t$, where $N_t$ is the total number of time steps. 

The above system of ODEs will be solved in Matlab\textregistered\, by means of \textit{ode45} solver, which is a variable step solver (which means that it automatically chooses the value of the time stepping) and is based on an explicit Runge-Kutta (4,5) formula, the Dormand-Prince pair, namely a combination of $4^{th}$ and $5^{th}$ order method. During the iteration procedure, once every few time steps (e.g. 100 iterations), the solution vectors are smoothed in order to avoid the checkerboard instability by means of a Gaussian-weighted moving average filter (with a window containing 12 points). It is remarkable that LBM schemes do not need such additional smoothing because they are more robust against the checkerboard instability. The reason is due to the third term of the left hand side of Eq. (\ref{taylor-zeroth}) in the LBM asymptotic equations \cite{asinariOhwada}.

\subsection{LBM solver}

Before discussing practical details of the LBM solver, it is worth to realize that the one-dimensional case is actually degenerate and this may require some modifications in the formulas involving the transport coefficients. For example, recalling Eq. (\ref{lbm-sigma}) in this case yields
\begin{equation}\label{lbm-sigma1D}
{\sigma}_\varphi=\nu_\varphi^{ef}\,\frac{\partial u_\varphi}{\partial x}
+O(\varepsilon^2),
\end{equation}
where
\begin{equation}\label{lbm-ef-viscosity}
\nu_\varphi^{ef}=\left(\frac{1}{\omega_\varphi}-\frac{1}{2}\right)
-\frac{c_s^2\psi_\varphi }{\omega_\varphi}.
\end{equation}
Using Eq. (\ref{effective-viscosity}) and Eq. (\ref{effective-bulk}) leads to an effective viscosity $\nu_\varphi^{ef}=\nu_\varphi\,d_\varphi$ where $d_\varphi$ is a corrective factor for one-dimensional degeneracy, namely
\begin{equation}\label{1D-correction}
d_\varphi = \frac{4}{3}+\frac{\xi_\varphi}{\nu_\varphi}.
\end{equation}
Hence two strategies are possible for recovering the given kinematic viscosity:
\begin{itemize}
\item $\xi_\varphi=-(1/3)\,\nu_\varphi$, as reported in the previous section, and consequently $d_\varphi=1$ (with $\psi_\varphi=2\,{\omega_\varphi}{\nu_\varphi}/{c_s^2}$), which ensures the maximum consistency with the multidimensional cases but also some small oscillations during the transient dynamics (up to 10\% in the tested cases) without impacting on the steady state solution;
\item $\xi_\varphi=(5/3)\,\nu_\varphi$ and consequently $d_\varphi=3$ (with $\psi_\varphi=0$) which ensures the best performance also during the transient dynamics. 
\end{itemize}
In the following, the second strategy is adopted. This has an impact also on the way one computes the stress in the one-dimensional case. Recalling Eq. (\ref{taylor-second3}) in this case yields
\begin{equation}\label{taylor-second31D}
\frac{\partial u_\varphi}{\partial x}
= \omega_\varphi({{\Pi}}^{eq}_\varphi-{{\Pi}}_\varphi)+\psi_\varphi c_s^2\,S_\varphi+O(\varepsilon^2),
\end{equation}
and substituting it into Eq. (\ref{lbm-sigma1D}) yields
\begin{equation}\label{taylor-second31Db}
\sigma_\varphi
= \nu_\varphi^{ef}\left[\omega_\varphi({{\Pi}}^{eq}_\varphi-{{\Pi}}_\varphi)+\psi_\varphi c_s^2\,S_\varphi
\right]+O(\varepsilon^2),
\end{equation}
where the effective viscosity is used. 

\begin{table}
 \centering
  \begin{tabular}{ccc|c}
  \multicolumn{3}{c|}{INLET ($\hat{x}=1$): compute $\hat{y}^{B(in)}$} & $f_{*\varphi}$ \\
  \hline\hline
  Type & generic density $\hat{y}^{B(in)}_{\epsilon}$ & generic velocity $\hat{y}^{B(in)}_{u}$ &  \\
  \hline
  BB(I) & $\hat{\epsilon}_g$: EX(hw) & $\hat{u}_g=\hat{u}_g^{IN}$  & $f_g$ \\
  BB(I) & $\hat{\epsilon}_l$: EX(hw) & $\hat{u}_l=\hat{u}_l^{IN}$ & $f_l$ \\
  EQ(S) & $\hat{\alpha}_g=\hat{\alpha}_g^{IN}$ & $\hat{u}_g=\hat{u}_g^{IN}$ & $f_{\alpha g}$ \\
  EQ(S) & $\hat{\alpha}_l=1-\hat{\alpha}_g^{IN}$ & $\hat{u}_l=\hat{u}_l^{IN}$ & $f_{\alpha l}$ \\
  EQ(S) & $\hat{\beta}_g$: EX(fw) & $\hat{c}_g$: EX(fw) & $f_{\beta g}$ \\
  EQ(S) & $\hat{\beta}_l$: EX(fw) & $\hat{c}_l$: EX(fw) & $f_{\beta l}$ \\
  \end{tabular}
  \caption{BCs for LBM schemes: INLET ($\hat{x}=1$). Acronyms stand for: bounce back (BB); anti bounce back (ABB); extrapolation (EX) and equilibrium (EQ). Uppercase letter in parentheses stand for: incompressible equilibrium (I) and standard equilibrium (S). Lowercase letter in parentheses stand for: half-way (hw) and full-way (fw).}
  \label{tab:BCforLBMINLET}
\end{table}

\begin{table}
 \centering
  \begin{tabular}{c|ccc}
  $f_{*\varphi}$ & \multicolumn{3}{|c}{OUTLET ($\hat{x}=N_x$): compute $\hat{y}^{B(out)}$} \\
  \hline\hline
  & Type &  generic density $\hat{y}^{B(out)}_{\epsilon}$ & generic velocity $\hat{y}^{B(out)}_{u}$ \\
  \hline
  $f_g$ & ABB(I) & $\hat{\epsilon}_g=1$ & $\hat{u}_g$: EX(hw) \\
  $f_l$ & ABB(I) & $\hat{\epsilon}_l=1$ & $\hat{u}_l$: EX(hw) \\
  $f_{\alpha g}$ & ABB(S) & $\hat{\alpha}_g$: EX(hw) & $\hat{u}_g$: EX(hw) \\
  $f_{\alpha l}$ & ABB(S) & $\hat{\alpha}_l$: EX(hw) & $\hat{u}_l$: EX(hw) \\
  $f_{\beta g}$ & EQ(S) & $\hat{\beta}_g$: EX(fw) & $\hat{c}_g$: EX(fw) \\
  $f_{\beta l}$ & EQ(S) & $\hat{\beta}_l$: EX(fw) & $\hat{c}_l$: EX(fw) \\
  \end{tabular}
  \caption{BCs for LBM schemes: OUTLET ($\hat{x}=N_x$). Acronyms stand for: bounce back (BB); anti bounce back (ABB); extrapolation (EX) and equilibrium (EQ). Uppercase letter in parentheses stand for: incompressible equilibrium (I) and standard equilibrium (S). Lowercase letter in parentheses stand for: half-way (hw) and full-way (fw).}
  \label{tab:BCforLBMOUTLET}
\end{table}

For solving the one-dimensional test case, let us adopt the standard implementation for the D1Q3 lattice \cite{KruegerLBM}, with the equilibrium distributions reported in \ref{D1Q3}. Aside from the bulk implementation discussed so far, BCs in LBM are quite critical for stability and efficiency. The proposed BCs are those which showed the best compromise between stability and accuracy. They can be formulated by the two steps described below.
\begin{itemize}
    \item STEP \#1. First of all, the generic flow quantity $\hat{y}\in\{\hat{\epsilon}_g, \hat{\epsilon}_l, \hat{\alpha}_g, \hat{\alpha}_l, \hat{\beta}_g, \hat{\beta}_l, \hat{u}_g, \hat{u}_l, \hat{c}_g, \hat{c}_l\}$ at the boundary $\hat{x}_B$ (boundary) can be imposed (Dirichlet) or extrapolated. Concerning extrapolation, there are two cases: half-way (hw), where the computed generic quantity is intended at the position $\hat{x}_B$ mid-way between the last mesh node at $\hat{x}_E$ (edge) and the missing node outside the mesh, namely $|\hat{x}_B-\hat{x}_E|=1/2$, and full-way (fw), where the computed quantity is intended in the missing node, namely $|\hat{x}_B-\hat{x}_E|=1$. These two cases can be summarized as
\begin{equation}\label{boundary}
|\hat{x}_B-\hat{x}_E| = 
    \begin{cases}
            1/2, & \text{if half-way (hw)},\\
            1, & \text{if full-way (fw)}.\\
    \end{cases}
\end{equation}
    The choice between these two options depends on the LBM boundary condition (see the following step). At the position $\hat{x}_B$, the generic flow quantity can be imposed or extrapolated
\begin{equation}\label{BCIN}
\hat{y}^{B(in)} = 
    \begin{cases}
            \hat{y}^{IN}, & \text{if Dirichlet},\\
            (3/2)\,\hat{y}(1)-(1/2)\,\hat{y}(2), & \text{if EX(hw)},\\
            2\,\hat{y}(1)-\hat{y}(2), & \text{if EX(fw)},\\
    \end{cases}
\end{equation}
and
\begin{equation}\label{BCOUT}
\hat{y}^{B(out)} = 
    \begin{cases}
            \hat{y}^{OUT}, & \text{if Dirichlet},\\
            (3/2)\,\hat{y}(N_x)-(1/2)\,\hat{y}(N_x-1), & \text{if EX(hw)},\\
            2\,\hat{y}(N_x)-\hat{y}(N_x-1), & \text{if EX(fw)}.\\
    \end{cases}
\end{equation}    
    The proposed BCs for the LBM framework are described in Table \ref{tab:BCforLBMINLET} for the inlet at $\hat{x}=1$ and in Table \ref{tab:BCforLBMOUTLET} for the outlet at $\hat{x}=N_x$, respectively. For the sake of the following step, the generic flow quantities can be divided in two subgroups: generic densities $\hat{y}_\epsilon\in\{\hat{\epsilon}_g, \hat{\epsilon}_l, \hat{\alpha}_g, \hat{\alpha}_l, \hat{\beta}_g, \hat{\beta}_l\}$ and generic velocities $\hat{y}_u\in\{ \hat{u}_g, \hat{u}_l, \hat{c}_g, \hat{c}_l\}$. All previous options are summarized in Table \ref{tab:BCforLBMINLET} for the inlet and in Table \ref{tab:BCforLBMOUTLET} for the outlet by using the same acronyms (see columns for generic density and for generic velocity). 
    \item STEP \#2. Secondly, once the flow quantities are imposed or extrapolated, the usual LBM techniques can be used to transfer the boundary condition at the position $\hat{x}_B$ to the generic distribution function $f_{*\varphi}\in\{f_g,f_l,f_{\alpha g},f_{\alpha l},f_{\beta g},f_{\beta l}\}$ in the mesh node $\hat{x}_E$ during streaming. Let us consider some popular approaches: bounce back (BB) rule, anti bounce back (ABB) rule and equilibrium (EQ) rule \cite{KruegerLBM}. It is worth noticing that also BB and ABB are formulated in terms of the equilibrium distribution, which must be coherent with the one used in the bulk, namely incompressible equilibrium (I) $f_{I}^{eq}(\phi,\hat{\epsilon},\hat{{u}})$ or standard equilibrium (S) $f_S^{eq}(\hat{\alpha},\hat{{u}})$ (see \ref{D1Q3}). Let us refer to the generic equilibrium distribution function as $f_{Y}^{eq} (\hat{y}_\epsilon,\hat{y}_u)$, which can coincide with either $f_I^{eq}$ or $f_{S}^{eq}$ depending on the considered distribution function. The BB, ABB and EQ rule for computing the incoming distribution function during the streaming can be expressed as
\begin{equation}\label{BCINLET}
f_{*\varphi}(\hat{x}_E,
\hat{t}+1,BB(q^*))=
    \begin{cases}
            f^*_{*\varphi}(\hat{x}_E,
            \hat{t},q^*)+\Delta_Y^{BB}(\hat{y}_\epsilon^{B},\hat{y}_u^{B},q^*), & \text{if BB(Y)},\\
            -f^*_{*\varphi}(\hat{x}_E,
            \hat{t},q^*)+\Delta_Y^{ABB}(\hat{y}_\epsilon^{B},\hat{y}_u^{B},q^*), & \text{if ABB(Y)},\\
            f_{Y}^{eq} (\hat{y}_\epsilon^{B},\hat{y}_u^{B},BB(q^*)), & \text{if EQ(Y)},\\
    \end{cases}
\end{equation}
    where $\hat{x}_E$ is the edge node (in the one-dimensional case, it can be either $\hat{x}_E=1$ or $\hat{x}_E=N_x$), 
    $q^*$ is the identifier of the velocity $v_{q^*}$ leaving the computational domain ($v_{q^*}=-1$ for $\hat{x}_E=1$ and $v_{q^*}=+1$ for $\hat{x}_E=N_x$), $BB(q^*)$ is the opposite direction (bounce-back) and $f^*_{*\varphi}$ is the generic post-collision distribution function \cite{KruegerLBM}. The operators $\Delta_Y^{BB}$ and $\Delta_Y^{ABB}$ can be constructed by means of the proper equilibrium for the considered distribution function as
\begin{equation}\label{BB}
\Delta_Y^{BB} = 
f_{Y}^{eq} (\hat{y}_\epsilon^{B},-\hat{y}_u^{B})
-f_{Y}^{eq} (\hat{y}_\epsilon^{B},\hat{y}_u^{B}),
\end{equation}
\begin{equation}\label{ABB}
\Delta_Y^{ABB} = 
f_{Y}^{eq} (\hat{y}_\epsilon^{B},-\hat{y}_u^{B})
+f_{Y}^{eq} (\hat{y}_\epsilon^{B},\hat{y}_u^{B}).
\end{equation}
    The rationale behind the reported choices is the following: BB rule is typically used to transfer a Dirichlet condition for a generic velocity $\hat{y}_u$, ABB rule to transfer a Dirichlet condition for a generic density $\hat{y}_\epsilon$ (and hence also for a generic pressure), while EQ rule to transfer a Dirichlet condition for both a generic velocity and a generic density at the same time. The BB and the ABB rule imposes conditions half-way, namely at position $\hat{x}_B$ where $|\hat{x}_B-\hat{x}_E|=1/2$, while EQ rule streams the distribution function from a missing node located at full-way distance, namely from the position $\hat{x}_B$ where $|\hat{x}_B-\hat{x}_E|=1$. See Eq. (\ref{boundary}), Eq. (\ref{BCIN}) and Eq. (\ref{BCOUT}) from which we started. 
\end{itemize}

\begin{figure}
\begin{subfigure}{.5\textwidth}
  \centering
  \includegraphics[width=1.0\linewidth]{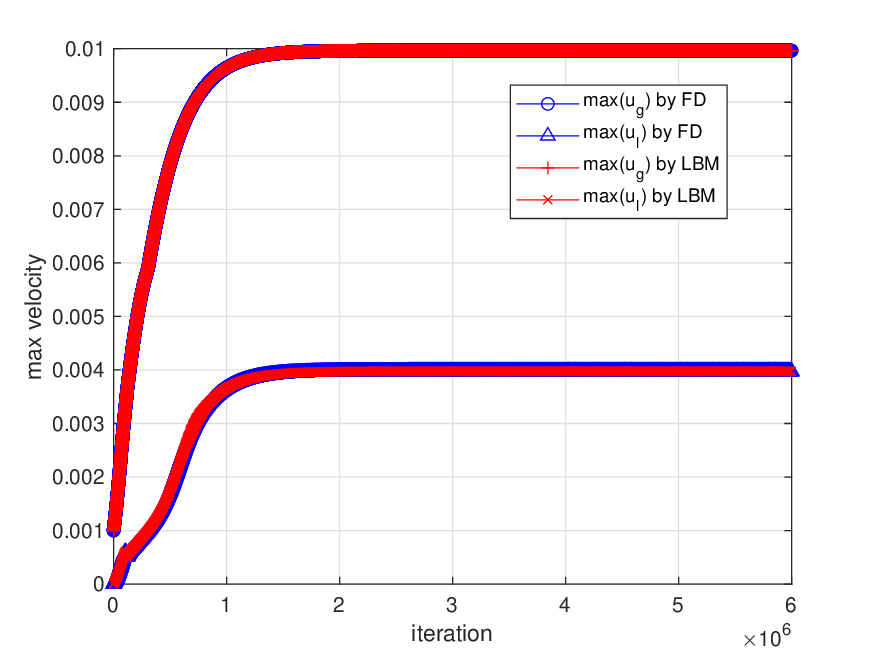}
  \caption{FD vs. LBM: max velocity}
  \label{fig1:sfig1}
\end{subfigure}
\begin{subfigure}{.5\textwidth}
  \centering
  \includegraphics[width=1.0\linewidth]{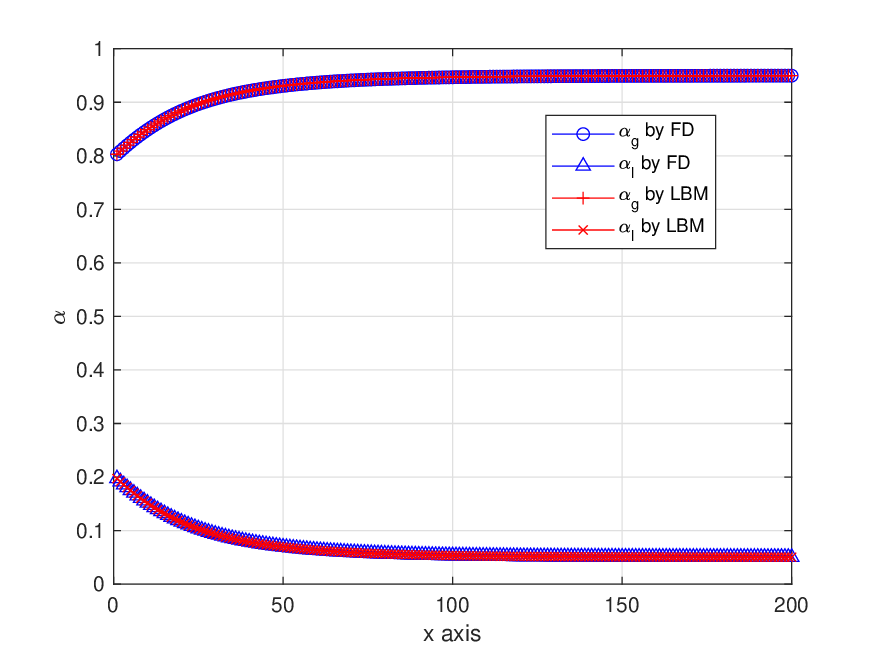}
  \caption{FD vs. LBM: volume fractions}
  \label{fig1:sfig2}
\end{subfigure}
\begin{subfigure}{.5\textwidth}
  \centering
  \includegraphics[width=1.0\linewidth]{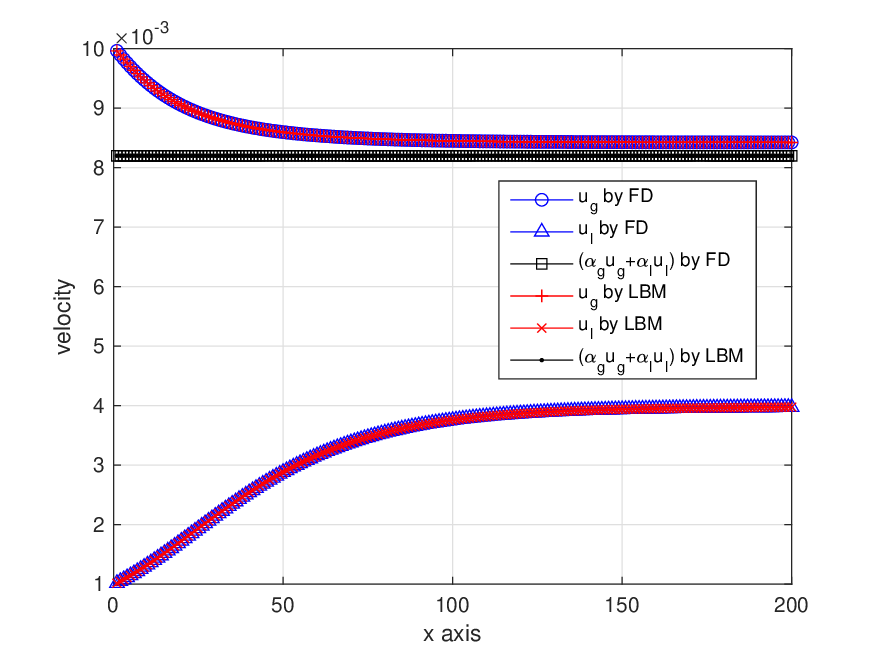}
  \caption{FD vs. LBM: phase velocities}
  \label{fig1:sfig3}
\end{subfigure}
\begin{subfigure}{.5\textwidth}
  \centering
  \includegraphics[width=1.0\linewidth]{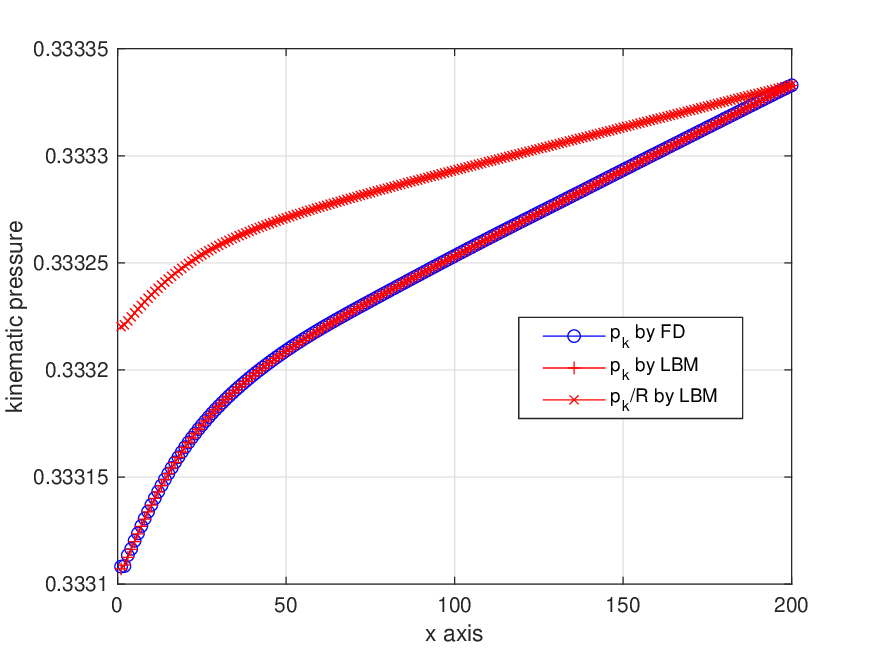}
  \caption{FD vs. LBM: kinematic pressure}
  \label{fig1:sfig4}
\end{subfigure}
\caption{TEST \#1: Comparison between numerical results by FD and by LBM ($R=2$ and $\hat{g} = 10^{-6}$).}
\label{fig:fig1}
\end{figure}

\begin{figure}
\begin{subfigure}{.5\textwidth}
  \centering
  \includegraphics[width=1.0\linewidth]{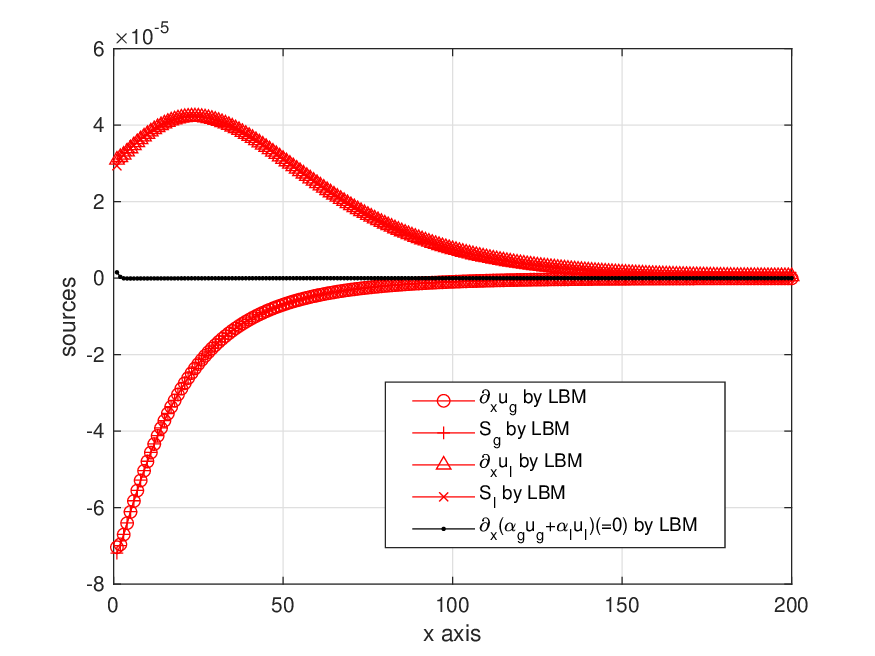}
  \caption{Phase sources computed by LBM}
  \label{fig2:sfig1}
\end{subfigure}
\begin{subfigure}{.5\textwidth}
  \centering
  \includegraphics[width=1.0\linewidth]{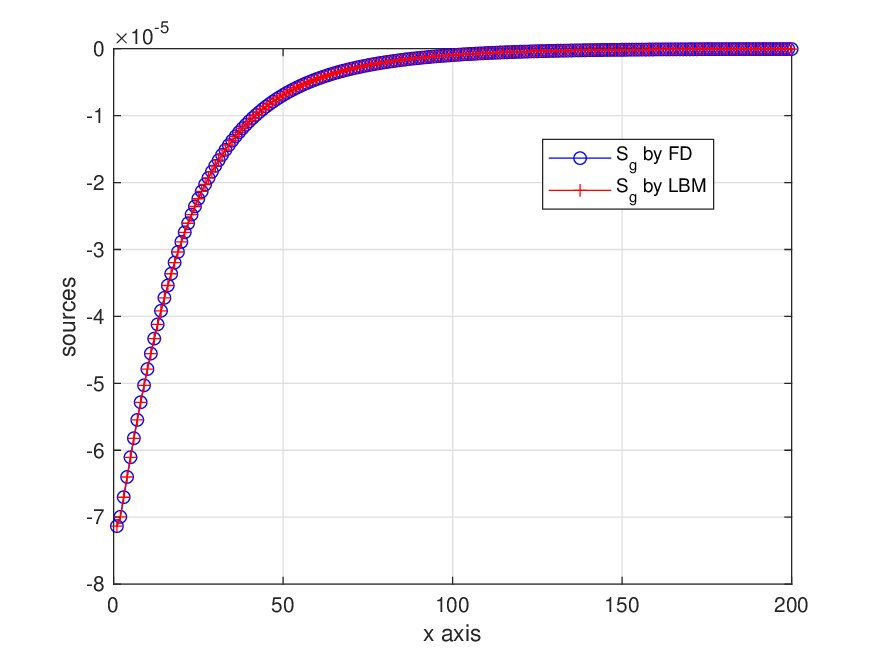}
  \caption{FD vs. LBM: dispersed phase source}
  \label{fig2:sfig2}
\end{subfigure}
\begin{subfigure}{.5\textwidth}
  \centering
  \includegraphics[width=1.0\linewidth]{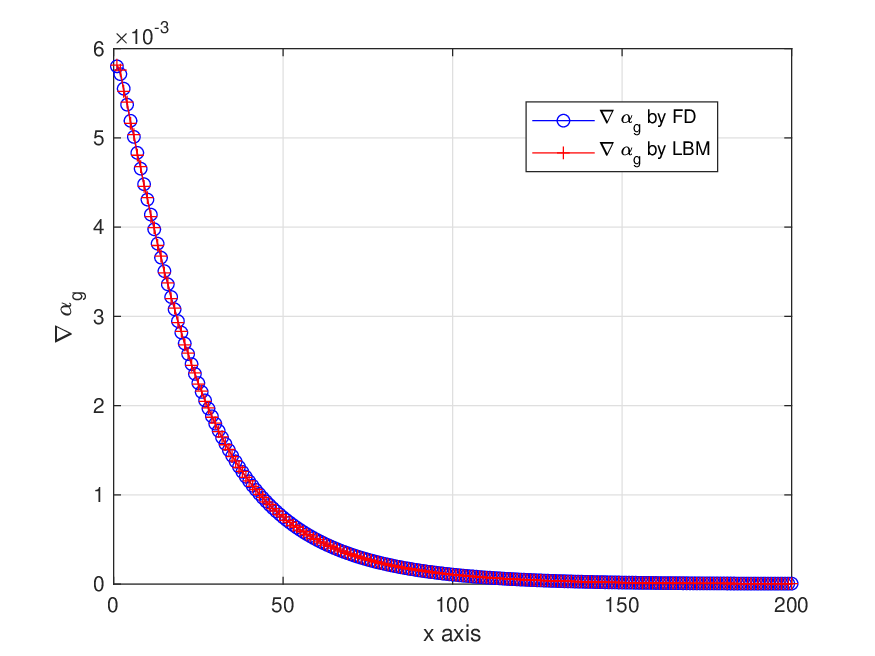}
  \caption{FD vs. LBM: dispersed phase gradient}
  \label{fig2:sfig3}
\end{subfigure}
\begin{subfigure}{.5\textwidth}
  \centering
  \includegraphics[width=1.0\linewidth]{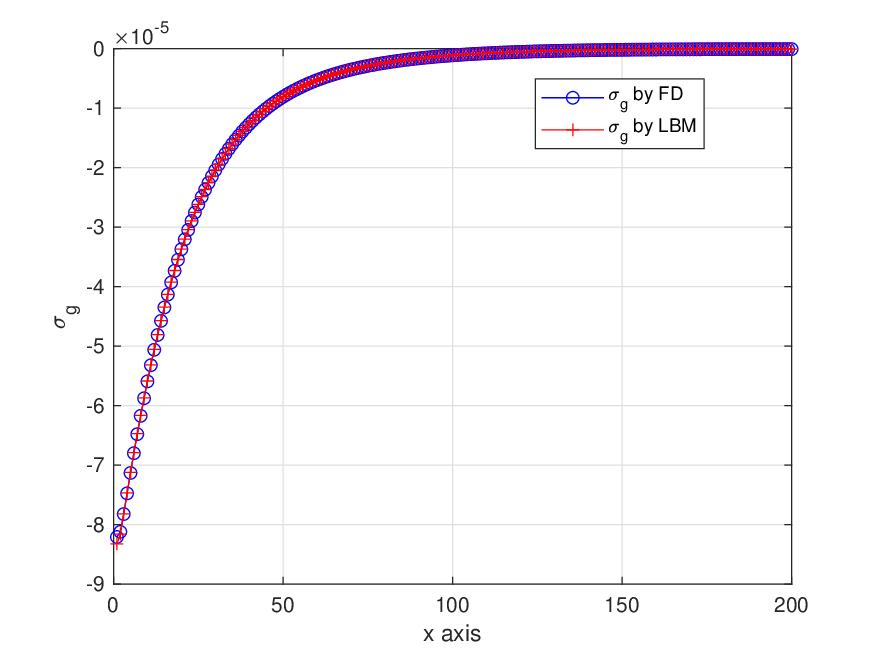}
  \caption{FD vs. LBM: dispersed phase stress}
  \label{fig2:sfig4}
\end{subfigure}
\begin{subfigure}{.5\textwidth}
  \centering
  \includegraphics[width=1.0\linewidth]{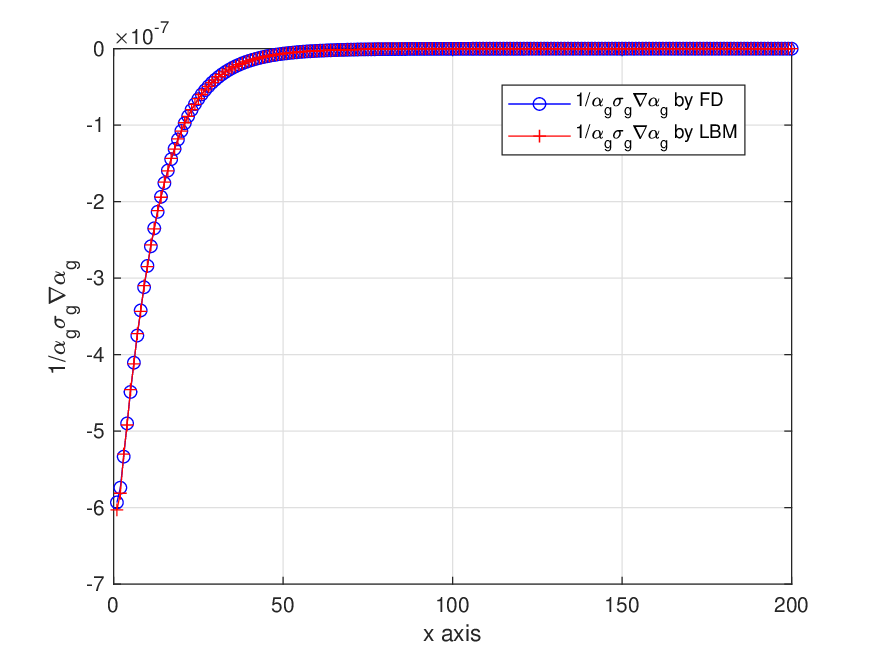}
  \caption{FD vs. LBM: dispersed phase force}
  \label{fig2:sfig5}
\end{subfigure}
\begin{subfigure}{.5\textwidth}
  \centering
  \includegraphics[width=1.0\linewidth]{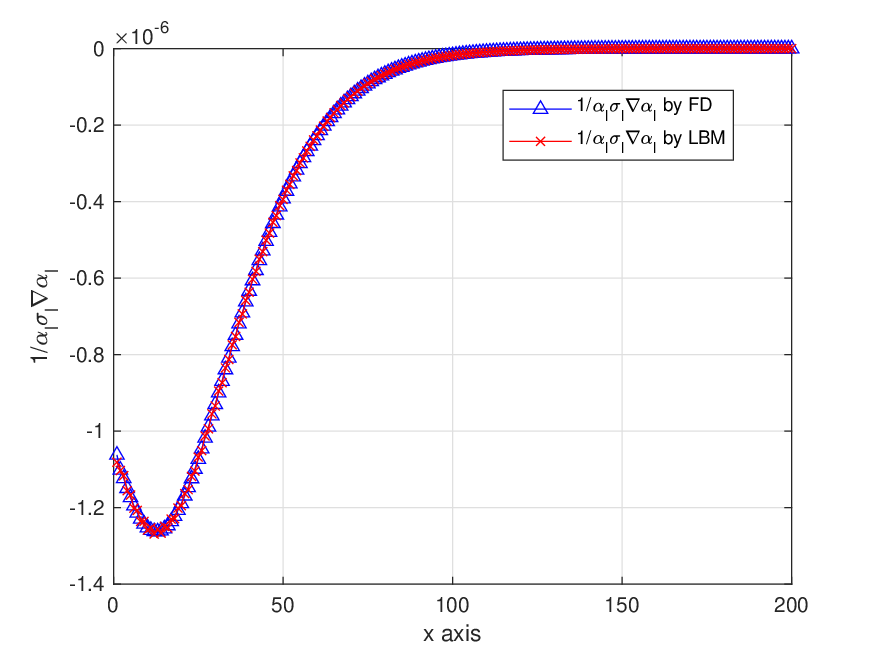}
  \caption{FD vs. LBM: liquid phase force}
  \label{fig2:sfig6}
\end{subfigure}
\caption{TEST \#1: Computational details of the LBM schemes without FD corrections and comparison with numerical results by FD operators ($R=2$ and $\hat{g} = 10^{-6}$).}
\label{fig:fig2}
\end{figure}

\begin{figure}
\begin{subfigure}{.5\textwidth}
  \centering
  \includegraphics[width=1.0\linewidth]{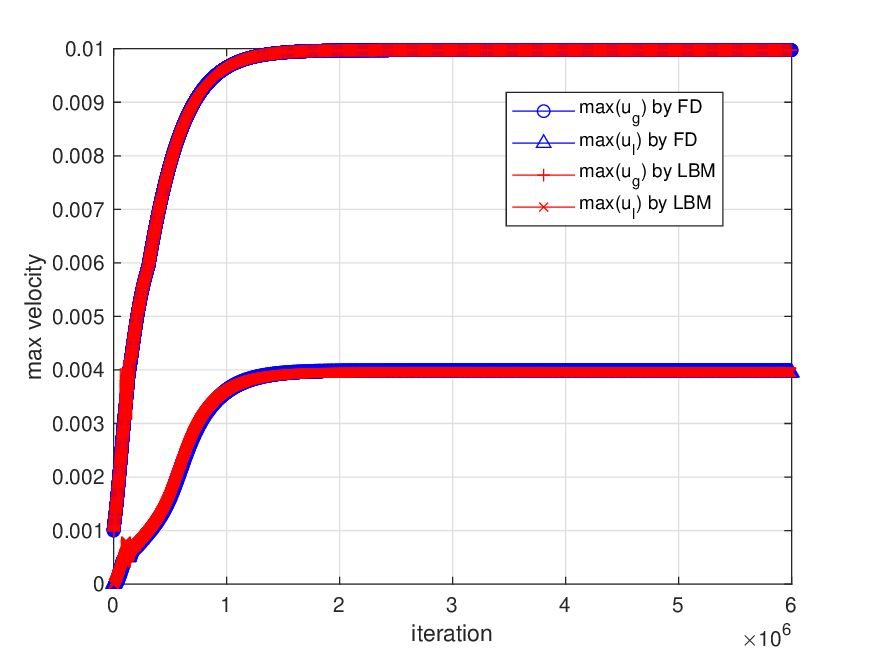}
  \caption{FD vs. LBM: max velocity}
  \label{fig3:sfig1}
\end{subfigure}
\begin{subfigure}{.5\textwidth}
  \centering
  \includegraphics[width=1.0\linewidth]{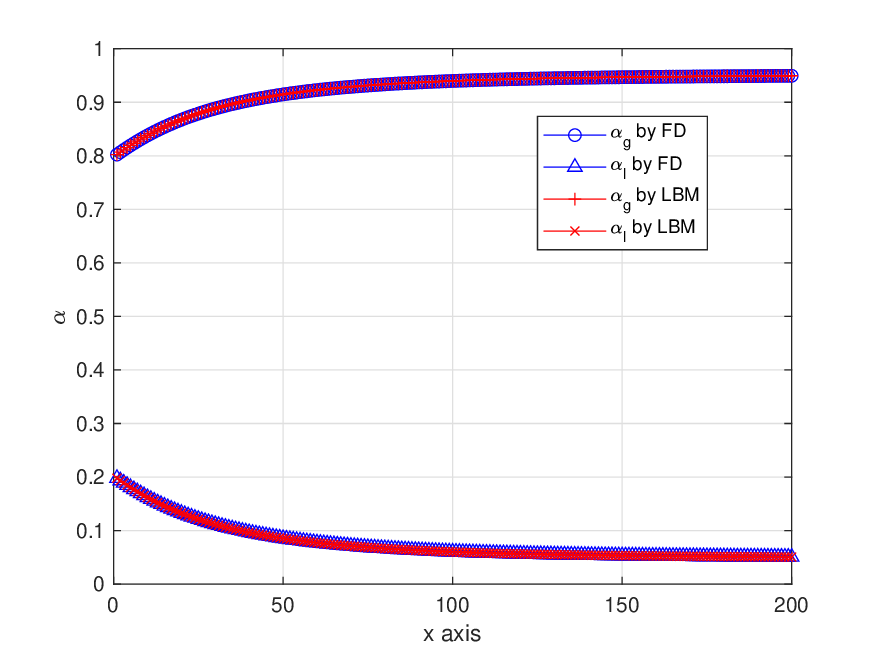}
  \caption{FD vs. LBM: volume fractions}
  \label{fig3:sfig2}
\end{subfigure}
\begin{subfigure}{.5\textwidth}
  \centering
  \includegraphics[width=1.0\linewidth]{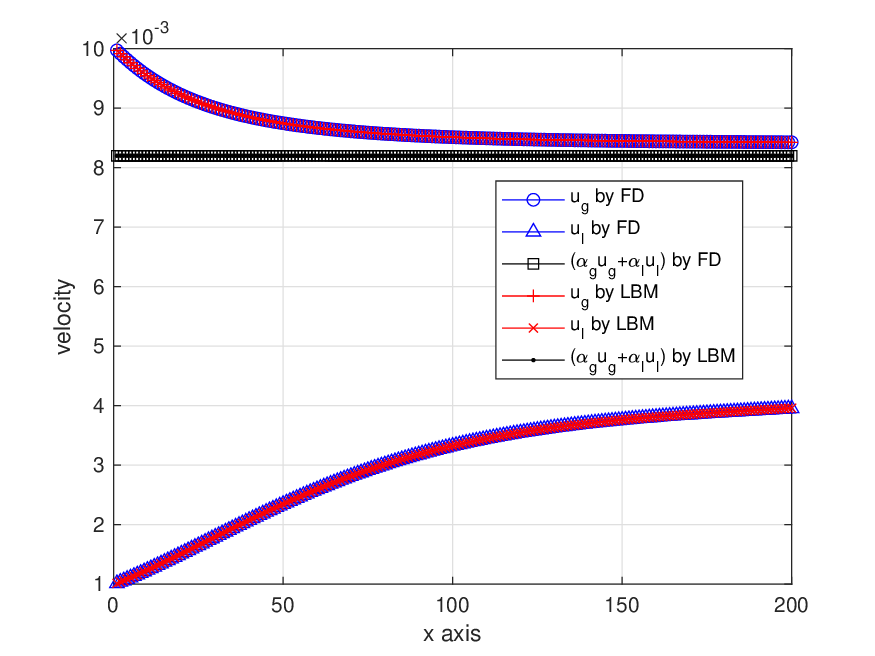}
  \caption{FD vs. LBM: phase velocities}
  \label{fig3:sfig3}
\end{subfigure}
\begin{subfigure}{.5\textwidth}
  \centering
  \includegraphics[width=1.0\linewidth]{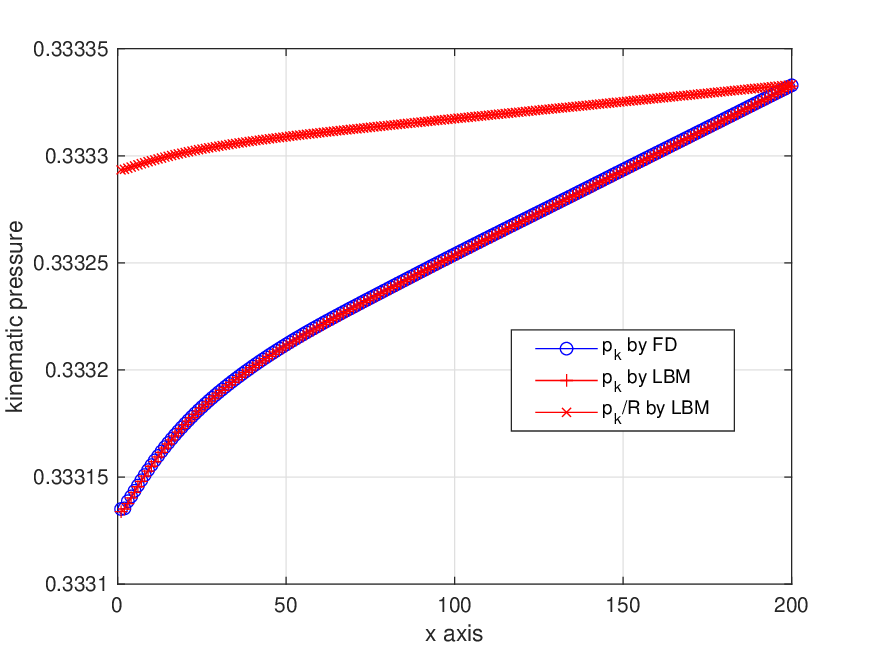}
  \caption{FD vs. LBM: kinematic pressure}
  \label{fig3:sfig4}
\end{subfigure}
\caption{TEST \#2: Comparison between numerical results by FD and by LBM ($R=5$ and $\hat{g} = 2.5\times 10^{-7}$).}
\label{fig:fig3}
\end{figure}

\begin{figure}
\begin{subfigure}{.5\textwidth}
  \centering
  \includegraphics[width=1.0\linewidth]{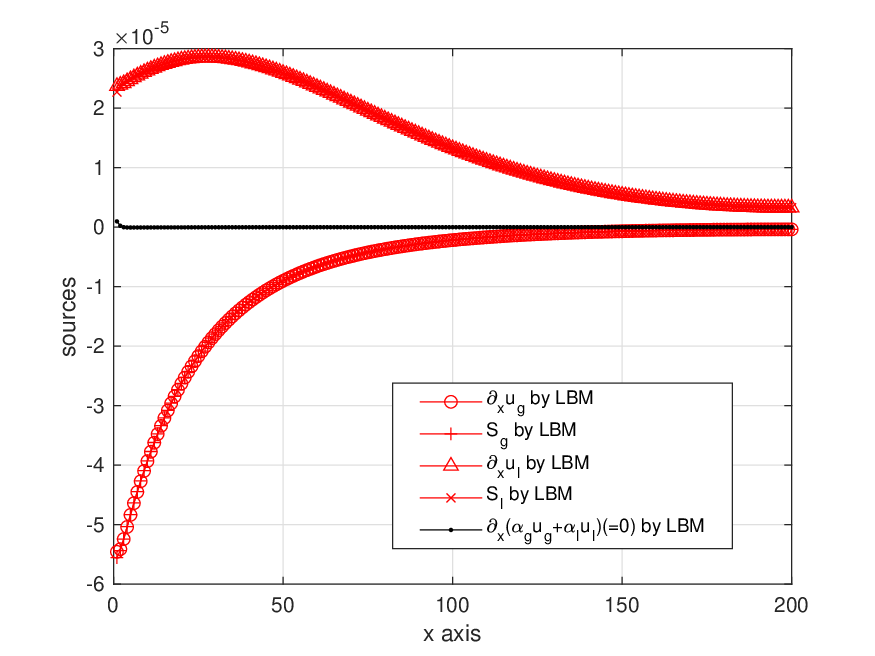}
  \caption{Phase sources computed by LBM}
  \label{fig4:sfig1}
\end{subfigure}
\begin{subfigure}{.5\textwidth}
  \centering
  \includegraphics[width=1.0\linewidth]{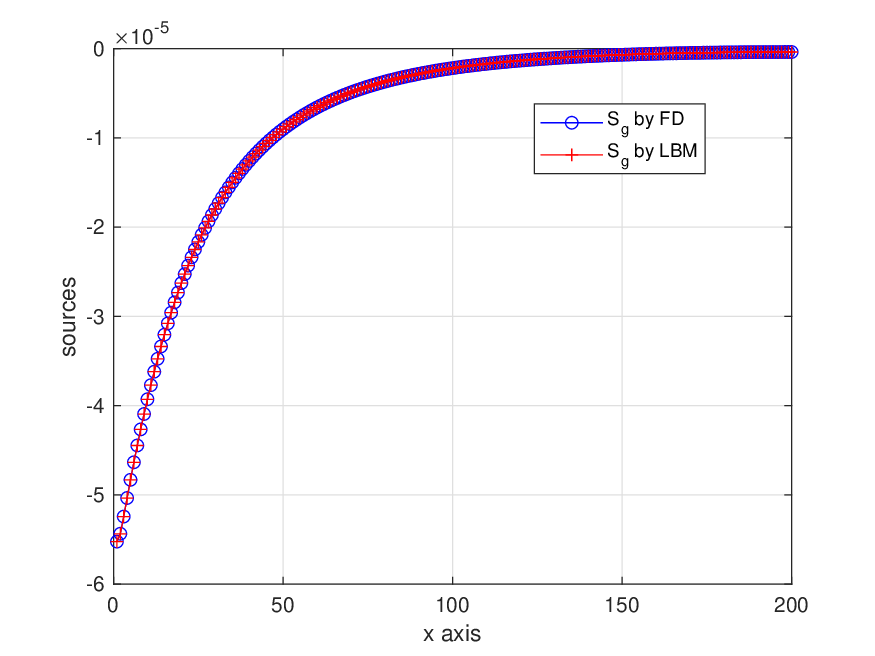}
  \caption{FD vs. LBM: dispersed phase source}
  \label{fig4:sfig2}
\end{subfigure}
\begin{subfigure}{.5\textwidth}
  \centering
  \includegraphics[width=1.0\linewidth]{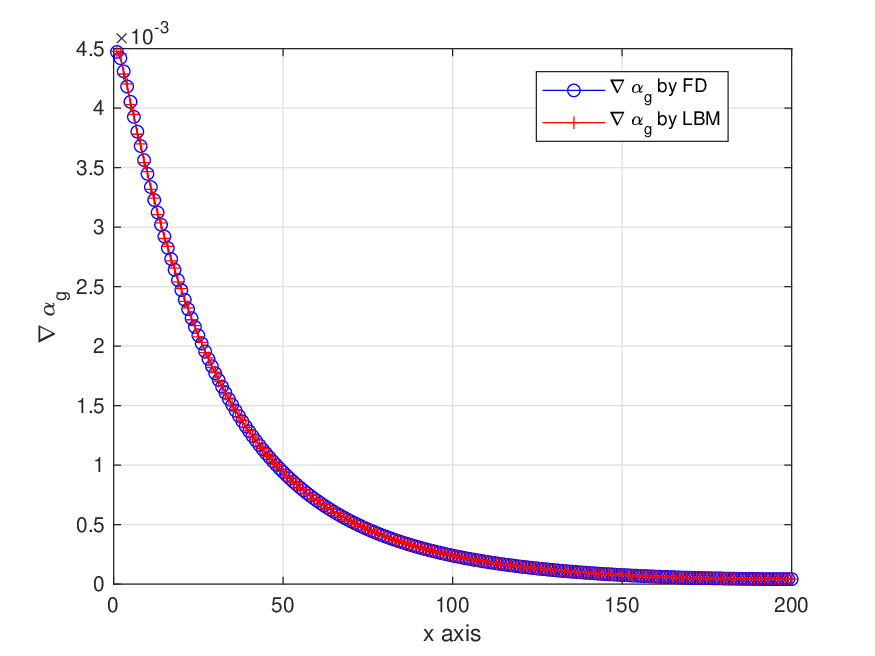}
  \caption{FD vs. LBM: dispersed phase gradient}
  \label{fig4:sfig3}
\end{subfigure}
\begin{subfigure}{.5\textwidth}
  \centering
  \includegraphics[width=1.0\linewidth]{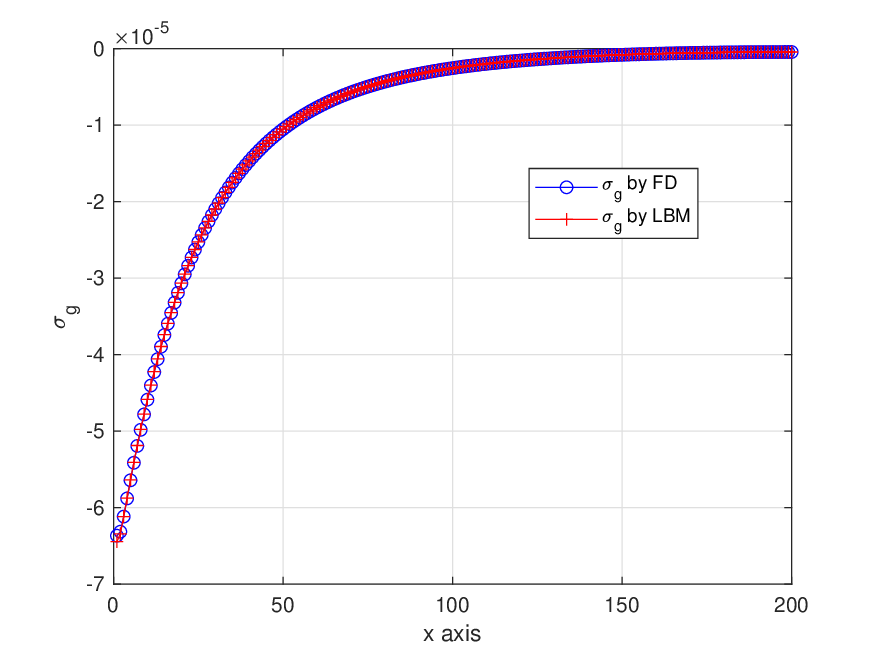}
  \caption{FD vs. LBM: dispersed phase stress}
  \label{fig4:sfig4}
\end{subfigure}
\begin{subfigure}{.5\textwidth}
  \centering
  \includegraphics[width=1.0\linewidth]{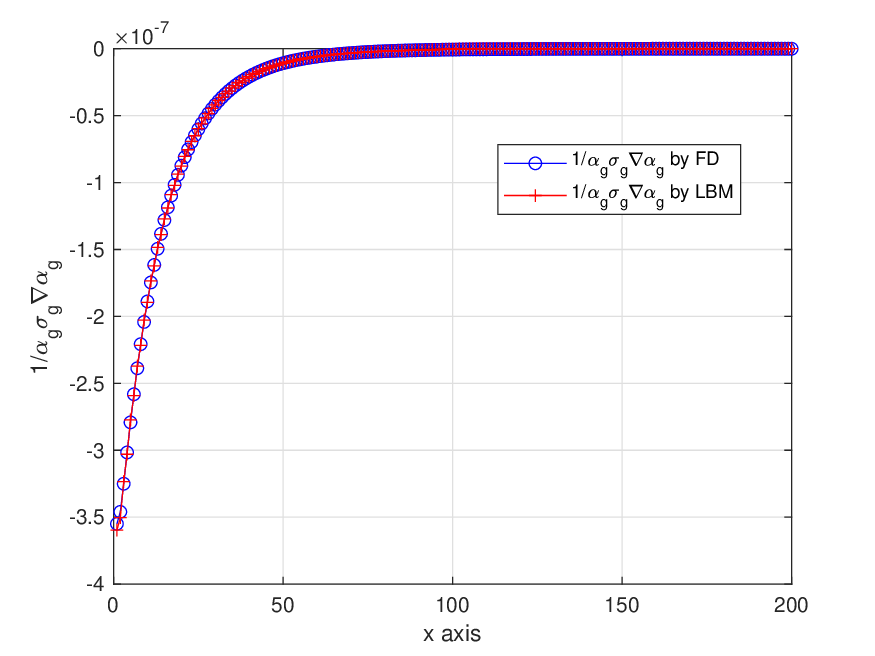}
  \caption{FD vs. LBM: dispersed phase force}
  \label{fig4:sfig5}
\end{subfigure}
\begin{subfigure}{.5\textwidth}
  \centering
  \includegraphics[width=1.0\linewidth]{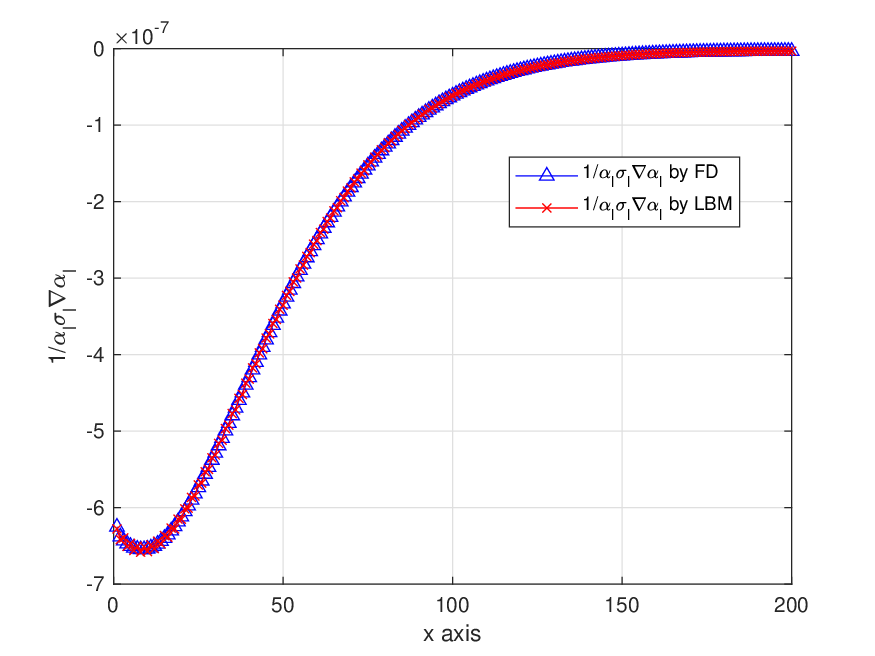}
  \caption{FD vs. LBM: liquid phase force}
  \label{fig4:sfig6}
\end{subfigure}
\caption{TEST \#2: Computational details of the LBM schemes without FD corrections and comparison with numerical results by FD operators ($R=5$ and $\hat{g} = 2.5\times 10^{-7}$).}
\label{fig:fig4}
\end{figure}

\subsection{Numerical results and discussion}

In this section, we report the numerical results for some meaningful test cases of the proposed methodology, which are grouped into a preliminary category and an advanced category. The preliminary category (TEST \#1 and \#2) is intended to realize the basic essential features of the methodology, while the advanced category (TEST \#3 and \#4) discusses challenging features, i.e. very high density ratios and realistic phenomenological relations for the interphase momentum exchange.

First of all, let us consider first the preliminary category with preliminary numerical results. In particular, let us consider the following test case (corresponding to a natural circulation loop driven by gas injection), called TEST \#1: $\rho_g^0 = 1.2$, $\rho_l^0 = 2.4$ ($R=2$), $\nu_g = 1.1667$, $\nu_l = 1.1667$, $\hat{g} = 10^{-6}$, $\hat{K}_I = 10^{-2}$ and $\hat{K}_W = 10^{-2}$. The progressively increased quantities at the inlet (see Eq. (\ref{transient})) for this test case are: $\hat{\alpha}_g^{max} = 0.8$ and $\hat{\alpha}_g^{min} = 10^{-2}$ for the dispersed phase volume fraction; $\hat{u}_l^{max} = 10^{-3}$ and $\hat{u}_l^{min} = 0$ for the liquid phase velocity; while, for the dispersed phase velocity, $\hat{u}_g^{max} = 10^{-2}$ and $\hat{u}_g^{min}$ is given by Eq. (\ref{analytical0}), namely
\begin{equation}\label{analytical0bis}
\hat{u}_g^{min}=\frac{1}{\sqrt{\hat{K}_I}}\,\sqrt{\hat{\alpha}_g^{min}\,(1-\hat{\alpha}_g^{min})\,(R-1)\,\hat{g}}.
\end{equation}
For this test case, we used $N_x = 200$ which is the number of nodal values, $N_t = 6\times10^6$ which is number of time steps and $n_t = 5\times 10^5 = N_t/12$ which is the number of time steps for the transient smooth ramp-up (see Eq. (\ref{transient})). Similarly, for the purpose to test the capability of the proposed methodology to deal with phases with different average densities, let us consider also another test case, called TEST \#2, where $\rho_g^0 = 1.2$, but this time $\rho_l^0 = 6.0$ ($R=5$), $\hat{g} = 2.5\times 10^{-7}$, in such a way that $(R-1)\,\hat{g}$ is equal in the two tests discussed so far, and all remaining parameters are the same as in TEST \#1. 

Next, let us consider the advanced category with advanced numerical results for testing the proposed methodology in challenging setups, i.e. very high density ratios and stiff (but realistic) phenomenological relations for the interphase momentum exchange. Let us consider a test case, called TEST \#3, where $\rho_g^0 = 1.2$, but this time we take the very challenging value $\rho_l^0 = 1000.0$ ($R=833.3$), $\hat{g} = 1.2\times 10^{-9}$, in such a way that $(R-1)\,\hat{g}$ is equal in all cases discussed so far. This third test case is intended to explore how the solvers deal with a very large density ratio -- typically very challenging for LBM -- which could trigger some numerical instabilities. For this reason, this test case requires also the stabilization ingredients discussed in section \ref{largedensityratios}: $\gamma = 1$ with the third stabilization strategy (for better consistency) and $n_\gamma = 200\ll N_t$. All remaining parameters are the same as in other test cases discussed so far. 

Moreover, for the purpose to test also the additional dependency of the drag force on the volume fraction, let us also consider a test case, called TEST \#4, where the interphase drag force is modeled by the CGW model. In this case, the simulation parameters must be derived from physical values, which are consistent with a realistic setup (despite the limitations of the one-dimensional configuration). For this aim, let us focus on the experimental setup described in Ref. \cite{PFLEGER19995091} for extracting the most relevant physical quantities: $g = 9.81$, $\kappa_I = 3\times0.66/(4\times0.002)=247.5$, $u_g^{max} = 0.3$ and $\alpha_g^{max}=0.1$, where we assumed $U = 1\;m/s$ as the characteristic flow speed and $L=1\;m$ as the characteristic length scale of the flow field. 
%
%ugstar = 0.3; % [m/s], Fig. 2 p. 5
%gstar = 9.81; % [m/s^2]
%alphagmax = 0.1; % Fig. 5 p. 7
%kappaIstar = 3*(0.66)/(4*0.002); % [1/m]
%kappaIstar = 247.5 1/m
%
Assuming $\hat{g} = 1.2\times 10^{-9}$ and $\hat{\kappa}_I=1.45\times10^{-4}$ ($\hat{\kappa}_W = \hat{\kappa}_I$) for stability reasons, applying Eq. (\ref{scaling-grav-acceleration}) and (\ref{scaling-kappaI}) yields $c/\tau=8.18\times10^{9}$ and $c\,\tau=5.86\times10^{-7}$, respectively. Therefore the lattice speed is $c = \sqrt{(c/\tau)\,(c\tau)}=69.2$. This means that $\hat{u}_g^{max}=u_g^{max}/c = 0.0043$ which ensures a Mach number small enough to be consistent with the incompressible limit. Moreover, let us assume $\hat{u}_l^{max}=\hat{u}_g^{max}/4=0.0011$ for having a balanced inlet condition, as we prove below. Adopting $u_g^{\dagger}=u_g^{max}$, $u_l^{\dagger}=u_l^{max}$ and $\alpha_g^{\dagger}=\alpha_g^{max}$ as reference values, it is possible to estimate the ratio between the momentum exchange force and the buoyancy force in the momentum equation of the gas phase by Eq. (\ref{Mex}), which is equal to $M_{ex} = 1.04$. This means that the drag force of this setup is expected to almost balance the buoyancy force at the inlet. Finally, let us assume $\hat{u}_g^{min} = \hat{u}_g^{max}/10$, $\hat{u}_l^{min} = \hat{u}_l^{max}/10$ and $\alpha_g^{min}=10^{-2}$. All remaining parameters are the same as in other test cases discussed
so far. 

\begin{figure}
\begin{subfigure}{.5\textwidth}
  \centering
  \includegraphics[width=1.0\linewidth]{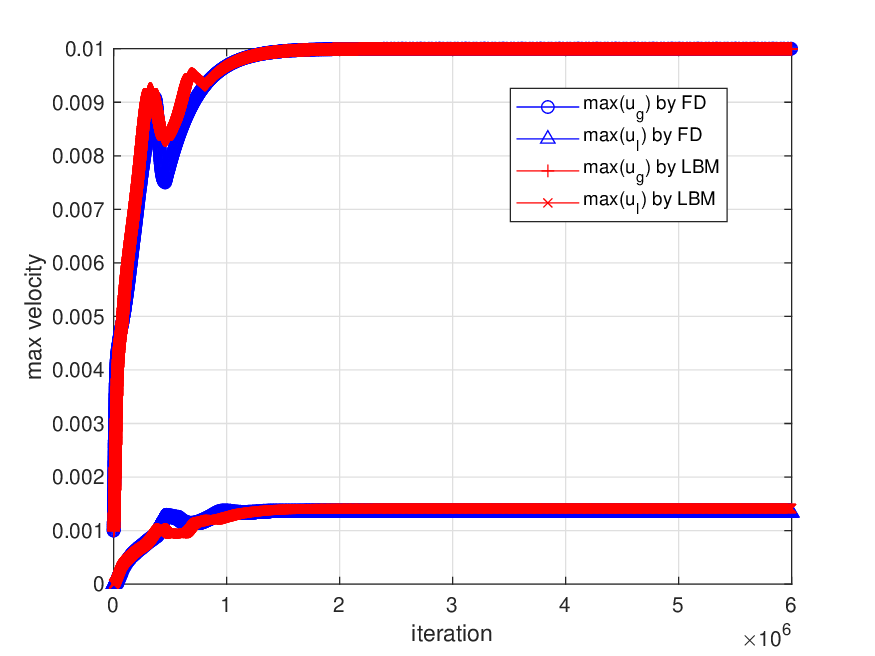}
  \caption{FD vs. LBM: max velocity}
  \label{fig5:sfig1}
\end{subfigure}
\begin{subfigure}{.5\textwidth}
  \centering
  \includegraphics[width=1.0\linewidth]{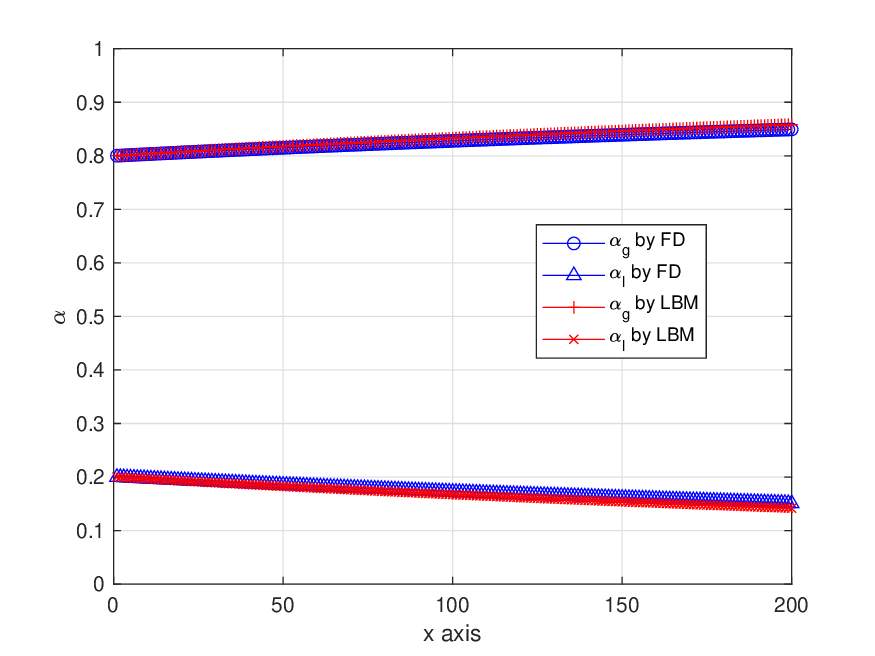}
  \caption{FD vs. LBM: volume fractions}
  \label{fig5:sfig2}
\end{subfigure}
\begin{subfigure}{.5\textwidth}
  \centering
  \includegraphics[width=1.0\linewidth]{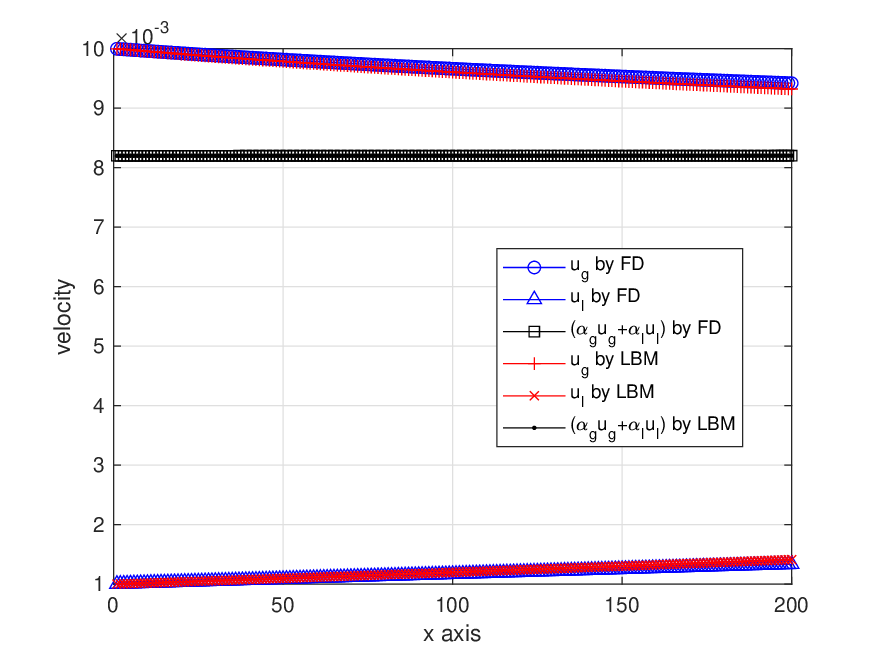}
  \caption{FD vs. LBM: phase velocities}
  \label{fig5:sfig3}
\end{subfigure}
\begin{subfigure}{.5\textwidth}
  \centering
  \includegraphics[width=1.0\linewidth]{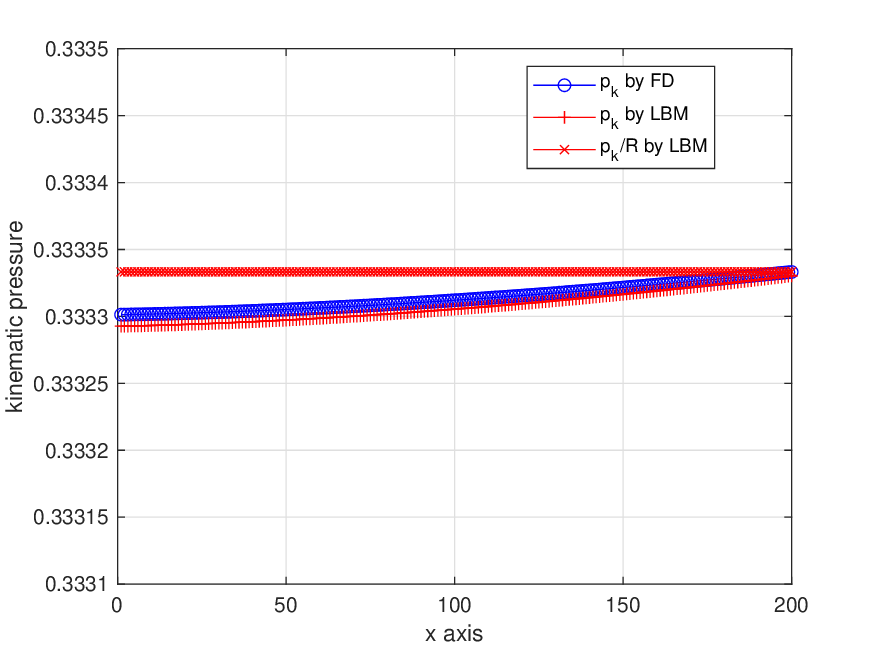}
  \caption{FD vs. LBM: kinematic pressure}
  \label{fig5:sfig4}
\end{subfigure}
\caption{TEST \#3: Comparison between numerical results by FD and by LBM ($R=833.3$ and $\hat{g} = 1.2\times 10^{-9}$).}
\label{fig:fig5}
\end{figure}

\begin{figure}
\begin{subfigure}{.5\textwidth}
  \centering
  \includegraphics[width=1.0\linewidth]{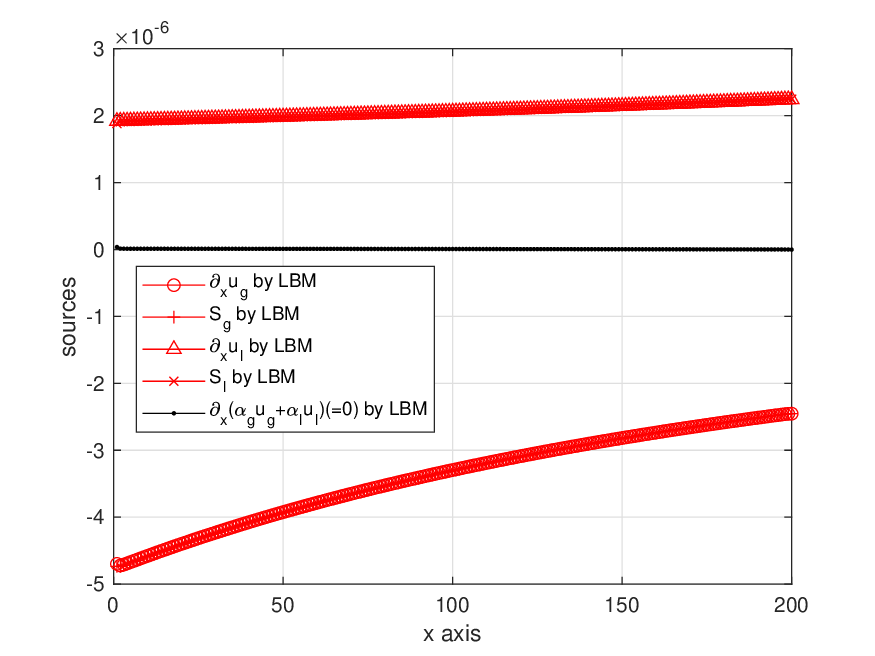}
  \caption{Phase sources computed by LBM}
  \label{fig6:sfig1}
\end{subfigure}
\begin{subfigure}{.5\textwidth}
  \centering
  \includegraphics[width=1.0\linewidth]{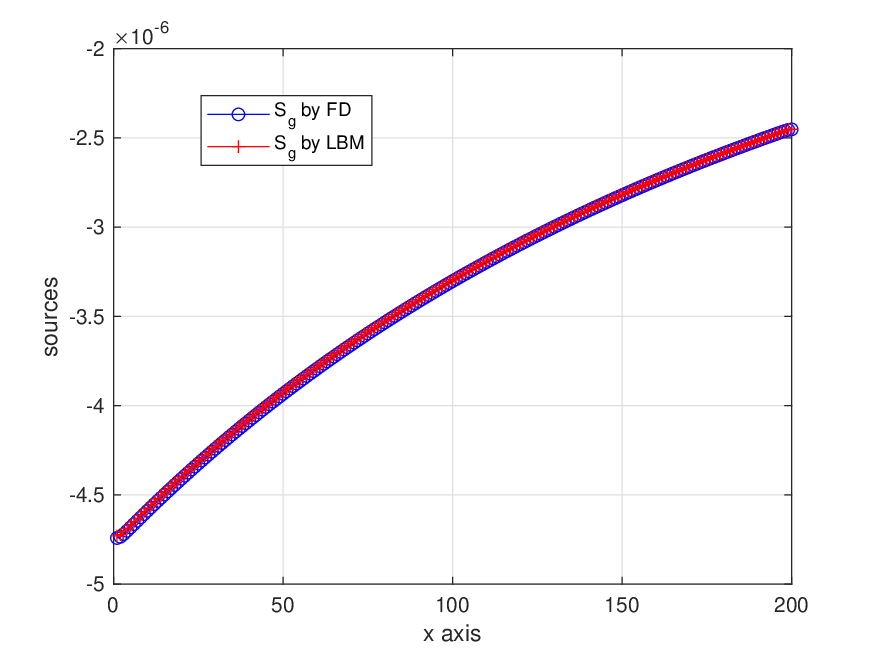}
  \caption{FD vs. LBM: dispersed phase source}
  \label{fig6:sfig2}
\end{subfigure}
\begin{subfigure}{.5\textwidth}
  \centering
  \includegraphics[width=1.0\linewidth]{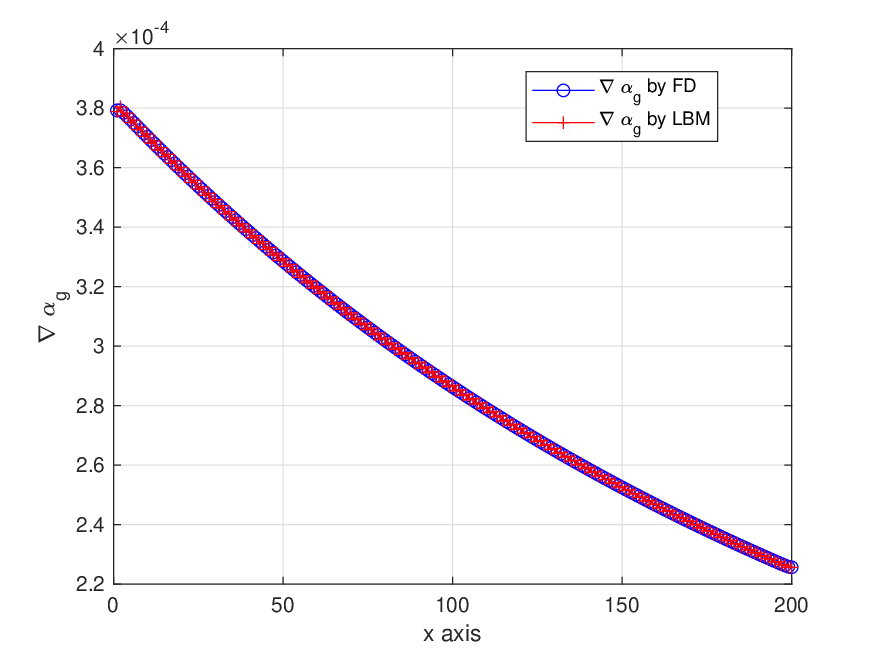}
  \caption{FD vs. LBM: dispersed phase gradient}
  \label{fig6:sfig3}
\end{subfigure}
\begin{subfigure}{.5\textwidth}
  \centering
  \includegraphics[width=1.0\linewidth]{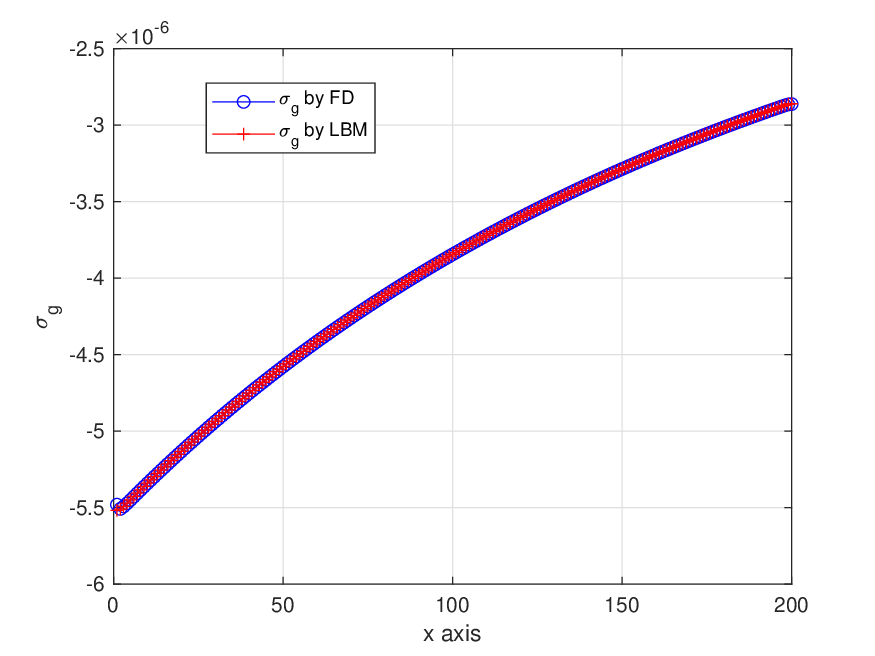}
  \caption{FD vs. LBM: dispersed phase stress}
  \label{fig6:sfig4}
\end{subfigure}
\begin{subfigure}{.5\textwidth}
  \centering
  \includegraphics[width=1.0\linewidth]{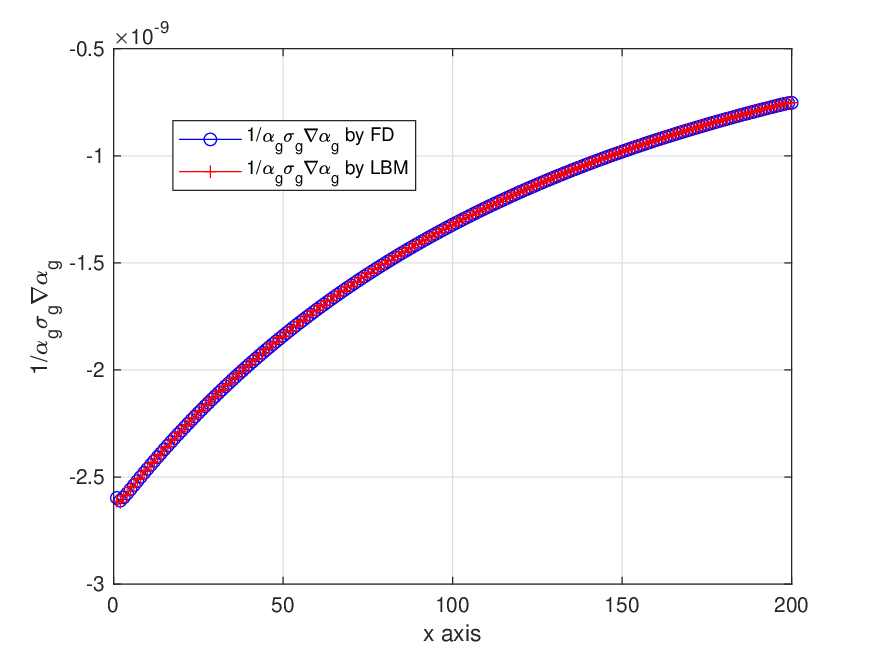}
  \caption{FD vs. LBM: dispersed phase force}
  \label{fig6:sfig5}
\end{subfigure}
\begin{subfigure}{.5\textwidth}
  \centering
  \includegraphics[width=1.0\linewidth]{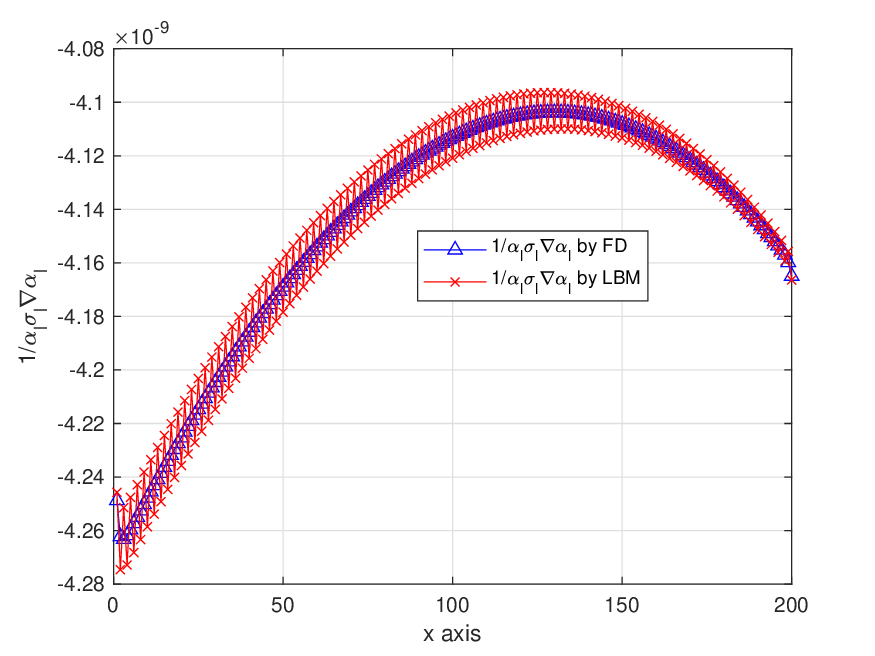}
  \caption{FD vs. LBM: liquid phase force}
  \label{fig6:sfig6}
\end{subfigure}
\caption{TEST \#3: Computational details of the LBM schemes without FD corrections and comparison with numerical results by FD operators ($R=833.3$ and $\hat{g} = 1.2\times 10^{-9}$).}
\label{fig:fig6}
\end{figure}

\begin{figure}
\begin{subfigure}{.5\textwidth}
  \centering
  \includegraphics[width=1.0\linewidth]{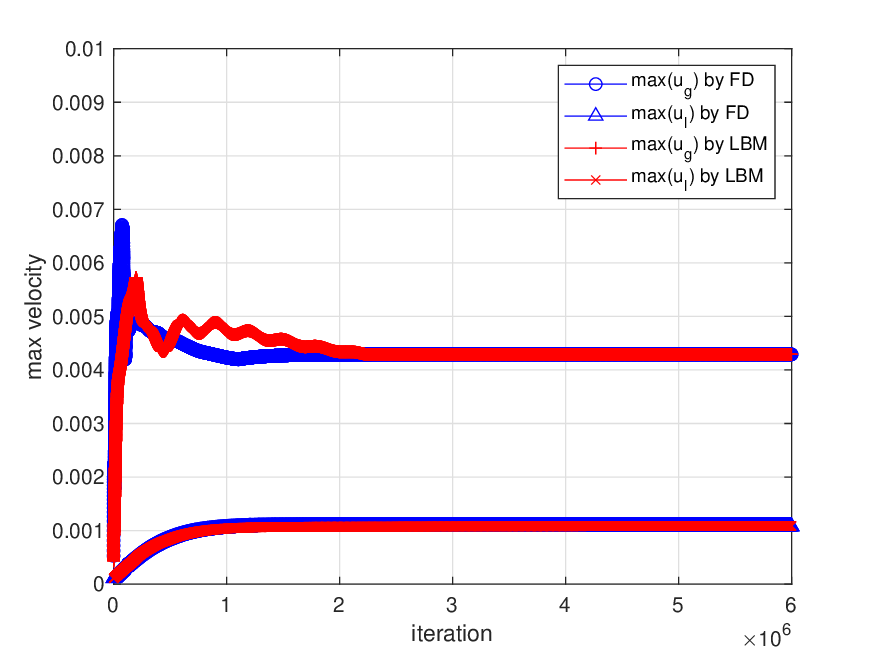}
  \caption{FD vs. LBM: max velocity}
  \label{fig7:sfig1}
\end{subfigure}
\begin{subfigure}{.5\textwidth}
  \centering
  \includegraphics[width=1.0\linewidth]{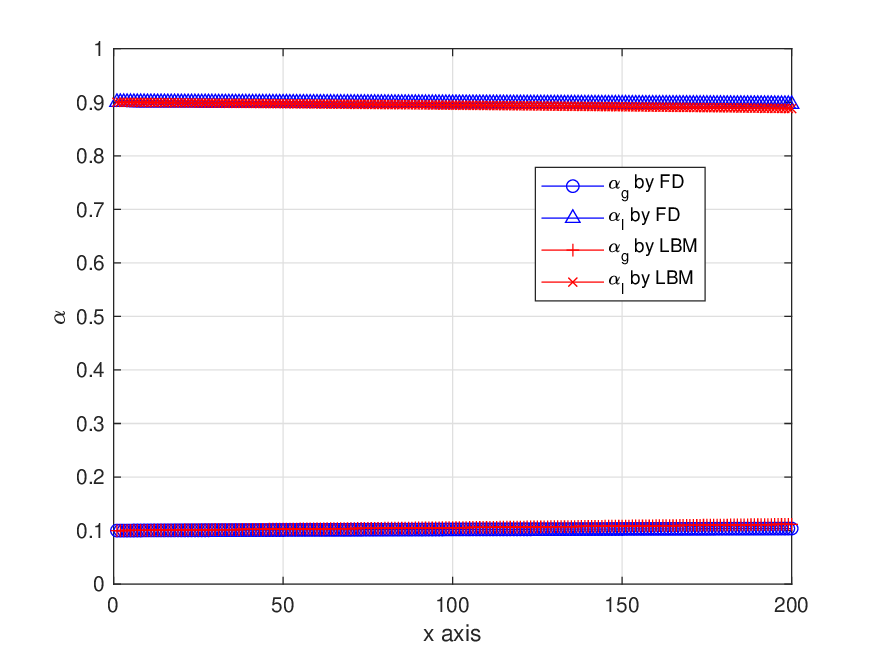}
  \caption{FD vs. LBM: volume fractions}
  \label{fig7:sfig2}
\end{subfigure}
\begin{subfigure}{.5\textwidth}
  \centering
  \includegraphics[width=1.0\linewidth]{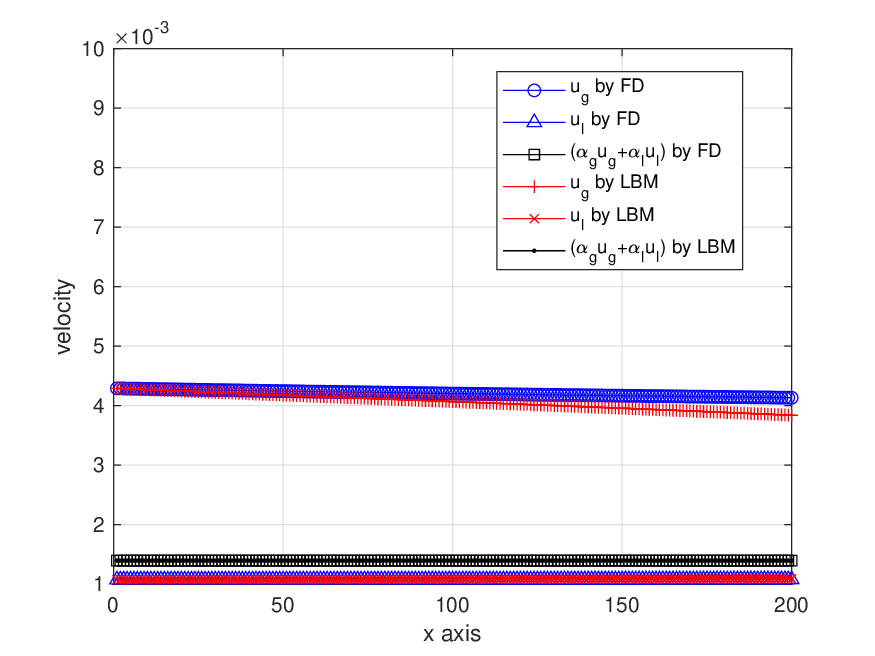}
  \caption{FD vs. LBM: phase velocities}
  \label{fig7:sfig3}
\end{subfigure}
\begin{subfigure}{.5\textwidth}
  \centering
  \includegraphics[width=1.0\linewidth]{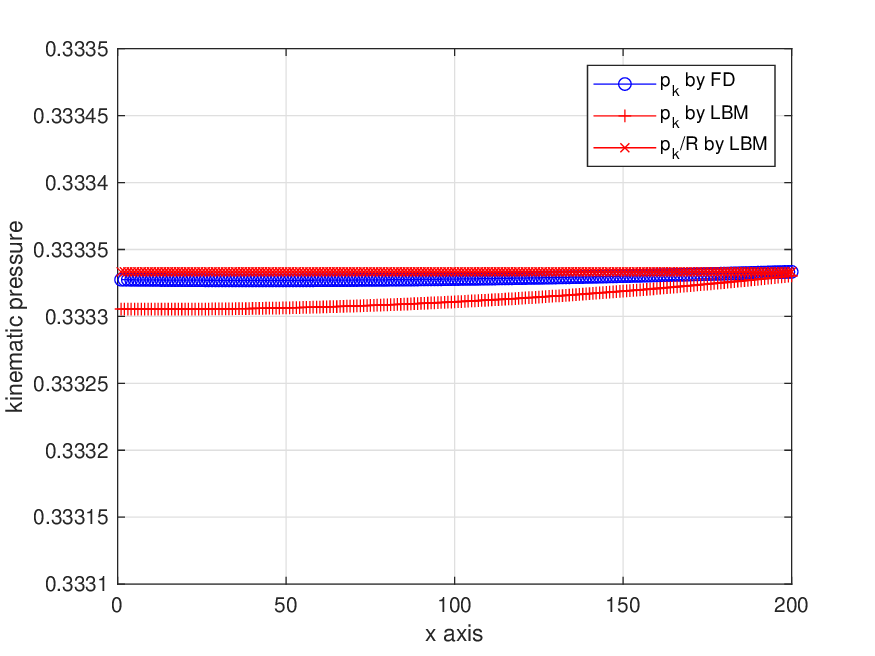}
  \caption{FD vs. LBM: kinematic pressure}
  \label{fig7:sfig4}
\end{subfigure}
\caption{TEST \#4: Comparison between numerical results by FD and by LBM ($R=833.3$, $\hat{g} = 1.2\times 10^{-9}$, model for drag force by Clift, Grace \& Weber with $\hat{\kappa}_I=1.45\times10^{-4}$).}
\label{fig:fig7}
\end{figure}

\begin{figure}
\begin{subfigure}{.5\textwidth}
  \centering
  \includegraphics[width=1.0\linewidth]{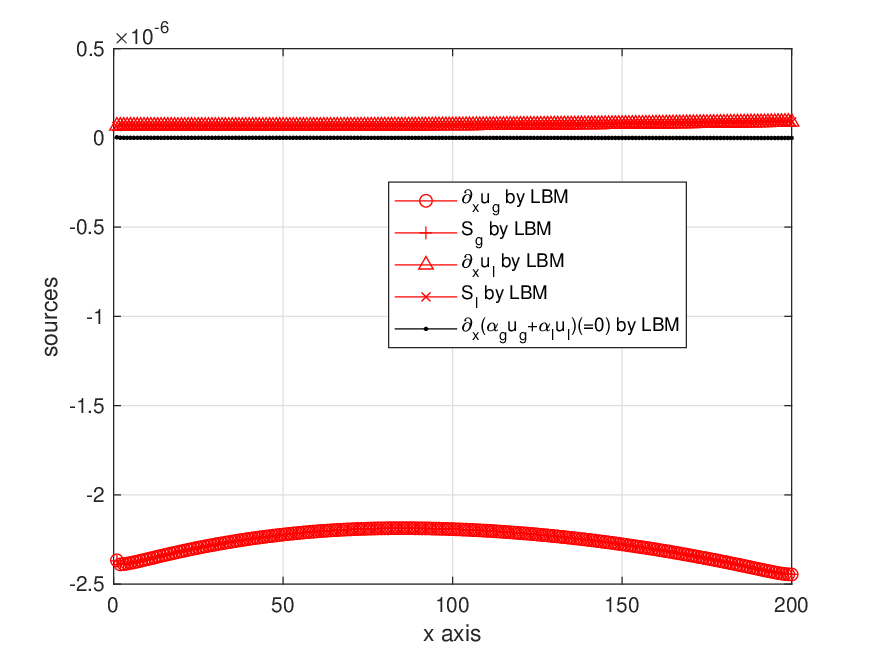}
  \caption{Phase sources computed by LBM}
  \label{fig8:sfig1}
\end{subfigure}
\begin{subfigure}{.5\textwidth}
  \centering
  \includegraphics[width=1.0\linewidth]{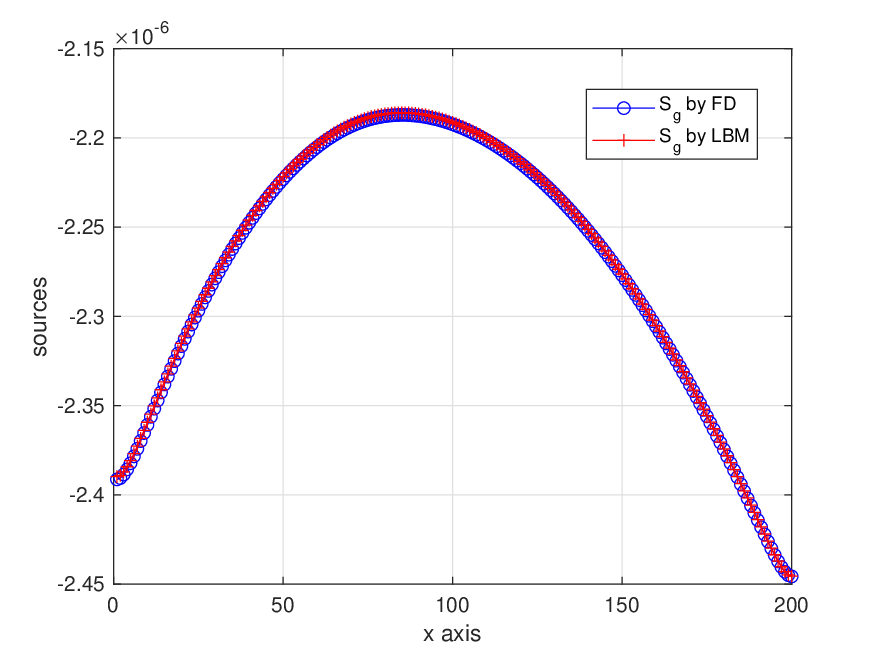}
  \caption{FD vs. LBM: dispersed phase source}
  \label{fig8:sfig2}
\end{subfigure}
\begin{subfigure}{.5\textwidth}
  \centering
  \includegraphics[width=1.0\linewidth]{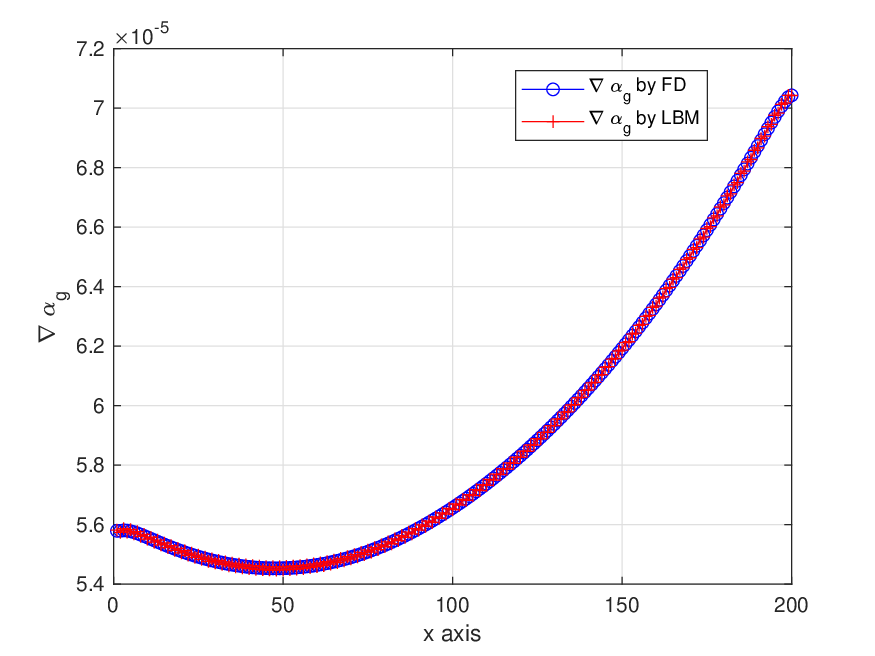}
  \caption{FD vs. LBM: dispersed phase gradient}
  \label{fig8:sfig3}
\end{subfigure}
\begin{subfigure}{.5\textwidth}
  \centering
  \includegraphics[width=1.0\linewidth]{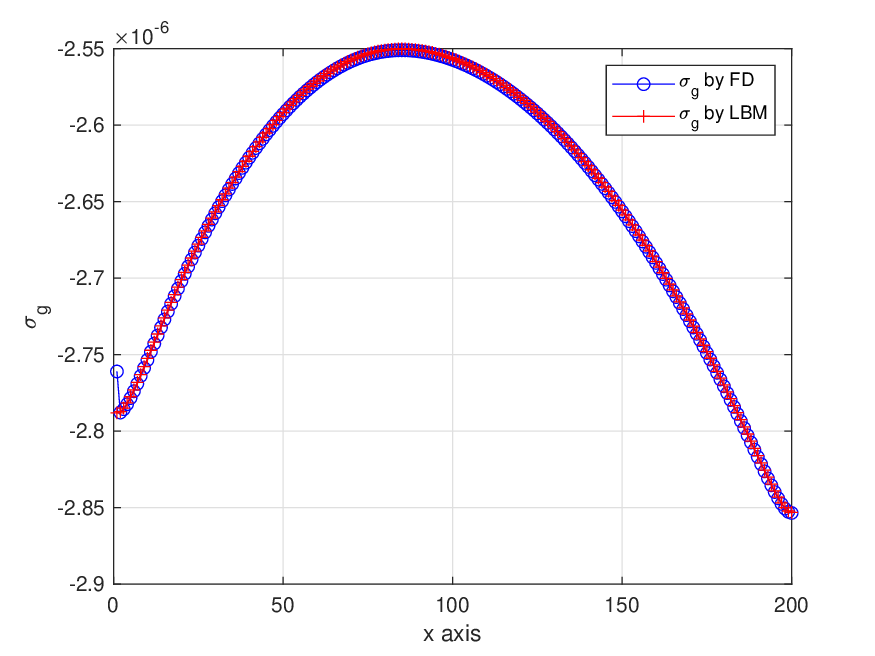}
  \caption{FD vs. LBM: dispersed phase stress}
  \label{fig8:sfig4}
\end{subfigure}
\begin{subfigure}{.5\textwidth}
  \centering
  \includegraphics[width=1.0\linewidth]{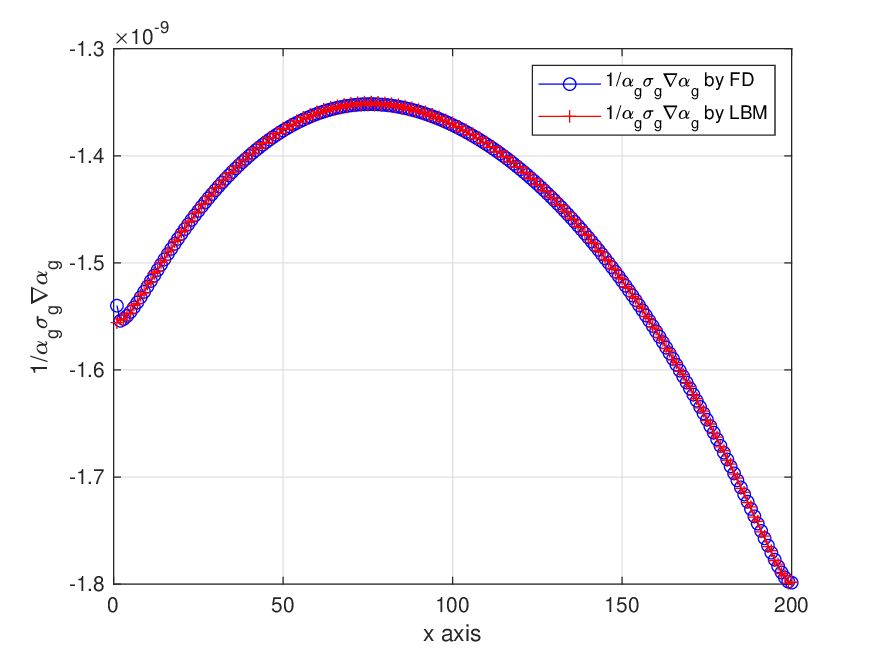}
  \caption{FD vs. LBM: dispersed phase force}
  \label{fig8:sfig5}
\end{subfigure}
\begin{subfigure}{.5\textwidth}
  \centering
  \includegraphics[width=1.0\linewidth]{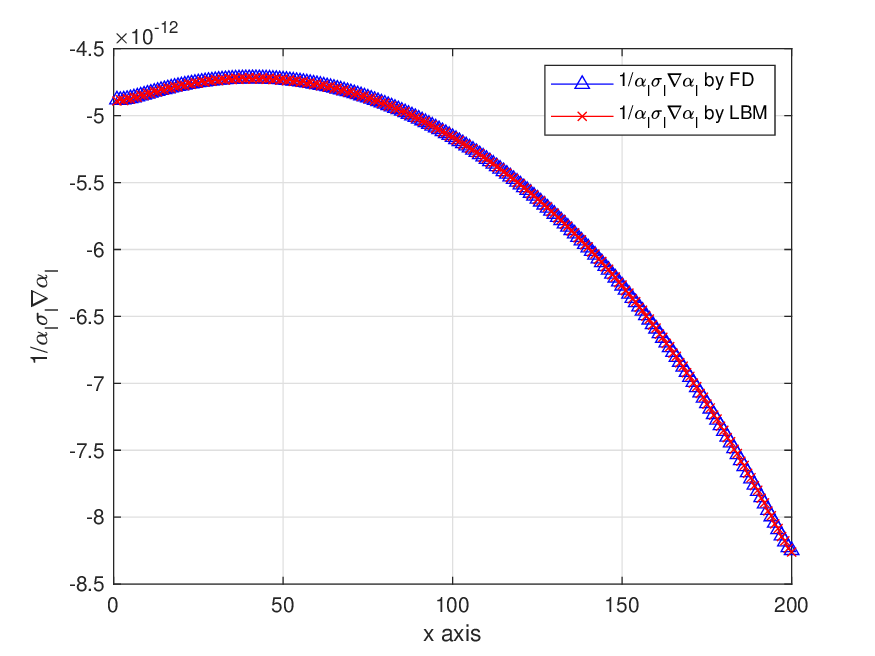}
  \caption{FD vs. LBM: liquid phase force}
  \label{fig8:sfig6}
\end{subfigure}
\caption{TEST \#4: Computational details of the LBM schemes without FD corrections and comparison with numerical results by FD operators ($R=833.3$, $\hat{g} = 1.2\times 10^{-9}$, model for drag force by Clift, Grace \& Weber with $\hat{\kappa}_I=1.45\times10^{-4}$).}
\label{fig:fig8}
\end{figure}

Numerical results are reported below for all test cases, comparing those obtained by the FD engine and those by the LBM engine. Fig. \ref{fig:fig1} and \ref{fig:fig3} show the most relevant features for TEST \#1 and \#2 respectively, while Fig. \ref{fig:fig2} and Fig. \ref{fig:fig4} show additional computational details again for the same test cases. Moreover, Fig. \ref{fig:fig5} shows the main features of the TEST \#3 about a very large density ratio, together with additional computational details for the same test case in Fig. \ref{fig:fig6}. Even more remarkably, Fig. \ref{fig:fig7} shows the main features of the TEST \#4 about a very large density ratio with a stiff (but realistic) relation for the drag force. More computational details about this very challenging test case can be found in Fig. \ref{fig:fig8}.

Concerning these numerical results, let us comment first the preliminary ones below. 
\begin{itemize}
    \item It is important to highlight that the proposed LBM schemes agree with the FD numerical results in an excellent way. This opens a promising venue for simulating multiphase flows on large HPC facilities. In fact, it is worth to highlight that the proposed methodology and all LBM formulas reported in section \ref{LBMdesign} can be automatically applied to any dimension, including all boundary conditions.
    \item The key point, namely the artificially compressible continuity equation for each phase, works well and sub Fig. \ref{fig2:sfig1} and Fig. \ref{fig4:sfig1} prove that the individual phase sources can be computed in LBM schemes without FD corrections.
    \item LBM schemes without FD corrections can be used effectively to compute all terms in multiphase equations, including gradients of the volume fractions and stress tensors, as clearly reported in sub Figs. \ref{fig2:sfig2}-\ref{fig2:sfig6} and sub Figs. \ref{fig4:sfig2}-\ref{fig4:sfig6}, where the LBM results show an excellent agreement with those obtained by FD operators.
    \item It is important to highlight that incompressible multiphase flows has divergence-free mixture velocity, as prescribed by Eq. (\ref{continuitysum}), but this condition does not hold for individual phase velocity. The one-dimensional test case is enough to show this important difference: indeed $\hat{\alpha}_g \hat{{u}}_g + \hat{\alpha}_l \hat{{u}}_l$ is constant (the only way to ensure the divergence free condition in one dimension), but $\hat{{u}}_g$ and $\hat{{u}}_l$ are not. It is clear that sub Fig. \ref{fig1:sfig3} and sub Fig. \ref{fig3:sfig3} confirm this expectation. Hence the proposed LBM framework ensures, as it should be, the divergence free condition for the mixture velocity, but it does not force it also for the individual phase velocities. This is one of the most important features of the proposed method in comparison with what is already available in the literature. 
%    
% ug(200) = 0.008423711190781
% ul(200) = 0.003975027343349
% alphag(200) = 0.949686618689092
% alphal(200) = 0.050313381310908
% alphag(200)*ug(200)+alphal(200)*ul(200) = 0.008199882864033
% r=(Kw/Ki)*alphag(200) = 0.949686618689092
% ug0 = 1/sqrt(Ki)*sqrt(alphag(200)*(1-alphag(200))*(R-1)*gg) = 0.002185908163029
% (ug(200)-sqrt(ug(200)^2-(1-r)*(ug(200)^2-ug0^2)))/(1-r) = 0.003975437432568
% (ug(200)-ug0)*(1+ug0/ug(200))/2 = 0.003928239835722
%
    \item Another consideration is to investigate the relation between the numerical solution at the outlet and the analytical solution derived in the previous section. For the sake of simplicity, let us consider TEST \#1. At the outlet, the LBM framework provides $\hat{u}_g(N_x) = 8.4237\times 10^{-3}$, $\hat{u}_l(N_x) = 3.9750\times 10^{-3}$, $\hat{\alpha}_g(N_x) = 0.950$ and $\hat{\alpha}_l(N_x) = 0.050$, which corresponds to the following mixture velocity $\hat{\alpha}_g(N_x)\,\hat{u}_g(N_x)+\hat{\alpha}_l(N_x)\,\hat{u}_l(N_x) = 8.1999\times 10^{-3}$. In this test case, $r = 0.9496$. Using $\hat{\alpha}_g(N_x)$ instead of $\overline{\alpha}_g$ into Eq. (\ref{analytical0}) leads to $\overline{u}_g^0 = 2.1859\times 10^{-3}$. Using the latter velocity into Eq. (\ref{analyticalsolution}) and assuming $\overline{u}_g=\hat{u}_g(N_x)$ leads to $\overline{u}_l = 3.9754\times 10^{-3}$ which is in excellent agreement with $\hat{u}_l(N_x)$ (mismatch $0.01$ \%). Because in this case $r$ is pretty close to $1$, then it is also possible to use the approximated formula $\overline{{u}}_l\approx(\overline{{u}}_g-\overline{{u}}_g^0)\,(1+\overline{{u}}_g^0/\overline{{u}}_g)/2 = 3.9282\times 10^{-3}$ (mismatch $1$ \%).
    \item It is also worth to compare the numerical results between TEST \#1 and TEST \#2. It is clear in sub Fig. (\ref{fig1:sfig3}) and sub Fig. (\ref{fig3:sfig3}) the difference in the slope of the kinematic pressure profiles of the two phases for the TEST \#1 ($R=2$) and TEST \#2 ($R=5$), respectively. This proves that the proposed methodology can describe properly the effects due to the density ratio between the phases. 
\end{itemize}

Secondly, advanced results are discussed below. 
\begin{itemize}
    \item The previous preliminary results highlight some challenges, mainly about large density ratios, which are discussed here in the category about advanced results. In the limit of very large density ratios, namely $R\gg 1$, some terms proportional to $1/R$ in the momentum equation of the liquid phase become very small and hence comparable with the numerical errors. For this reason, TEST \#3 ($R=833.3$) also requires the stabilization ingredients discussed in section \ref{largedensityratios}: $\gamma = 1$ with the third stabilization strategy and $n_\gamma = 200\ll N_t$. The transient profiles of the LBM results are slightly smoother than those of the FD results, as shown in Fig. (\ref{fig5:sfig1}). The kinematic pressure profiles computed by the two engines in Fig. (\ref{fig5:sfig4}) reveal a small discrepancy, but this is acceptable because the pressure gradient is what really matters. The numerical details of the proposed methodology reported in Fig. (\ref{fig:fig6}) are excellent, with the exception of the computed second term of the force ${G}_l$ (in the one-dimensional case) given by Eq. (\ref{Gl}) and reported in Fig. (\ref{fig6:sfig6}), which shows some numerical oscillations (but the magnitude of these oscillations is extremely small).
    \item An important feature to be discussed in the category about the advanced numerical results is the need to include realistic phenomenological relations. In particular, TEST \#4 adopts the model by Clift, Grace and Weber (CGW) \cite{PFLEGER19995091}, which introduces a further dependence on the volume fraction into the effective drag coefficient (stiff coupling). The transient profiles of the LBM results show some discrepancies with regards to those of the FD results, as shown in Fig. (\ref{fig7:sfig1}). In particular, the LBM solution is not necessarily smoother than the reference solution, likely because of the stiffness of the CGW model. At steady state, in spite of the excellent numerical details reported in Fig. (\ref{fig:fig8}), there is a small mismatch for this mesh in the volume fraction profiles reported in Fig. (\ref{fig7:sfig2}) and in the velocity profiles reported in Fig. (\ref{fig7:sfig3}) (e.g. looking at the gas phase, there is 8\% discrepancy at the outlet). Moreover, the kinematic pressure profiles computed by the two engines in Fig. (\ref{fig7:sfig4}) reveal a negligible discrepancy in absolute terms, but which may have an impact on the effective pressure gradient and hence on the momentum equations. It is important to recall that a small discrepancy in $\alpha_g$ at small values, but with very large density ratios, may lead to a significant impact on $\Lambda(\alpha_g)$ and hence on the drag force acting in the momentum equation.
\end{itemize}

\section{Conclusions}
In spite of its apparent simplicity, LBM requires a careful design of the overall framework to be used for solving complex systems of partial differential equations, as needed by multiphase flows. By framework we mean (at least) the following information: minimum set of distribution functions, the definitions of the relevant target quantities as functions of the distribution functions (which are just auxiliaries), the definitions of the equilibrium distribution functions and how to compute them, including the proper equation of states, the proper additional forces, sources, etc. The framework is essential because LBM is defined in an auxiliary space of variables and also, more importantly, because only a subset of this auxiliary space ensures stable calculations (e.g. zeroth-order moment of the distribution function must be usually not too small for avoiding singularities, lattice sound speed cannot be too large, etc.).

For the first time to our knowledge, this work proposes a novel LBM framework to solve Eulerian-Eulerian multiphase flow equations, without any finite difference correction, including very large density ratios and also a realistic model for the drag coefficient. The proposed methodology and all reported LBM formulas can be applied to any dimension. This opens a promising venue for simulating multiphase flows in large HPC facilities and on novel parallel hardware. This LBM framework consists of six LBM schemes for the two phases, namely: two LBM schemes ($f_g$ and $f_l$) for artificially compressible continuity equations and momentum equations; two LBM schemes ($f_{\alpha g}$ and $f_{\alpha l}$) for phase volume fractions; two LBM schemes ($f_{\beta g}$ and $f_{\beta l}$) for phase continuity sources. All of these schemes are run on the same lattice and are coupled with each other, ensuring the best synergy for efficient implementation in large codes with minimum effort. 

\section*{Acknowledgments}
This publication is part of the project PNRR-NGEU which has received funding from the MUR – DM 352/2022.

\appendix
\section{Computational details for D1Q3 lattice}\label{D1Q3}

In this Appendix, we report the computational formulas of the equilibrium functions used in the proposed LBM framework for the one-dimensional case. Let us consider the D1Q3 lattice \cite{KruegerLBM}, where the generic velocity becomes a scalar and it belongs to the lattice $\mathbb{L}$, namely ${v}_q\in \mathbb{L}$, where
\begin{equation}\label{D1Q3lattice}
\mathbb{L} = 
\begin{bmatrix}
0\\
1\\
-1
\end{bmatrix}.
\end{equation}
It is worth to highlight that the D1Q3 lattice is just a projection of the lattices with higher dimensionality on the axis $\mathbf{e}_x$ \cite{KruegerLBM} and this makes one confident that it is easy, in the LBM context, to extend the proposed framework to higher dimensionality. 

Concerning the equilibrium functions, standard (S) equilibrium given by Eq. (\ref{equilibrium}) becomes
\begin{equation}\label{equilibriumD1Q3}
f_S^{eq}(\hat{\alpha},\hat{{u}})
\equiv f^{eq}(\hat{\alpha},\hat{{u}}) = 
\begin{bmatrix}
\hat{\alpha}\,(2/3-\hat{{u}}^2)\\
\hat{\alpha}\,(1+3\,\hat{{u}}+3\,\hat{{u}}^2)/6\\
\hat{\alpha}\,(1-3\,\hat{{u}}+3\,\hat{{u}}^2)/6
\end{bmatrix}.
\end{equation}
Incompressible (I) equilibrium given by Eq. (\ref{equilibriumI}) becomes
\begin{equation}\label{equilibriumID1Q3}
f_{I}^{eq}(\phi,\hat{\epsilon},\hat{{u}}) = 
\begin{bmatrix}
\hat{\epsilon}\,(3-\phi)/3-\hat{{u}}^2\\
(\hat{\epsilon}\,\phi+3\,\hat{{u}}+3\,\hat{{u}}^2)/6\\
(\hat{\epsilon}\,\phi-3\,\hat{{u}}+3\,\hat{{u}}^2)/6
\end{bmatrix}.
\end{equation}
Finally linearized equilibrium given by Eq. (\ref{forcingoperator}) becomes
\begin{equation}\label{forcingoperatorD1Q3}
f^{eq}_L(\psi,\hat{S},\hat{{G}})=
\begin{bmatrix}
\hat{S}\,(3-\psi)/3\\
(\hat{S}\,\psi+3\,\hat{{G}})/6\\
(\hat{S}\,\psi-3\,\hat{{G}})/6
\end{bmatrix}.
\end{equation}

%% \printbibliography
\bibliography{cas-refs.bib}

\end{document}

\endinput
%%
%% End of file `elsarticle-template-harv.tex'.